\documentclass[stsy,nonblindrev]{informs-stsy}

\OneAndAHalfSpacedXI
\usepackage{natbib}
 \bibpunct[, ]{(}{)}{,}{a}{}{,}%
 \def\bibfont{\small}%
\TheoremsNumberedThrough     % Preferred (Theorem 1, Lemma 1, Theorem 2)
\ECRepeatTheorems

\EquationsNumberedThrough    % Default: (1), (2), ...
\MANUSCRIPTNO{}

\usepackage[letterpaper]{geometry}
\usepackage{authblk}
\usepackage{amsmath,amsfonts}
\usepackage{color,graphicx}
\usepackage[colorlinks,linkcolor=blue,citecolor=cyan,anchorcolor=blue]{hyperref}
\usepackage{enumerate}
\usepackage{wrapfig}
\usepackage{accents}
\usepackage{tikz}
\usetikzlibrary{arrows,
                petri,
                topaths,decorations.pathreplacing}
\usepackage{multirow}
\usepackage{verbatim}
\usetikzlibrary{shapes,calc,arrows.meta,positioning}
\usetikzlibrary{plotmarks}
\usetikzlibrary{shapes}

\newcommand{\R}{\mathbb{R}}

\newcommand{\Z}{\mathbb{Z}}
\newcommand{\N}{\mathbb{N}}
\newcommand{\E}{\mathbb{E}}

\newcommand{\Prob}{\mathbb{P}}

\newcommand{\blue}[1]{#1}

\newcommand{\startproof}{\proof}
\newcommand{\finishproof}{\hfill $\square$\endproof}
\providecommand{\abs}[1]{\left\lvert#1\right\rvert}
\providecommand{\norm}[1]{\lVert#1\rVert}
\usepackage[normalem]{ulem}
\begin{document}
%%%%%%%%%%%%%%%%

% Outcomment only when entries are known. Otherwise leave as is and
%   default values will be used.
%\setcounter{page}{1}
%\VOLUME{00}%
%\NO{0}%
%\MONTH{Xxxxx}% (month or a similar seasonal id)
%\YEAR{0000}% e.g., 2005
%\FIRSTPAGE{000}%
%\LASTPAGE{000}%
%\SHORTYEAR{00}% shortened year (two-digit)
%\ISSUE{0000} %
%\LONGFIRSTPAGE{0001} %
%\DOI{10.1287/xxxx.0000.0000}%

% Author's names for the running heads
% Sample depending on the number of authors;
% \RUNAUTHOR{Jones}
% \RUNAUTHOR{Jones and Wilson}
% \RUNAUTHOR{Jones, Miller, and Wilson}
% \RUNAUTHOR{Jones et al.} % for four or more authors
% Enter authors following the given pattern:
\RUNAUTHOR{Braverman}

% Title or shortened title suitable for running heads. Sample:
% \RUNTITLE{Bundling Information Goods of Decreasing Value}
% Enter the (shortened) title:
\RUNTITLE{Convergence rates for the join-the-shortest queue system}

% Full title. Sample:
% \TITLE{Bundling Information Goods of Decreasing Value}
% Enter the full title:
\TITLE{The join-the-shortest-queue system in the Halfin-Whitt regime: rates of convergence to the diffusion limit}

% Block of authors and their affiliations starts here:
% NOTE: Authors with same affiliation, if the order of authors allows,
%   should be entered in ONE field, separated by a comma.
%   \EMAIL field can be repeated if more than one author
\ARTICLEAUTHORS{%
\AUTHOR{Anton Braverman}
\AFF{Kellogg School of Management, Northwestern University, Evanston, IL 60208, \EMAIL{anton.braverman@kellogg.northwestern.edu}} %, \URL{}}
% Enter all authors
} % end of the block

\ABSTRACT{We show that the steady-state distribution of the join-the-shortest-queue (JSQ) system converges, in the Halfin-Whitt regime,  to its diffusion limit at a rate of at least $1/\sqrt{n}$, where $n$ is the number of servers. Our proof uses Stein's method, and, specifically, the recently proposed \emph{prelimit} generator comparison approach.  \blue{The JSQ system is non-trivial, high-dimensional, and has a state-space collapse component, and our analysis may serve as a helpful example to readers wishing to apply the approach to their own setting. }}
%For example, the approach has the potential to be applied to other heavy-traffic regimes or load-balancing policies.   

% Sample
%\KEYWORDS{deterministic inventory theory; infinite linear programming duality;
%  existence of optimal policies; semi-Markov decision process; cyclic schedule}

% Fill in data. If unknown, outcomment the field

\maketitle
%%%%%%%%%%%%%%%%%%%%%%%%%%%%%%%%%%%%%%%%%%%%%%%%%%%%%%%%%%%%%%%%%%%%%%

% Samples of sectioning (and labeling) in MNSC
% NOTE: (1) \section and \subsection do NOT end with a period
%       (2) \subsubsection and lower need end punctuation
%       (3) capitalization is as shown (title style).
%
%\section{Introduction.}\label{intro} %%1.
%\subsection{Duality and the Classical EOQ Problem.}\label{class-EOQ} %% 1.1.
%\subsection{Outline.}\label{outline1} %% 1.2.
%\subsubsection{Cyclic Schedules for the General Deterministic SMDP.}
%  \label{cyclic-schedules} %% 1.2.1
%\section{Problem Description.}\label{problemdescription} %% 2.

% Text of your paper here

\section{Introduction}
\label{sec:introduction}
Consider a queueing system with $n$ identical servers, each with a finite buffer of length $b$. Customers arrive according to a Poisson process with rate $n\lambda$, and service times are i.i.d., exponentially distributed with rate $1$. Customers cannot change servers after the initial routing decision, and a customer arriving to a system where all servers are busy and all buffers are full is blocked. This is known as a parallel-server system.  A  load-balancing policy  specifies the manner in which arriving customers are assigned to the servers. In this paper, we consider the classical join-the-shortest-queue (JSQ) policy. Under JSQ, an arriving customer enters service immediately if at least one server is idle; if not, they get routed to the server with the smallest number of customers in its buffer. Ties are broken arbitrarily. We refer to this   as the JSQ system.

Parallel-server systems have generated immense interest in  recent years, and the JSQ policy is fundamental because it minimizes the expected customer delay  and  maximizes, with respect to stochastic order, the number of customers served in a given time interval; see, for instance, \cite{Wins1977,Webe1978}. For a  sample of recent work on the JSQ policy, we refer readers to  \cite{EryiSrik2012, MukhBorsLeeuWhit2016, GamaEsch2018,GuptWalt2019, BaneMukh2019, LiuYing2019,  BaneMukh2020, Brav2020, ZhouShro2020a, ZhouShro2020, ZhaoBaneMukh2021, HurtMagu2021, Caoetal2021}. Other popular load-balancing policies include the join-the-idle-queue policy (\cite{Stol2015a, MukhBorsLeeuWhit2016}), the idle-one-first policy (\cite{GuptWalt2019}),  and of course the     power-of-$d$ policy (\cite{VvedDobrKarp1996, Mitz2001}), but in this paper we focus on the JSQ policy. We make no attempt to give a comprehensive review of the literature on parallel-server systems, instead referring the reader to  \cite{BoorBorsLeeuMukh2021} for a recent survey.

Understanding the exact performance of the system is known to be  difficult  \blue{and much attention has been devoted over the past decade to heavy-traffic asymptotics}. The term ``heavy traffic'' refers to parameter regimes where the system utilization tends to one. ``Conventional heavy traffic'' assumes that the number of servers $n$ is fixed and $\lambda \uparrow 1$, while ``many-server heavy traffic'' assumes that  $n \to \infty $ and $\lambda \uparrow 1$ jointly. For two   examples of work  in the conventional heavy-traffic setting, see  \cite{EryiSrik2012} and \cite{ZhouShro2020}. In this paper, we use the term ``heavy-traffic'' to refer to the many-server setting -- the setting considered in most of the papers mentioned in the previous paragraph.   

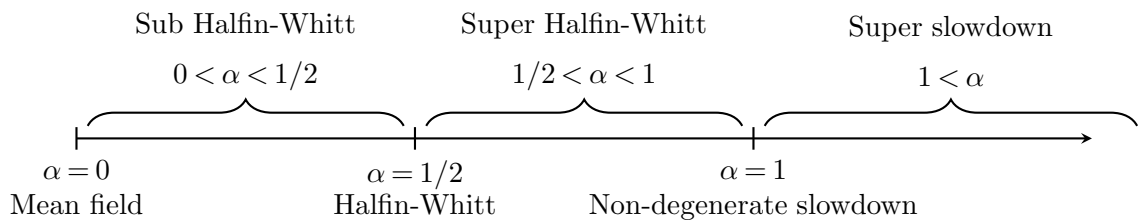
\begin{figure}[h!]
\centering
\begin{tikzpicture}[scale=1.5,>=stealth,
every node/.style={align=center,scale=1}] 
\draw[thick, ->] (0,0)--(9,0);
\draw[thick] (0,.1)--(0,-.1) node[below]{$\alpha = 0$};
\draw[thick] (3,.1)--(3,-.1) node[below]{$\alpha = 1/2$};
\draw[thick] (6,.1)--(6,-.1) node[below]{$\alpha = 1$}; 
\draw [thick, decorate,decoration={brace,amplitude=10pt ,raise=1ex}]
  (.1,0) -- (2.9,0) node[midway,yshift=3em]{Sub Halfin-Whitt \\ $0 < \alpha < 1/2$};
\draw [thick, decorate,decoration={brace,amplitude=10pt ,raise=1ex}]
  (3.1,0) -- (5.9,0) node[midway,yshift=3em]{Super Halfin-Whitt \\ $1/2 < \alpha < 1$};
\draw [thick, decorate,decoration={brace,amplitude=10pt ,raise=1ex}]
  (6.1,0) -- (9.4,0) node[midway,yshift=3em]{Super slowdown \\ $1 < \alpha$};
\draw (0,-.4)  node[below]{Mean field}; 
\draw (3,-.4)  node[below]{Halfin-Whitt}; 
\draw (6,-.4)  node[below]{Non-degenerate  slowdown};   
\end{tikzpicture}
\caption{The various many-server heavy-traffic regimes. Higher values of $\alpha$ represent heavier   loads. Existing work across the different parameter regimes is reviewed in Section~\ref{sec:litrev}. }
\label{fig:regimes}
\end{figure}

\blue{There are multiple many-server heavy-traffic regimes, depending on how $n$ and $\lambda$ jointly converge to their limit. For example, assuming that $\lambda = 1 - 1/n$ yields entirely different asymptotic behavior compared to when  $\lambda = 1 - 1/\sqrt{n}$. To capture all the possible heavy-traffic regimes, it is common practice to assume that} the per-server load $\lambda$ is related to the number of servers $n$ through $\lambda = 1 - \beta/n^{\alpha} \in (0,1)$ for some  $\alpha \geq 0$  and $\beta > 0$. In this paper we focus on the case when $\alpha = 1/2$; i.e.,  $\lambda = 1-\beta/\sqrt{n}$. This regime is known as the Halfin-Whitt regime and is ubiquitous across the queueing theory literature. It derives from the work of  \cite{HalfWhit1981}    and is  also known as the quality-and-efficiency-driven regime because it achieves reasonable customer wait times while maintaining high utilization of servers.  \blue{The full list of  parameter regimes} is found in   Figure~\ref{fig:regimes}. 

We now state and discuss our main results. Let $Q_i(t)$ be the number of servers with $i$ or more customers at time $t \geq 0$, noting that $Q_{i}(t) = 0$ for $i > b+1$. The process $\{Q(t) = (Q_1(t), \ldots, Q_{b+1}(t)) \}$ is an irreducible continuous-time Markov chain (CTMC) on a finite state space  and therefore possesses a unique stationary distribution. We let $Q = (Q_1,\ldots, Q_{b+1})$ be the random vector having the   stationary distribution of the CTMC. To describe the asymptotic behavior of $Q$, we let $\delta = 1/\sqrt{n}$ and define the diffusion-scaled random vector $X = (X_1, \ldots, X_{b+1})$ by  $X_1 = \delta(n-Q_1)$, and $X_i = \delta Q_i$ for $2 \leq i \leq b+1$. \blue{The results of \cite{GamaEsch2018} and \cite{Brav2020} imply that $X$  converges in distribution to some limiting $\R^{b+1}_+$-valued random vector $Y$ as $n \to \infty$. In this paper we establish an upper bound of order $1/\sqrt{n}$ on the rate of convergence to $Y$.}

The random variable $Y$ is distributed according to  the stationary distribution of  the diffusion process  $\{Y(t) \in \R^{b+1}_+\}$, which satisfies 
\begin{align}
&Y_1(t) = Y_1(0) + \sqrt{2} W(t) + \beta t - \int_{0}^{t} (Y_1(s) + Y_2(s)) ds + U(t), \notag \\
&Y_2(t) = Y_2(0) + U(t)  - \int_{0}^{t} Y_2(s) ds, \quad Y_3(t) = \cdots = Y_{b+1}(t) = 0, \label{eq:diffusion2}
\end{align} 
where $\{W(t)\}$ is standard Brownian motion and $\{U(t)\}$ is the unique nondecreasing, nonnegative process in the space of c\`{a}dl\`{a}g functions $D[0,\infty)$ satisfying $\int_{0}^{\infty} 1(Y_1(t)>0) d U(t) = 0$.  
The  diffusion $\{Y(t)\}$ was shown to be positive recurrent; see \cite{BaneMukh2019} or \cite{Brav2020}. Furthermore, \eqref{eq:diffusion2} implies that $Y_3 = \ldots = Y_{b+1} = 0$. 

 Our main result is   that there exists a constant $C(b,\beta)$ such that for all $n \geq 1$, and any  function $h: \R^{b+1}_+ \to \R$  whose first-order  and second-order  partial derivatives are bounded in magnitude by one, 
\begin{align}
\abs{\E h(X) - \E h(Y) } \leq C(b,\beta)/\sqrt{n}. \label{eq:intromain} 
\end{align} 
The assumption that $b$ is finite is used frequently in the proof of \eqref{eq:intromain} and, specifically, in the proof of Proposition~\ref{prop:diffbounds}. We deem the finite buffer assumption to be acceptable because    it was shown by \cite{Brav2020} that even with infinite-sized buffers, $\E Q_3 \leq C(\beta)$ for all $n \geq 1$ in the Halfin-Whitt regime, implying that  $X_3 \Rightarrow 0$, or that the mass concentrates on those states with at most one customer waiting. Moreover, \cite{LiuYing2020} showed that assuming finite buffers,  $\E Q_3 \to 0$ as $n \to \infty$ in the even busier super-Halfin-Whitt regime ($1/2 < \alpha < 1$).  
%I believe that the dependence of $C(b,\beta)$ on $b$ is merely a result of the proof technique, and conjecture that the right-hand side in \eqref{eq:intromain} does not depend on $b$, and that a similar convergence rate holds in the infinite-sized buffer setting. 

In addition to the novelty of our result, this paper  makes a methodological contribution.
 We prove \eqref{eq:intromain} using Stein's method, a framework introduced by \cite{Stei1972} that allows one to study the rate of convergence of a sequence of random variables to its limit. Popularized in the area  of queueing systems by \cite{Gurv2014, Ying2017, BravDai2017, Gast2017}, the generator comparison approach of Stein's method, attributed to \cite{Barb1988, Barb1990} and \cite{Gotz1991},  is used to study  convergence rates of steady-state Markov chain distributions to their diffusion, fluid, or mean-field  limits. For a few recent applications of the generator comparison approach in queueing, we refer the reader to \cite{ GaunWalt2020, HurtMagu2021, Lu2021, LiuGongYing2022}; this list is by no means comprehensive.   In this paper, we restrict our attention to the case when the limit is the stationary distribution of a diffusion process, referring the reader to \cite{Ying2017} for a treatment of fluid and mean-field limits.
 
  The generator approach requires bounds on various moments of the prelimit, known as moment bounds, and  bounds on  the derivatives of the solution to the Poisson equation for the limiting distribution. The  latter are  called gradient bounds in \cite{BravDai2017}, but in this paper we stick with the original term ``Stein factors'', or ``Stein factor bounds''; e.g.,  \cite{Ross2011}. While moment bounds can be difficult to obtain in some applications, Stein factor bounds are typically the bigger problem.   When the limit is one-dimensional, Stein factors are bounded using the explicit form of the solution to  the Poisson equation --- an ordinary differential equation.   When the limit is multidimensional,  the Poisson equation is a partial differential equation (PDE) that generally does not have an explicit solution, making Stein factor bounds harder to establish. Techniques proposed to obtain multidimensional Stein factor bounds include using a priori Schauder estimates from elliptic PDE theory as in \cite{Gurv2014}, using couplings to analyze and bound the sensitivity of the diffusion to its initial condition as in \cite{Barb1988} and \cite{GorhMack2016}, and bounding the Stein factors using Malliavin calculus as in  \cite{FangShaoXu2018} and \cite{Jinetal2021}. A detailed description of these techniques can be found in Section~1.1 of \cite{Brav2022}. However, despite  progress on multidimensional Stein factor bounds,  the JSQ system is not covered by existing results   because our limiting diffusion in \eqref{eq:diffusion2} is constrained to the nonnegative orthant via reflecting boundary conditions. 
   
To deal with the Stein factor bound problem, this paper promotes the use of the \emph{prelimit} generator comparison approach, which was recently proposed by \cite{Brav2022}  as an alternative to the generator comparison approach. The prelimit approach is the mirror image of the classical generator approach. Whereas the latter requires moment bounds on the prelimit $X$  and Stein factor bounds for limit $Y$, the former needs moment bounds on  $Y$ and Stein factor bounds for the prelimit $X$. For the moment bounds used in this paper, the result that all moments of $Y$ are finite,  proved by \cite{BaneMukh2019}, is sufficient because our limit  $Y$ does not depend on $n$. The  Stein factor bounds pose a bigger challenge, and we deal with them in Section~\ref{sec:perturbation}. It  was noted in \cite{Brav2022} that the prelimit and classical generator comparison approaches should be equivalent, in theory, in the sense that any bound on $\abs{\E h(X) - \E h(Y) }$ obtained using one of them should  be attainable using the other. However, in practice,  one approach could be more tractable, or convenient, to work with; see, for instance, the example in Section~4 of \cite{Brav2022}. \blue{In the case of the JSQ system, we discuss in Remark~\ref{rem:prelimgen} of Section~\ref{sec:proof_relate} how the discrete state space  simplifies the analysis of the couplings we use to establish Stein factor bounds, because the initial spacing of the coupled systems is preserved until coupling.}

The introduction of the prelimit approach in \cite{Brav2022} was intended to be gentle, with the only example used there being the $M/M/1$ system.  Our application of the approach to the  JSQ system   exposes   all of its moving pieces and can be useful to those who want to apply the prelimit approach to their own setting. \blue{For example, some of the technical components of this paper that could be useful in other settings include: the regenerative argument used to establish first-order Stein factor bounds in Section~\ref{sec:first_diff}, the approach we use to bound $\E|X|$ in Section~\ref{sec:moments}, and our treatment of reflecting boundary conditions in Appendix~\ref{sec:eps1}.}

%As our second contribution,  this paper is a step in the direction of obtaining convergence-rate results for parameter regimes beyond $\alpha = 1/2$  and even for non-JSQ load-balancing policies. \red{?? Though we do not obtain results when $\alpha \neq 1/2$, to generalize our appraoch we would need to establish Stein factor bounds. Seems like ZhaoBaneMukh could hold the answers for doing this when $1/2 < \alpha < 1$.  generalizing our results to such regimes  , } Our main technical hurdle, as discussed in Section~\ref{sec:perturbation}, is bounding the JSQ system's Stein factors, which correspond to the sensitivity of the JSQ system to a shift in its initial condition.  Obtaining similar bounds when $\alpha \neq 1/2$  or in non-JSQ policy settings would be the single biggest challenge to overcome for a convergence-rate result  and is likely to be difficult. However, that task may be feasible  for  at least some   load-balancing policies and/or parameter regimes. Indeed, in the setting of $\alpha \in (1/2,1)$, the contemporary work of \cite{ZhaoBaneMukh2021}  contained  several hitting-time estimates similar to the ones used in our proof of Lemma~\ref{lem:coupling_time_bound}, which is the key to our first-order Stein factor bounds. 

It should be noted  that \cite{HurtMagu2021} and \cite{ZhouShro2020a}  used the classical generator comparison approach to obtain rates of convergence of the steady-state total customer count to an  exponential random variable for  $\alpha > 2$. The former paper was in the continuous-time setting, while the latter considered the discrete-time system, and the results in both papers also hold for routing policies other than JSQ, such as the power-of-$d$ policy. Since the limiting random variable in both papers is one-dimensional, the Stein factors bounds do not pose a challenge there.

\subsection{Literature Review}
\label{sec:litrev}
Let us first review the literature on the analysis of the JSQ system in the various many-server heavy-traffic regimes. Most of the work has been done in the setting with infinite buffer sizes, so, unless otherwise noted,  we assume that $b = \infty$. In \cite{GamaEsch2018},   the authors established  the process-level convergence of $\{X(t)\}_{n=1}^{\infty}$ to its diffusion limit in the Halfin-Whitt regime ($\alpha = 1/2$). That paper triggered a wave of interest in the many-server heavy-traffic asymptotics of the JSQ system.  Convergence of the stationary distributions was later established by \cite{Brav2020}, and the behavior of the stationary distribution of the limiting diffusion was studied by \cite{BaneMukh2019,BaneMukh2020}. Our work fits with this group of papers, elevating the steady-state convergence result to one with rates of convergence.

Outside the Halfin-Whitt regime, \cite{MukhBorsLeeuWhit2016} studied the transient and steady-state behavior of the JSQ system's fluid limit when  $\lambda = 1-\beta < 1$ is a fixed constant ($\alpha = 0$), and \cite{GuptWalt2019} established process-level convergence to the diffusion limit when $\alpha = 1$; known as the non-degenerate slowdown (NDS) regime and introduced by \cite{Atar2012}. In the sub-Halfin-Whitt regime when $\alpha \in (0,1/2)$,  \cite{LiuYing2019} assumed finite buffers and  obtained bounds on the steady-state total customer count in the system. A similar result was obtained for Coxian-2 service times by \cite{LiuGongYing2022}, and by \cite{LiuYing2020} for the  super-Halfin-Whitt regime $\alpha \in (1/2,1)$. Another recent work in the super-Halfin-Whitt regime was by \cite{ZhaoBaneMukh2021}, who worked with infinite buffers and established transient and steady-state diffusion limits for the normalized total queue length process. Their analysis exploited the regenerative structure of the JSQ system  and contained several hitting-time estimates very close to our own estimates needed for the Stein factor bounds in Section~\ref{sec:perturbation}. Lastly,   both \cite{HurtMagu2021} and \cite{ZhouShro2020a}  established rates of convergence to the exponential distribution for the steady-state normalized total customer count. Their results covered the case when $\alpha > 2$. 

Other works have used Stein's method in the setting of parallel-server systems beyond  \cite{HurtMagu2021} and \cite{ZhouShro2020a}. In \cite{LiuYing2019, LiuYing2020,  LiuGongYing2022}, the authors used Stein's method for mean-field analysis to obtain bounds on steady-state performance metrics of interest, like $\E Q_2$ for instance, for the power-of-$d$ system. Another line of work on power-of-$d$ systems was by \cite{Gast2017, GastVanh2017, GastBortTrib2019}, where the authors showed how to derive refined mean-field models for improved steady-state approximations. More recently, \cite{HairLiuYing2021} provide calculable error bounds for the mean-field approximation of the power-of-two-choices model.

%To my knowledge, there are three  approaches to date for bounding the Stein factors of a diffusion. One approach proposed by \cite{Gurv2014} uses apriori Schauder estimates from elliptic PDE theory to bound the diffusion Stein factors by a quantity involving the solution of the Poisson equation itself. The solution can be bounded by finding a Lyapunov function for the diffusion  satisfying an exponential ergodicity condition.  The second approach uses multiple coupled copies of the diffusion to analyze its sensitivity to its initial condition.  For a discussion on Stein factor bounds for fluid, or mean-field limits, we refer the reader to \cite{Ying2016}.

\subsection{Notation}
We use $\Z$ to denote the set of integers  and let $\N = \{0,1,2,\ldots\}$. For any   $k \in \N $ and  $B \subset \R^{d}$, we let $C^{k}(B)$ be the set of all $k$-times continuously differentiable functions $f: B \to \R$.     We let $e \in \R^{d}$ be the vector whose elements all equal  $1$  and let $e^{(i)}$ be the element with $1$ in the $i$th entry  and zeros otherwise. For any $\delta > 0$ and integer $d > 0$, we let $\delta \Z^{d} = \{\delta k :\ k \in \Z^{d}\}$ and define $\delta \N^{d}$ similarly. For any function $f: \delta \Z^{d} \to \R$, we define the forward difference operator in the $i$th direction as
\begin{align*}
\Delta_{i} f(\delta k) = f \big( \delta (k+e^{(i)}) \big) - f(\delta k), \quad k \in \Z^{d}, \  1 \leq i \leq d,
\end{align*}
and for $j \geq 0$, we define 
\begin{align}
\Delta_i^{j+1} f(\delta k) = \Delta_{i}^{j} f(\delta(k+e^{(i)})) - \Delta_{i}^{j} f(\delta k), \label{eq:diffdef}
\end{align}
with the convention that $\Delta_i^{0} f(\delta k) = f(\delta k)$. For a vector $a \in \N^{d}$, we also let 
\begin{align*}
\Delta^{a}  f(\delta k) =&\ \Delta_{1}^{a_1} \ldots \Delta_{d}^{a_d} f(\delta k),
\end{align*} 
and if  $f: \R^{d} \to \R$, then
\begin{align*}
\frac{\partial^{a}}{\partial x^{a}} f(x) =&\ \frac{\partial^{a_1}}{\partial x_1^{a_1}} \ldots \frac{\partial^{a_d}}{\partial x_d^{a_d}} f(x),
\end{align*}
and we adopt the convention that $\frac{\partial^{0}}{\partial x^{0}} f(x) = f(x)$. For any $x \in \R^{d}$, we define $\norm{x}_{1} = \sum_{i=1}^{d} \abs{x_i}$  and use $\abs{x}$ to denote the Euclidean norm. \blue{For any $f: \R^{d} \to \R$, we let $\norm{f}_{\infty} =  \sup_{x \in \R^{d}} \abs{f(x)}$}. Throughout the paper, we will often use $C$ to denote a generic positive constant  that may change from line to line  and that is independent of any parameters not explicitly specified.

\section{Main Result}
\label{sec:main}
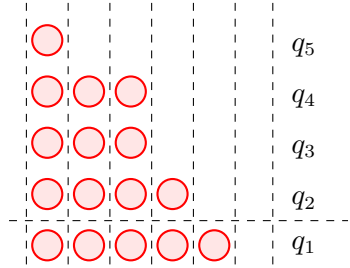
\begin{figure}[h!]
\centering
\begin{tikzpicture}
 [server/.style={circle, inner sep=1.0mm, minimum width=.8cm,
    draw=green,fill=green!10,thick},
    customer/.style={circle, inner sep=0.5mm, minimum width=.4cm,
    draw=red,fill=red!10,thick},
    blue/.style={circle, inner sep=0.5mm, minimum width=.4cm,
    draw=blue,fill=blue!10,thick},
    green/.style={circle, inner sep=0.5mm, minimum width=.4cm,
    draw=green,fill=green!10,thick},
  buffer/.style={rectangle, rounded corners=3pt,
    inner sep=0.0mm, minimum width=.9cm,  minimum height=.8cm,
    draw=orange,fill=blue!10,thick},
 vbuffer/.style={rectangle,rounded corners=3pt,
     inner sep=0.0mm,  minimum width=.6cm, minimum height=.8cm,
     draw=orange,fill=blue!10,thick}
]
	\node (ref){};
  \node[customer]  (C1) {};
  \node[customer]  (C2) [right =.05in of C1]   {};
  \node[customer]  (C3) [right =.05in of C2]   {};
  \node[customer]  (C4) [right =.05in of C3]   {};
  \node[customer]  (C5) [right =.05in of C4]   {};
  \node  (C6) [right =.05in of C5]   {};
  
  \node[customer]  (C11) [above =.1in of C1]   {};
  \node[customer]  (C21) [above =.1in of C2]   {};
  \node[customer]  (C31) [above =.1in of C3]   {};
  \node[customer]  (C41) [above =.1in of C4]   {};
  \node  (C61) [above =.1in of C6]   {};
  
  \node[customer]  (D1) [above =.1in of C11]   {};
  \node[customer]  (D2) [above =.1in of C21]   {};
  \node[customer]  (D3) [above =.1in of C31]   {};
  \node (D6) [above =.1in of C61]   {};
  
  \node[customer]  (E1) [above =.1in of D1]   {};
  \node[customer]  (E2) [above =.1in of D2]   {};
  \node[customer]  (E3) [above =.1in of D3]   {};
  \node (E6) [above =.1in of D6]   {};
  
  \node[customer]  (F1) [above =.1in of E1]   {};
  \node (F6) [above =.1in of E6]   {};

 \draw[dashed]  ($(C1.north)+(-0.2in,.045in)$) edge ($(C5.north)+(0.7in,.045in)$);
 \draw[dashed]  ($(C1.west)+(-0.025in,-.1in)$) edge ($(C1.west)+(-0.025in,1.3in)$);
 \draw[dashed]  ($(C2.west)+(-0.025in,-.1in)$) edge ($(C2.west)+(-0.025in,1.3in)$);
 \draw[dashed]  ($(C3.west)+(-0.025in,-.1in)$) edge ($(C3.west)+(-0.025in,1.3in)$);
 \draw[dashed]  ($(C4.west)+(-0.025in,-.1in)$) edge ($(C4.west)+(-0.025in,1.3in)$);
 \draw[dashed]  ($(C5.west)+(-0.025in,-.1in)$) edge ($(C5.west)+(-0.025in,1.3in)$);
 \draw[dashed]  ($(C6.west)+(-0.025in,-.1in)$) edge ($(C6.west)+(-0.025in,1.3in)$);
 \draw[dashed]  ($(C6.east)+(0.065in,-.1in)$) edge ($(C6.east)+(0.065in,1.3in)$);
  
  \node[left, align=left] (X1) at ($(C6)+(.4in,0)$) {$q_1$};
  \node[left, align=left] (X2) at ($(C61)+(.4in,.025in)$) {$q_2$}; 
  \node[left, align=left] (X3) at ($(D6)+(.4in,.09in)$) {$q_3$}; 
  \node[left, align=left] (X4) at ($(E6)+(.4in,.15in)$) {$q_4$}; 
  \node[left, align=left] (X5) at ($(F6)+(.4in,.21in)$) {$q_5$}; 
 
 \end{tikzpicture}
\caption{An example of a state   $Q(t) = q$ in a system where the number of servers $n = 5$.  Customers below the dashed horizontal line are in service, while those above are waiting in buffers. Each vertical column corresponds to a server and its buffer.  }
\label{fig:stateexample}
\end{figure}

Recall  that $Q_i(t)$ is the number of servers with $i$ or more customers at time $t \geq 0$ and that   $\{Q(t) = (Q_{i}(t) )_{i=1}^{b+1}\}_{t \geq 0}$ is an irreducible CTMC with state space given by
\begin{align}
S_{Q} = \big\{  q \in \{0,\ldots,n\}^{b+1} :  q_{i} \geq q_{i+1}  \big\}. \label{eq:sq}
\end{align} 
Figure~\ref{fig:stateexample} gives an example of a  state  $q \in S_Q$.
 We assume that $\lambda = 1 - \beta/\sqrt{n}$ for some fixed $\beta > 0$. 
 Let $\delta = 1/\sqrt{n}$ and  define the diffusion-scaled CTMC $\{X(t)\}$ by
\begin{align*}
&X_1(t) = \delta(n - Q_1(t)), \quad X_i(t) = \delta Q_i(t), \quad 2 \leq i \leq b+1,
\end{align*}
which takes values on the state space 
\begin{align*}
S = \big\{ (x^q_1,x^q_2,\ldots, x^q_{b+1}) =  \big(\delta(n-q_1),\delta q_2, \ldots, \delta q_{b+1} \big) :  q \in S_{Q}  \big\}.
\end{align*}
We will often use $x^q \in S$ and $q \in S_{Q}$ interchangeably. Recalling that $\Delta_{i} f(x^q) = f(x^q + \delta e^{(i)}) - f(x^q)$, for any $f:S \to \R$, the infinitesimal generator of $\{X(t)\}$ satisfies  
\begin{align}
 G_{X} f(x^q)  =&\ -1(q_1 < n) n \lambda  \Delta_1 f(x^q- \delta e^{(1)})  + n \lambda \sum_{j=1}^{b} 1(q_1=\ldots = q_j = n, q_{j+1} < n)  \Delta_{j+1} f(x^q) \notag \\
&+ (q_1 - q_{2})  \Delta_1 f(x^q) - \sum_{j=2}^{b} (q_j - q_{j+1})   \Delta_j f(x^q-\delta e^{(j)}) - q_{b+1}   \Delta_{b+1} f(x^q-\delta e^{(b+1)}). \label{eq:gxdef}
\end{align}
The first line of transitions in \eqref{eq:gxdef} correspond to arrivals. We see that for $j \geq 2$,  the $j$th component of $x^q_j$ only grows provided the preceding $j-1$ horizontal levels, as depicted in Figure~\ref{fig:stateexample}, are full. The transitions in the second line of \eqref{eq:gxdef} correspond to service completions. Using Figure~\ref{fig:stateexample} again, we interpret $(q_{j} - q_{j+1})$ as the number of servers (vertical columns) with exactly $j$ customers.

Recall that $X = (X_1, \ldots, X_{b+1})$ and $Y = (Y_1,Y_2, 0, \ldots, 0)$ are distributed according to the stationary distributions of the  scaled CTMC  and  the  diffusion   $\{Y(t) \in \R^{b+1}_+\}$ defined in \eqref{eq:diffusion2}, respectively. \blue{Going forward, we note that unless explicitly stated, all expectations are with respect to the stationary distribution at hand; i.e., either $X$ or $Y$. } To state our main result,  we define
\begin{align*}
&\mathcal{M}_{j} =  \Big\{h^{\ast}: \R^{b+1}  \to \R,\  \Big{\|} \frac{\partial^{a}}{\partial x^a} h^{\ast}(x)  \Big{\|}_{\infty} \leq 1,\ 1 \leq \norm{a}_{1} \leq j  \Big\}, 
\end{align*}
and $d_{\mathcal{M}_{j}}(X,Y) =  \sup_{h^{\ast} \in \mathcal{M}_{j}} \big| \E h^{\ast}(X) - \E h^{\ast}(Y) \big|$. We use an asterisk to emphasize that $h^{\ast}(x)$ is defined on the continuum $\R^{b+1}$. Later we will drop the asterisk to refer to functions defined only on the grid $\delta \Z^{b+1}$.  It was shown in Lemma 2.2 of \cite{GorhMack2016} that $\mathcal{M}_{3}$ is a convergence-determining class; i.e., $d_{\mathcal{M}_{3}}(U,V) \to 0$ implies $U$ and $V$ converge in distribution.  The following is our main result.
%The following are our main results, which bound the rate at which $d_{\mathcal{M}_{2}}(X ,Y )$ converges to zero in the cases when $b < \infty$, and $b = \infty$. 
\begin{theorem}
\label{thm:main}
For any $0 < b < \infty$, there exists a constant $C(b,\beta)$ such that for all $n \geq 1$, 
\begin{align}
d_{\mathcal{M}_{2}}(X ,Y ) = \sup_{h^{\ast} \in \mathcal{M}_{2}} \big| \E h^{\ast}(X) - \E h^{\ast}(Y) \big|  \leq C(b,\beta)/\sqrt{n}. \label{eq:main}
\end{align}
\end{theorem}
Note that  $\mathcal{M}_{2}$ is also a convergence-determining class because $\mathcal{M}_{3} \subset \mathcal{M}_{2}$.
We prove Theorem~\ref{thm:main} in Section~\ref{sec:main_proof} using the prelimit generator approach of Stein's method. Multiple parts of the proof assume that $n$ is large enough, say, $n > N(\beta)$ for some $N(\beta) > 0$. We can make this assumption without loss of generality by redefining $C(b,\beta)$ to be larger than $\max_{1 \leq n \leq N(\beta)} d_{\mathcal{M}_{2}}(X ,Y )$. 
%The following is our second theorem, proved in Section~\ref{sec:t2pf}.
%\begin{theorem}
%\label{thm:2}
%Suppose that  $b = \infty$. For any  $\epsilon > 0$, There exists a constant $C(\beta, \epsilon)$ such that for all $n \geq 1$, 
%\begin{align}
%d_{\mathcal{M}_{2}}(X ,Y ) = \sup_{h^{\ast} \in \mathcal{M}_{2}} \big| \E h^{\ast}(X) - \E h^{\ast}(Y) \big|  \leq C(b,\epsilon)/n^{1/2-\epsilon}. \label{eq:main2}
%\end{align}
%\end{theorem}
 
%\begin{remark}
%The constant $C(b,\beta)$ in Theorem~\ref{thm:main} depends on the buffer size $b$, while most of the prior results about the JSQ model hold  for infinite-sized buffers. 
%\end{remark}

\subsection{Proving Theorem~\ref{thm:main}}
\label{sec:main_proof}
\blue{Central to our proof is the ability to extend any grid-valued function  to be defined on all of $\R^{b+1}_+$. Although there are infinitely many such extensions, we use a polynomial spline $A$ that extends grid-valued functions $f: \delta \N^{b+1} \to \R$ to functions $A f: \R^{b+1}_+ \to \R$. We leave the detailed construction to  Appendix~\ref{sec:prop1proof} because for this section, it suffices to know that $A$ is a linear operator, that $Af \in C^{2}(\R^{b+1}_+)$,  and that $A$ applied to a constant equals that constant.  
Recalling that $\delta = 1/\sqrt{n}$, the following auxiliary lemma is needed.  }
\begin{lemma}
\label{lem:convdet}
Define  
\begin{align*}
\mathcal{M}_{disc,j}(c) =  \Big\{h: \delta\N^{b+1}  \to \R,\  \abs{\Delta^{a} h(\delta k)}  \leq c \delta^{\norm{a}_{1}},\ 1 \leq \norm{a}_{1} \leq j,\  \delta k \in \delta \N^{b+1} \Big\}.
\end{align*}
There exist some   $C,C'> 0$ independent of any JSQ model parameters such that 
\begin{align}
d_{\mathcal{M}_{2}}(X ,Y ) \leq \sup_{\substack{h \in \mathcal{M}_{disc,2}(C)  }} \abs{\E h(X) - \E Ah(Y)} + C' \delta. \label{eq:mdisc}
\end{align}
\end{lemma}
\startproof{Proof of Lemma~\ref{lem:convdet}}
The result follows by repeating the arguments used in the proof of Lemma~1 in \cite{Brav2022}. 
\finishproof
 Going forward, when we write $\mathcal{M}_{disc,2}(C)$, the constant $C$ is assumed to be the one in Lemma~\ref{lem:convdet}. \blue{Furthermore, note that if $h(0) \neq 0$, then the linearity of $A$ and the fact that $A$ applied to a constant equals that constant implies that $\tilde h(x) = h(x) - h(0)$ satisfies $\E \tilde h(X) - \E A \tilde h(Y)  = \E h(X) - \E A h(Y)$. We therefore, without loss of generality, consider only those  $h \in \mathcal{M}_{disc,2}(C)$ such that $h(0) = 0$.}
 
  \blue{To prove Theorem~\ref{thm:main}, we bound the right-hand side of \eqref{eq:mdisc} with the help of the following two ingredients. The first ingredient is a rate-conservation law for  $\{Y(t)\}$, proved in Appendix~\ref{sec:supportmain}.}
\begin{lemma} \label{lem:ito}
Given $f \in C^{2}(\R^{b+1}_+)$, define 
\begin{align}
G_Y f(x) = \big(\beta - (x_1 + x_2)\big) \frac{\partial}{\partial x_1} f(x) - x_2 \frac{\partial}{\partial x_2} f (x) + \frac{\partial^2}{\partial x_1^2} f (x), \quad x  \in \R^{b+1}_+. \label{eq:gy}
\end{align}
If $\E \abs{f(Y)} < \infty$ and $\E \abs{G_Y f(Y)} < \infty$, and if $Y(0)$ is initialized according to $Y$, then
\begin{align}
  \E G_Y f(Y)  + \E \Big( \int_{0}^{1} \Big(\frac{\partial}{\partial x_1} f(Y(s)) + \frac{\partial}{\partial x_2} f(Y(s)) \Big)1(Y_1(s) = 0)  dU(s) \Big) = 0. \label{eq:ito}
\end{align}
\end{lemma}
The second ingredient is the Poisson equation. For  $h: \delta \N^{b+1} \to \R$ and $c \in \R$, let
\begin{align*} 
f_h^{(c)}(x^q) =&\ c + \int_{0}^{\infty} \big(\E_{x^q} h(X(t))- \E h(X) \big) dt, \quad x^q \in  S,
\end{align*}
which  is well defined because  the CTMC has a finite state space and is therefore exponentially ergodic. Furthermore, Lemma~2 of \cite{Brav2022} (see also Lemma~1 of \cite{Barb1988}) implies that 
\begin{align}
G_X f_h^{(c)}(x^q) = \E h(X) - h(x^q), \quad x^q \in S. \label{eq:discrete_poisson}
\end{align}
\blue{Most applications of Stein's method have $c = 0$, but we choose  $c = c^{\ast} = - f_h^{(0)}(0)$  and define
\begin{align} \label{eq:relval}
f_h(x^q) = f_h^{(c^{\ast})}(x^q) =&\ \int_{0}^{\infty} \big(\E_{x^q} h(X(t))- \E h(X) \big) dt - \int_{0}^{\infty} \big(\E_{0} h(X(t))- \E h(X) \big) dt \notag \\
=&\   \int_{0}^{\infty} \big(\E_{x^q} h(X(t))- \E_{0} h(X(t)) \big) dt, \quad x^q \in  S.
\end{align}
Our choice of $c$ yields $f_h(0) = 0$, which comes in handy later when we need to bound $\abs{f_h(x^q)}$ in Proposition~\ref{prop:diffbounds}.} Going forward, we assume that $c = c^{\ast}$ when referring to \eqref{eq:discrete_poisson}. 

\blue{Let us give an informal roadmap for bounding \eqref{eq:mdisc}, with the formal statement of the bounds left to Proposition~\ref{prop:main} below. We bound \eqref{eq:main} by  comparing the CTMC and diffusion generators. However, the former is defined only on a subset of $\R^{b+1}_+$, which requires the following workaround.  Suppose that we are given a set $B \subset \R^{b+1}_+$  such that (a)  $\E h(X) - A h(x) = A G_{X} f_h(x)$  for $x \in B$ and   (b) the probability that $Y \not \in B$ goes to zero rapidly (we will make this precise) as $n \to \infty$.   We decompose $\E h(X) - A h(x)$ as 
\begin{align*}
\E h(X) - A h(x) =&\ A G_{X} f_h(x) 1(x \in B) + \big(\E h(X) - A h(x)\big) 1(x \not \in B)   
\end{align*}
and take expected values with respect to $Y$ (we will show that these are finite) to get
\begin{align*}
\E h(X) - \E A h(Y) =&\ \E  \big(A G_{X} f_h(Y) 1(Y \in B)\big) + \E\Big(\big(\E h(X) - A h(Y)\big) 1(Y \not \in B)\Big). 
\end{align*}
Now extend $f_h(x^q)$ to  $\delta \N^{b+1}$ by defining  $f_h(x^q) = 0$ for $x^q \in \delta \N^{b+1} \setminus S$  and consider $Af_h(x)$. Provided that $\E \abs{Af_h(Y)} < \infty$  and $\E \abs{G_Y Af_h(Y)} < \infty$, we can invoke Lemma~\ref{lem:ito}  with $f(x) = A f_h(x)$ there to conclude that
\begin{align}
\E h(X) - \E Ah(Y) =&\  \E  \Big(\big(A G_{X} f_h(Y)-G_Y Af_h(Y)\big) 1(Y \in B)\Big)   \notag \\
&+  \E  \Big(\big(\E h(X) - A h(Y)-G_Y Af_h(Y)\big) 1(Y \not \in B)\Big) \notag \\
&-  \E  \Big( \int_{0}^{1}\Big(\frac{\partial}{\partial x_1} Af_h(Y(s)) + \frac{\partial}{\partial x_2} Af_h(Y(s)) \Big)1(Y_1(s) = 0)  dU(s) \Big), \label{eq:prefinalerror}
\end{align}
\blue{where $Y(0)$ in the third line is initialized according to $Y$.}  
We bound the first line  by showing that $G_X$ and $G_Y$ are close to one another. The middle term is small due to our choice of $B$ and the last term can be bounded because the JSQ system exhibits reflecting behavior similar to $\{Y(t)\}$ at the boundary $\{x \in S: x^q_1=0\}$. As a final remark, our choice of  $f_h(x^q) = 0$ for $x^q \in \delta \N^{b+1} \setminus S$ is made for convenience and is not essential to the proof, because the probability that $Y \not \in B$  shrinks  rapidly as $n \to \infty$.}

To state the following proposition,   define $k: \R^{b+1} \to \Z^{b+1}$ elementwise by  $k_{j}(x) = \lfloor x_j/\delta\rfloor$. For notational convenience, we also define $I = \big\{i = (i_1,i_2,0,\ldots,0) \in \N^{b+1} :\ 0 \leq i_1,i_2 \leq 4 \big\}$. The following proposition  is proved  in Appendix~\ref{sec:prop1proof}. 
\begin{proposition}
\label{prop:main}
 If $h \in \mathcal{M}_{disc,2}(C)$, then $Ah(Y)$, $Af_h(Y)$, and $G_Y Af_h(Y)$ are integrable, and  \eqref{eq:prefinalerror} holds. Furthermore, suppose that $n > 16$, define 
\begin{align*}
B =&\  \{ (x_1,x_2, 0, \ldots, 0) \in \R^{b+1}_+ :    x_2  + x_1 \leq \delta (n/2 - 8) = (n/2 - 8)/\sqrt{n} \},
\end{align*} 
and   let
\begin{align*} 
&\varepsilon_1(Y)  =   \big(A G_X f_h(Y) -  G_Y A f_h(Y)\big) 1(Y \in B), \\
&\varepsilon_2(Y)  = \big(\E h(X) - A h(Y)-G_Y Af_h(Y)\big) 1(Y \not \in B), \\
&\varepsilon_3(Y)  =   \Big(\frac{\partial}{\partial x_1} Af_h(Y) + \frac{\partial}{\partial x_2} Af_h(Y)\Big) 1(Y \in B), \text{ and }   \\
&\varepsilon_4(Y)  =   \Big(\frac{\partial}{\partial x_1} Af_h(Y) + \frac{\partial}{\partial x_2} Af_h(Y)\Big) 1(Y \not \in B).
\end{align*}
There exist   $C(\beta), C(b,\beta)> 0$   independent of $h(x)$ and $n$ such that 
\begin{align*}
\abs{\varepsilon_1(Y)} \leq&\ C(\beta)\big(1 + \delta^{-1}Y_2\big) \max_{\substack{i \in I \\ a_1+a_2 = 2}} \abs{\Delta_1^{a_1}\Delta_2^{a_2}  f_h\big(\delta(k(Y)   + i) \big) } + C(\beta) \delta^{-2} \max_{\substack{i \in I}} \abs{\Delta_{1}^{3}   f_h\big(\delta (k(Y)  +i )\big)} \\
&+   C(\beta) \delta^{-2} 1(Y_1 \leq \delta)\max_{\substack{i \in I \\ i_1 = 0}}  \abs{(\Delta_1^2 - (\Delta_1+\Delta_2))    f_h\big(\delta (k(Y)  +i )\big)}, \\ 
 \abs{\varepsilon_2(Y)} \leq&\  C(b,\beta) 1(Y \not \in B) \delta^{-2} (1 + Y_1+Y_2) \max_{\substack{ i \in I }} \abs{ f_h(\delta (k(Y)+i))}, \\ 
 \abs{ \varepsilon_3(Y) } \leq&\  C(\beta) \delta^{-1}  1(Y \in B)   \Big(\abs{(\Delta_1 + \Delta_2) f_h(\delta k(Y)}  +   \max_{\substack{i \in I \\  a_1+a_2 = 2 }} \abs{\Delta_{1}^{a_1}  \Delta_{2}^{a_2} f_h\big(\delta (k(Y)+i )\big)} \Big) , \\
  \abs{ \varepsilon_4(Y) } \leq&\ C(\beta) \delta^{-1}  1(Y \not \in B) \max_{\substack{ i \in I }} \abs{ f_h(\delta (k(Y)+i))}. 
\end{align*}  
\end{proposition}
Note that $\varepsilon_1(Y)$ and $\varepsilon_2(Y)$ are related to the first and second lines of \eqref{eq:prefinalerror}, respectively, while $\varepsilon_3(Y)$ and $\varepsilon_4(Y)$ are related to the last line there.  From the bounds in Proposition~\ref{prop:main}, we see  that the bound on \eqref{eq:prefinalerror} depends on the CTMC through the function   $f_h(x^q)$ and its differences, and on the diffusion through the distribution of $Y$. The differences of  $f_h(x^q)$ are commonly known as Stein factors, and the following proposition, proved in Section~\ref{sec:perturbation}, exhibits the Stein factor bounds we need to prove Theorem~\ref{thm:main}. 
\begin{proposition}
\label{prop:diffbounds}
There exists   $C(\beta,b) > 0$ such that for any $n \geq 1$ and  $h \in \mathcal{M}_{disc,2}(C)$, 
\begin{align*}
\abs{\Delta_1^{ a_1} \Delta_2^{ a_2} f_h(x^q)} \leq&\  C(\beta,b)\delta^{a_1+a_2}(1+x^q_2)^{a_1+a_2}, 
\end{align*} 
for all  $a_1,a_2 \geq 0$ with $1 \leq a_1+a_2 \leq 2$, and all $x^q \in S$ with $x^q_1 \leq \delta(n-a_1)$, $x^q_2 \leq \delta(n-a_2)$, and $x^q_3 = 0$.
Furthermore, 
\begin{align*}
&\abs{f_h(x^q)} \leq  C(\beta,b) (1 + x^q_2)(x^q_1 + x^q_2)/\delta , &\quad x^q \in S,\ x^q_3 = 0, \\
&\abs{\Delta_1^3 f_h(x^q)} \leq  C(\beta,b)\delta^{3}(1+x^q_2)^{3}, &\quad x^q \in S,\ x^q_1 \leq \delta(n-3),\ x^q_3=0,
\end{align*} 
and for all $x^q \in S$  with $ x^q_1 = 0$, $0 \leq x^q_2 \leq \delta(n-1)$, and $x^q_3 = 0$,
\begin{align*}
\abs{(\Delta_1 + \Delta_2) f_h(x^q)} \leq&\  C(\beta,b)\delta^2 (1+x^q_2)^2  \quad  \text{ and }   \\
 \abs{(\Delta_1^2 - (\Delta_1+\Delta_2))    f_h(x^q)} \leq&\ C(\beta,b) \delta^3(1+x^q_2)^{3}.
\end{align*}
\end{proposition}
The last component needed for the proof of Theorem~\ref{thm:main} is   the following lemma.
\begin{lemma}
\label{lem:ymoments}
All moments of $Y_1$ and $Y_2$ are finite. Furthermore, suppose that $Y(0)$ is initialized according to $Y$. Then for any $j > 0$, 
\begin{align}
\E Y_2^{j+1} =  \Big( \int_{0}^{1} (Y_2(s))^{j} 1(Y_1(s) = 0)  dU(s) \Big). \label{eq:ubound}
\end{align} 
\end{lemma}
\startproof{Proof of Lemma~\ref{lem:ymoments}}
The finiteness of the moments follows from
Theorem~2.1 of \cite{BaneMukh2019} and \eqref{eq:ubound} is implied by \eqref{eq:ito} of Lemma~\ref{lem:ito} with $f(y) = y_2^{j+1}$ there.
\finishproof  
\startproof{Proof of Theorem~\ref{thm:main}}
Initialize $Y(0)$ according to $Y$. 
Using \eqref{eq:prefinalerror} and the definitions of $\varepsilon_1(Y), \ldots, \varepsilon_4(Y)$, it follows that 
\begin{align*}
\E h(X) - \E Ah(Y) =&\ \E \varepsilon_1(Y) + \E \varepsilon_2(Y)  - \E \Big( \int_{0}^{1}\big(\varepsilon_3(Y(s)) + \varepsilon_4(Y(s)) \big)1(Y_1(s) = 0)  dU(s) \Big).
\end{align*}
We  argue that    $\abs{\E h(X) - \E Ah(Y)} \leq C(b,\beta) \delta$ for any $h \in \mathcal{M}_{disc,2}(C)$, which implies Theorem~\ref{thm:main} when combined with Lemma~\ref{lem:convdet}. Since $\delta (k_2(Y) + i_2) \leq Y_2 + 4\delta$ for $i \in I$, applying the Stein factor bounds in Proposition~\ref{prop:diffbounds} with the bounds on  $\varepsilon_1(Y)$ and $\varepsilon_2(Y)$ in Proposition~\ref{prop:main}   yields
\begin{align}
\abs{\varepsilon_1(Y)} \leq&\ C(b,\beta) 1(Y \in B)  \delta(1+Y_2)^{3},  \quad \abs{\varepsilon_2(Y)} \leq  1(Y \not \in B) C(\beta) \delta^{-3} (1 + Y_1+Y_2)^3. \label{eq:e12}
\end{align}
We point out that 
\begin{align}
\delta^{-1} \leq C(Y_1 + Y_2), \quad \text{ for any } Y \not \in B, \label{eq:deltaineq}
\end{align}
which follows from the facts that  $Y_1 + Y_2 \geq \delta(n/2-8) = \delta^{-1}/2 - \delta$ for $Y \not \in B$, that $\delta = 1/\sqrt{n}$,  and  that $n > 16$. Combining \eqref{eq:e12}, \eqref{eq:deltaineq}, and the fact that the moments of $Y_i$ are finite  yields 
\begin{align*}
\E \abs{\varepsilon_1(Y)} + \E \abs{\varepsilon_2(Y)} \leq C(b,\beta) \delta \E (1+Y_1+Y_2)^{7} \leq C(b,\beta)\delta.
\end{align*}
Furthermore, applying the Stein factor bounds in Proposition~\ref{prop:diffbounds} to  the bounds on $\varepsilon_3(Y)$ and $\varepsilon_4(Y)$ in Proposition~\ref{prop:main}, and using \eqref{eq:deltaineq}, we get 
\begin{align*}
\abs{\varepsilon_3(Y) +\varepsilon_4(Y)} \leq&\ C(b,\beta) 1(Y \in B)   \delta(1+Y_2)^2 + C(b,\beta) 1(Y \not \in B)   \delta^{-2}(1 + Y_1 + Y_2)^2 \\
\leq&\ C(b,\beta) \delta(1 + Y_1 + Y_2)^{5}.
\end{align*}
Thus, \eqref{eq:ubound} of Lemma~\ref{lem:ymoments} implies that 
\begin{align*}
\E \Big( \int_{0}^{1}\big|\varepsilon_3(Y(s)) + \varepsilon_4(Y(s)) \big|1(Y_1(s) = 0)  dU(s) \Big) \leq C(b,\beta)\delta.
\end{align*}
\finishproof

\section{Stein Factor Bounds}
\label{sec:perturbation}
\blue{In this section we prove Proposition~\ref{prop:diffbounds}. We bound the first-order differences in Section~\ref{sec:first_diff}. This requires the most effort. The second-order differences are bounded at the start of Section~\ref{sec:high_diff}, with Section~\ref{sec:moments}  showing how they can be used to bound $\E |h(X)|$, which may be of independent interest. Section~\ref{sec:third_bounds} contains the third-order bounds and Section~\ref{sec:proof_relate} proves two technical lemmas needed for the second-order bounds. }

\subsection{First-Order Differences}
\label{sec:first_diff}
In this section we bound
\begin{align*}
\Delta_{i} f_h(x^q) = \int_{0}^{\infty} \E_{x^q+\delta e^{(i)}}  h(X(t)) - \E_{x^q}  h(X(t)) dt
\end{align*}
by coupling two copies of the JSQ model initialized one customer apart. The coupling is introduced in the following lemma, which is stated in terms of the unscaled CTMC $\{Q(t)\}$.
\begin{lemma}
\label{lem:coupling}
For $1 \leq i \leq b+1$, define $\Theta_{i}^{Q} =  \{ (q,\widetilde q) \in S_Q \times S_Q : q_i <n,\ \widetilde q_i = q_i+1 \}$. There exists a coupling $\{\widetilde Q(t)\}$ of $\{Q(t)\}$  whose transient distribution satisfies
\begin{align}
&\{\widetilde Q(t) |  (Q(0),\widetilde Q(0)) \in \Theta_i^{Q},\ Q(0) = q \}_{t \geq 0} \stackrel{d}{=} \{Q(t) | Q(0) = (q+   e^{(i)})\}. \label{eq:marginal}
\end{align}
Furthermore, if $(Q(0),\widetilde Q(0)) \in \bigcup_{i=1}^{b+1} \Theta_i^{Q}$, then
\begin{enumerate}[(a)]
\item   $\widetilde Q(t) =  Q(t)$ for all times  $t \geq \tau_C$, where $\tau_C = \inf \{t \geq 0: Q(t) = \widetilde Q(t) \}$.
\item The pair $(Q(t),\widetilde Q(t))$  belongs to $\bigcup_{i=1}^{b+1} \Theta_i^{Q}$ for all times $t < \tau_C$. 
\item Let $V$ be a unit-mean exponentially distributed random variable  independent of $\{Q(t)\}$. Then
\begin{align}
\tau_C \stackrel{d}{=}  \min \bigg\{ \inf_{t \geq 0} \Big\{ \int_{0}^{t} 1\big( (Q(s),\widetilde Q(s)) \in \Theta_1^{Q}\big) ds = V  \Big\},\  \inf_{t \geq 0} \Big\{Q_{b+1}(t) =n \Big\} \bigg\}. \label{eq:coup_time_alternative}
\end{align}
\end{enumerate} 
\end{lemma}

\startproof{Proof of Lemma~\ref{lem:coupling}}
Let us construct a joint CTMC $\{ (Q(t),\widetilde Q(t))\}$ by specifying its transitions.  For simplicity, we refer to $\{Q(t)\}$ as system $1$ and to $\{\widetilde Q(t)\}$ as system $2$. We think of system 2 as a copy of system 1  but with an additional low-priority customer following a preemptive resume rule. That is, service is interrupted, and the extra customer moves to the back of its buffer when a regular customer joins, even if the low-priority customer is currently in service.

\begin{figure}
\centering 
\scalebox{0.95}{
\begin{tikzpicture}
 [server/.style={circle, inner sep=1.0mm, minimum width=.8cm,
    draw=green,fill=green!10,thick},
    customer/.style={circle, inner sep=0.5mm, minimum width=.4cm,
    draw=red,fill=red!10,thick},
    blue/.style={circle, inner sep=0.5mm, minimum width=.4cm,
    draw=blue,fill=blue!10,thick},
    green/.style={circle, inner sep=0.5mm, minimum width=.4cm,
    draw=green,fill=green!10,thick},
  buffer/.style={rectangle, rounded corners=3pt,
    inner sep=0.0mm, minimum width=.9cm,  minimum height=.8cm,
    draw=orange,fill=blue!10,thick},
 vbuffer/.style={rectangle,rounded corners=3pt,
     inner sep=0.0mm,  minimum width=.6cm, minimum height=.8cm,
     draw=orange,fill=blue!10,thick}
]
	\node (ref){};
  \node[customer]  (C1) {};
  \node[customer]  (C2) [right =.05in of C1]   {};
  \node[customer]  (C3) [right =.05in of C2]   {};
  \node[blue]  (C4) [right =.05in of C3]   {};
  \node[customer]  (C5) [right =.05in of C4]   {};
  \node  (C6) [right =.05in of C5]   {};
  
  \node[customer]  (C11) [above =.1in of C1]   {};
  \node[customer]  (C21) [above =.1in of C2]   {};
  \node[customer]  (C31) [above =.1in of C3]   {};
%  \node[customer]  (C41) [above =.1in of C4]   {};
  \node  (C61) [above =.1in of C6]   {};
  
%  \node[customer]  (D1) [above =.1in of C11]   {};
%  \node[customer]  (D2) [above =.1in of C21]   {};
%  \node[customer]  (D3) [above =.1in of C31]   {};
%  \node[blue]  (D4) [above =.1in of C41]   {};
  \node (D6) [above =.1in of C61]   {};
  
%  \node[customer]  (E1) [above =.1in of D1]   {};
%  \node[customer]  (E2) [above =.1in of D2]   {};
%  \node[customer]  (E3) [above =.1in of D3]   {};
  \node (E6) [above =.1in of D6]   {};
  
%  \node[customer]  (F1) [above =.1in of E1]   {};
  \node (F6) [above =.1in of E6]   {};

 \draw[dashed]  ($(C1.north)+(-0.2in,.045in)$) edge ($(C5.north)+(0.7in,.045in)$);
 \draw[dashed]  ($(C1.west)+(-0.025in,-.1in)$) edge ($(C1.west)+(-0.025in,1.3in)$);
 \draw[dashed]  ($(C2.west)+(-0.025in,-.1in)$) edge ($(C2.west)+(-0.025in,1.3in)$);
 \draw[dashed]  ($(C3.west)+(-0.025in,-.1in)$) edge ($(C3.west)+(-0.025in,1.3in)$);
 \draw[dashed]  ($(C4.west)+(-0.025in,-.1in)$) edge ($(C4.west)+(-0.025in,1.3in)$);
 \draw[dashed]  ($(C5.west)+(-0.025in,-.1in)$) edge ($(C5.west)+(-0.025in,1.3in)$);
 \draw[dashed]  ($(C6.west)+(-0.025in,-.1in)$) edge ($(C6.west)+(-0.025in,1.3in)$);
 \draw[dashed]  ($(C6.east)+(0.065in,-.1in)$) edge ($(C6.east)+(0.065in,1.3in)$);
  
  \node[left, align=left] (X1) at ($(C6)+(.4in,0)$) {$q_1$};
  \node[left, align=left] (X2) at ($(C61)+(.4in,.025in)$) {$q_2$}; 
  \node[left, align=left] (X3) at ($(D6)+(.4in,.09in)$) {$q_3$}; 
  \node[left, align=left] (X4) at ($(E6)+(.4in,.15in)$) {$q_4$}; 
  \node[left, align=left] (X5) at ($(F6)+(.4in,.21in)$) {$q_5$}; 
 \end{tikzpicture}}
 \hspace{2cm}
\scalebox{0.95}{
\begin{tikzpicture}
 [server/.style={circle, inner sep=1.0mm, minimum width=.8cm,
    draw=green,fill=green!10,thick},
    customer/.style={circle, inner sep=0.5mm, minimum width=.4cm,
    draw=red,fill=red!10,thick},
    blue/.style={circle, inner sep=0.5mm, minimum width=.4cm,
    draw=blue,fill=blue!10,thick},
    green/.style={circle, inner sep=0.5mm, minimum width=.4cm,
    draw=green,fill=green!10,thick},
  buffer/.style={rectangle, rounded corners=3pt,
    inner sep=0.0mm, minimum width=.9cm,  minimum height=.8cm,
    draw=orange,fill=blue!10,thick},
 vbuffer/.style={rectangle,rounded corners=3pt,
     inner sep=0.0mm,  minimum width=.6cm, minimum height=.8cm,
     draw=orange,fill=blue!10,thick}
]
	\node (ref){};
  \node[customer]  (C1) {};
  \node[customer]  (C2) [right =.05in of C1]   {};
  \node[customer]  (C3) [right =.05in of C2]   {};
  \node[customer]  (C4) [right =.05in of C3]   {};
  \node[customer]  (C5) [right =.05in of C4]   {};
  \node  (C6) [right =.05in of C5]   {};
  
  \node[customer]  (C11) [above =.1in of C1]   {};
  \node[customer]  (C21) [above =.1in of C2]   {};
  \node[customer]  (C31) [above =.1in of C3]   {};
  \node[customer]  (C41) [above =.1in of C4]   {};
  \node  (C61) [above =.1in of C6]   {};
  
  \node[customer]  (D1) [above =.1in of C11]   {};
  \node[customer]  (D2) [above =.1in of C21]   {};
  \node[customer]  (D3) [above =.1in of C31]   {};
  \node[blue]  (D4) [above =.1in of C41]   {};
  \node (D6) [above =.1in of C61]   {};
  
  \node[customer]  (E1) [above =.1in of D1]   {};
  \node[customer]  (E2) [above =.1in of D2]   {};
  \node[customer]  (E3) [above =.1in of D3]   {};
  \node (E6) [above =.1in of D6]   {};
  
  \node[customer]  (F1) [above =.1in of E1]   {};
  \node (F6) [above =.1in of E6]   {};

 \draw[dashed]  ($(C1.north)+(-0.2in,.045in)$) edge ($(C5.north)+(0.7in,.045in)$);
 \draw[dashed]  ($(C1.west)+(-0.025in,-.1in)$) edge ($(C1.west)+(-0.025in,1.3in)$);
 \draw[dashed]  ($(C2.west)+(-0.025in,-.1in)$) edge ($(C2.west)+(-0.025in,1.3in)$);
 \draw[dashed]  ($(C3.west)+(-0.025in,-.1in)$) edge ($(C3.west)+(-0.025in,1.3in)$);
 \draw[dashed]  ($(C4.west)+(-0.025in,-.1in)$) edge ($(C4.west)+(-0.025in,1.3in)$);
 \draw[dashed]  ($(C5.west)+(-0.025in,-.1in)$) edge ($(C5.west)+(-0.025in,1.3in)$);
 \draw[dashed]  ($(C6.west)+(-0.025in,-.1in)$) edge ($(C6.west)+(-0.025in,1.3in)$);
 \draw[dashed]  ($(C6.east)+(0.065in,-.1in)$) edge ($(C6.east)+(0.065in,1.3in)$);
  
  \node[left, align=left] (X1) at ($(C6)+(.4in,0)$) {$q_1$};
  \node[left, align=left] (X2) at ($(C61)+(.4in,.025in)$) {$q_2$}; 
  \node[left, align=left] (X3) at ($(D6)+(.4in,.09in)$) {$q_3$}; 
  \node[left, align=left] (X4) at ($(E6)+(.4in,.15in)$) {$q_4$}; 
  \node[left, align=left] (X5) at ($(F6)+(.4in,.21in)$) {$q_5$}; 
\end{tikzpicture}}
 
\caption{Two possible states of the joint chain $(Q(t),\widetilde Q(t))$ are depicted. The red circles correspond to customers in $Q(t)$, while the blue circle is the extra customer in $\widetilde Q(t)$. \blue{In the figure on the left, the joint chain is in $\Theta_1^{Q}$, meaning the blue customer is in service and will leave the system  after an exponentially distributed amount of time, coupling the joint chain.}  In the figure on the right, the joint chain is in $\Theta_3^{Q}$  because the blue customer is assigned to a server with a total of three customers.  }
\label{fig:jointex}
\end{figure}
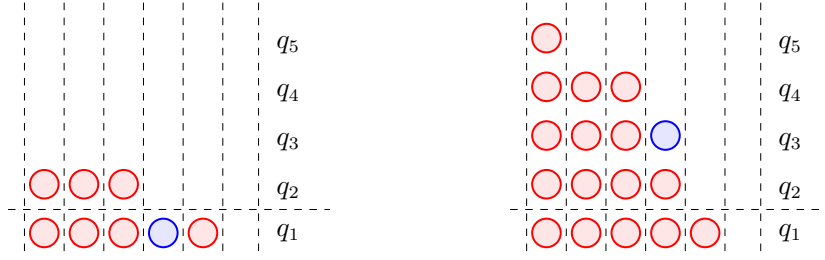

Any state in $\Theta_1^{Q}$  is one where   the low-priority customer is in service. The remaining $\Theta_i^{Q}$ correspond to states where the low-priority customer is assigned to a server with a total of $i$ customers;  Figure~\ref{fig:jointex} contains an example of a states in $\Theta_1^{Q}$ and $\Theta_3^{Q}$. Assuming $(Q(0),\widetilde Q(0)) = (q,\widetilde q) \in \Theta_i^{Q}$ for some $1 \leq i \leq b+1$, we now describe the possible transitions of the joint chain.

If $i = 1$, then the low-priority customer is in service. After a unit-mean exponentially distributed amount of time, he leaves system 2 and both systems couple. After coupling, systems 1 and 2 are identical in terms of current and future customers,   so they coincide on every sample path. All other transitions of the joint chain are based on the standard transitions of the JSQ model. In other words, a service completion by any of the $q_1$ servers working in system 1 results in a customer departure from both systems. 

\begin{figure}
\centering 
\scalebox{0.92}{
\begin{tikzpicture}
 [server/.style={circle, inner sep=1.0mm, minimum width=.8cm,
    draw=green,fill=green!10,thick},
    customer/.style={circle, inner sep=0.5mm, minimum width=.4cm,
    draw=red,fill=red!10,thick},
    blue/.style={circle, inner sep=0.5mm, minimum width=.4cm,
    draw=blue,fill=blue!10,thick},
    green/.style={circle, inner sep=0.5mm, minimum width=.4cm,
    draw=green,fill=green!10,thick},
  buffer/.style={rectangle, rounded corners=3pt,
    inner sep=0.0mm, minimum width=.9cm,  minimum height=.8cm,
    draw=orange,fill=blue!10,thick},
 vbuffer/.style={rectangle,rounded corners=3pt,
     inner sep=0.0mm,  minimum width=.6cm, minimum height=.8cm,
     draw=orange,fill=blue!10,thick}
]
	\node (ref){};
  \node[customer]  (C1) {};
  \node[customer]  (C2) [right =.05in of C1]   {};
  \node[customer]  (C3) [right =.05in of C2]   {};
  \node[blue]  (C4) [right =.05in of C3]   {};
  \node[customer]  (C5) [right =.05in of C4]   {};
  \node  (C6) [right =.05in of C5]   {};
  
  \node[customer]  (C11) [above =.1in of C1]   {};
  \node[customer]  (C21) [above =.1in of C2]   {};
  \node[customer]  (C31) [above =.1in of C3]   {};
%  \node[customer]  (C41) [above =.1in of C4]   {};
  \node  (C61) [above =.1in of C6]   {};
  
%  \node[customer]  (D1) [above =.1in of C11]   {};
%  \node[customer]  (D2) [above =.1in of C21]   {};
%  \node[customer]  (D3) [above =.1in of C31]   {};
%  \node[blue]  (D4) [above =.1in of C41]   {};
  \node (D6) [above =.1in of C61]   {};
  
%  \node[customer]  (E1) [above =.1in of D1]   {};
%  \node[customer]  (E2) [above =.1in of D2]   {};
%  \node[customer]  (E3) [above =.1in of D3]   {};
  \node (E6) [above =.1in of D6]   {};
  
%  \node[customer]  (F1) [above =.1in of E1]   {};
  \node (F6) [above =.1in of E6]   {};

 \draw[dashed]  ($(C1.north)+(-0.2in,.045in)$) edge ($(C5.north)+(0.7in,.045in)$);
 \draw[dashed]  ($(C1.west)+(-0.025in,-.1in)$) edge ($(C1.west)+(-0.025in,1.3in)$);
 \draw[dashed]  ($(C2.west)+(-0.025in,-.1in)$) edge ($(C2.west)+(-0.025in,1.3in)$);
 \draw[dashed]  ($(C3.west)+(-0.025in,-.1in)$) edge ($(C3.west)+(-0.025in,1.3in)$);
 \draw[dashed]  ($(C4.west)+(-0.025in,-.1in)$) edge ($(C4.west)+(-0.025in,1.3in)$);
 \draw[dashed]  ($(C5.west)+(-0.025in,-.1in)$) edge ($(C5.west)+(-0.025in,1.3in)$);
 \draw[dashed]  ($(C6.west)+(-0.025in,-.1in)$) edge ($(C6.west)+(-0.025in,1.3in)$);
 \draw[dashed]  ($(C6.east)+(0.065in,-.1in)$) edge ($(C6.east)+(0.065in,1.3in)$);
  
  \node[left, align=left] (X1) at ($(C6)+(.4in,0)$) {$q_1$};
  \node[left, align=left] (X2) at ($(C61)+(.4in,.025in)$) {$q_2$}; 
  \node[left, align=left] (X3) at ($(D6)+(.4in,.09in)$) {$q_3$}; 
  \node[left, align=left] (X4) at ($(E6)+(.4in,.15in)$) {$q_4$}; 
  \node[left, align=left] (X5) at ($(F6)+(.4in,.21in)$) {$q_5$}; 
 \end{tikzpicture}}
 \hspace{.75cm} 
\scalebox{0.92}{
\begin{tikzpicture}
 [server/.style={circle, inner sep=1.0mm, minimum width=.8cm,
    draw=green,fill=green!10,thick},
    customer/.style={circle, inner sep=0.5mm, minimum width=.4cm,
    draw=red,fill=red!10,thick},
    blue/.style={circle, inner sep=0.5mm, minimum width=.4cm,
    draw=blue,fill=blue!10,thick},
    green/.style={circle, inner sep=0.5mm, minimum width=.4cm,
    draw=green,fill=green!10,thick},
  buffer/.style={rectangle, rounded corners=3pt,
    inner sep=0.0mm, minimum width=.9cm,  minimum height=.8cm,
    draw=orange,fill=blue!10,thick},
 vbuffer/.style={rectangle,rounded corners=3pt,
     inner sep=0.0mm,  minimum width=.6cm, minimum height=.8cm,
     draw=orange,fill=blue!10,thick}
]
	\node (ref){};
  \node[customer]  (C1) {};
  \node[customer]  (C2) [right =.05in of C1]   {};
  \node[customer]  (C3) [right =.05in of C2]   {};
  \node[blue]  (C4) [right =.05in of C3]   {};
  \node[customer]  (C5) [right =.05in of C4]   {};
  \node[customer]  (C6) [right =.05in of C5]   {};
  
  \node[customer]  (C11) [above =.1in of C1]   {};
  \node[customer]  (C21) [above =.1in of C2]   {};
  \node[customer]  (C31) [above =.1in of C3]   {};
%  \node[customer]  (C41) [above =.1in of C4]   {};
  \node  (C61) [above =.1in of C6]   {};
  
%  \node[customer]  (D1) [above =.1in of C11]   {};
%  \node[customer]  (D2) [above =.1in of C21]   {};
%  \node[customer]  (D3) [above =.1in of C31]   {};
%  \node[blue]  (D4) [above =.1in of C41]   {};
  \node (D6) [above =.1in of C61]   {};
  
%  \node[customer]  (E1) [above =.1in of D1]   {};
%  \node[customer]  (E2) [above =.1in of D2]   {};
%  \node[customer]  (E3) [above =.1in of D3]   {};
  \node (E6) [above =.1in of D6]   {};
  
%  \node[customer]  (F1) [above =.1in of E1]   {};
  \node (F6) [above =.1in of E6]   {};

 \draw[dashed]  ($(C1.north)+(-0.2in,.045in)$) edge ($(C5.north)+(0.7in,.045in)$);
 \draw[dashed]  ($(C1.west)+(-0.025in,-.1in)$) edge ($(C1.west)+(-0.025in,1.3in)$);
 \draw[dashed]  ($(C2.west)+(-0.025in,-.1in)$) edge ($(C2.west)+(-0.025in,1.3in)$);
 \draw[dashed]  ($(C3.west)+(-0.025in,-.1in)$) edge ($(C3.west)+(-0.025in,1.3in)$);
 \draw[dashed]  ($(C4.west)+(-0.025in,-.1in)$) edge ($(C4.west)+(-0.025in,1.3in)$);
 \draw[dashed]  ($(C5.west)+(-0.025in,-.1in)$) edge ($(C5.west)+(-0.025in,1.3in)$);
 \draw[dashed]  ($(C6.west)+(-0.025in,-.1in)$) edge ($(C6.west)+(-0.025in,1.3in)$);
 \draw[dashed]  ($(C6.east)+(0.065in,-.1in)$) edge ($(C6.east)+(0.065in,1.3in)$);
  
  \node[left, align=left] (X1) at ($(C6)+(.4in,0)$) {$q_1$};
  \node[left, align=left] (X2) at ($(C61)+(.4in,.025in)$) {$q_2$}; 
  \node[left, align=left] (X3) at ($(D6)+(.4in,.09in)$) {$q_3$}; 
  \node[left, align=left] (X4) at ($(E6)+(.4in,.15in)$) {$q_4$}; 
  \node[left, align=left] (X5) at ($(F6)+(.4in,.21in)$) {$q_5$}; 
 \end{tikzpicture}}
 \hspace{.75cm}
\scalebox{0.92}{
\begin{tikzpicture}
 [server/.style={circle, inner sep=1.0mm, minimum width=.8cm,
    draw=green,fill=green!10,thick},
    customer/.style={circle, inner sep=0.5mm, minimum width=.4cm,
    draw=red,fill=red!10,thick},
    blue/.style={circle, inner sep=0.5mm, minimum width=.4cm,
    draw=blue,fill=blue!10,thick},
    green/.style={circle, inner sep=0.5mm, minimum width=.4cm,
    draw=green,fill=green!10,thick},
  buffer/.style={rectangle, rounded corners=3pt,
    inner sep=0.0mm, minimum width=.9cm,  minimum height=.8cm,
    draw=orange,fill=blue!10,thick},
 vbuffer/.style={rectangle,rounded corners=3pt,
     inner sep=0.0mm,  minimum width=.6cm, minimum height=.8cm,
     draw=orange,fill=blue!10,thick}
]
	\node (ref){};
  \node[customer]  (C1) {};
  \node[customer]  (C2) [right =.05in of C1]   {};
  \node[customer]  (C3) [right =.05in of C2]   {};
  \node[customer]  (C4) [right =.05in of C3]   {};
  \node[customer]  (C5) [right =.05in of C4]   {};
  \node[customer]  (C6) [right =.05in of C5]   {};
  
  \node[customer]  (C11) [above =.1in of C1]   {};
  \node[customer]  (C21) [above =.1in of C2]   {};
  \node[customer]  (C31) [above =.1in of C3]   {};
  \node[blue]  (C41) [above =.1in of C4]   {};
  \node  (C61) [above =.1in of C6]   {};
  
%  \node[customer]  (D1) [above =.1in of C11]   {};
%  \node[customer]  (D2) [above =.1in of C21]   {};
%  \node[customer]  (D3) [above =.1in of C31]   {};
%  \node[blue]  (D4) [above =.1in of C41]   {};
  \node (D6) [above =.1in of C61]   {};
  
%  \node[customer]  (E1) [above =.1in of D1]   {};
%  \node[customer]  (E2) [above =.1in of D2]   {};
%  \node[customer]  (E3) [above =.1in of D3]   {};
  \node (E6) [above =.1in of D6]   {};
  
%  \node[customer]  (F1) [above =.1in of E1]   {};
  \node (F6) [above =.1in of E6]   {};

 \draw[dashed]  ($(C1.north)+(-0.2in,.045in)$) edge ($(C5.north)+(0.7in,.045in)$);
 \draw[dashed]  ($(C1.west)+(-0.025in,-.1in)$) edge ($(C1.west)+(-0.025in,1.3in)$);
 \draw[dashed]  ($(C2.west)+(-0.025in,-.1in)$) edge ($(C2.west)+(-0.025in,1.3in)$);
 \draw[dashed]  ($(C3.west)+(-0.025in,-.1in)$) edge ($(C3.west)+(-0.025in,1.3in)$);
 \draw[dashed]  ($(C4.west)+(-0.025in,-.1in)$) edge ($(C4.west)+(-0.025in,1.3in)$);
 \draw[dashed]  ($(C5.west)+(-0.025in,-.1in)$) edge ($(C5.west)+(-0.025in,1.3in)$);
 \draw[dashed]  ($(C6.west)+(-0.025in,-.1in)$) edge ($(C6.west)+(-0.025in,1.3in)$);
 \draw[dashed]  ($(C6.east)+(0.065in,-.1in)$) edge ($(C6.east)+(0.065in,1.3in)$);
  
  \node[left, align=left] (X1) at ($(C6)+(.4in,0)$) {$q_1$};
  \node[left, align=left] (X2) at ($(C61)+(.4in,.025in)$) {$q_2$}; 
  \node[left, align=left] (X3) at ($(D6)+(.4in,.09in)$) {$q_3$}; 
  \node[left, align=left] (X4) at ($(E6)+(.4in,.15in)$) {$q_4$}; 
  \node[left, align=left] (X5) at ($(F6)+(.4in,.21in)$) {$q_5$}; 
 \end{tikzpicture}}
 
\caption{\blue{From left to right, the figures depict the arrival of two customers. The second arrival results in a transition from $\Theta_1^{Q}$ to $\Theta_2^{Q}$.} }
\label{fig:arrivals}
\end{figure}
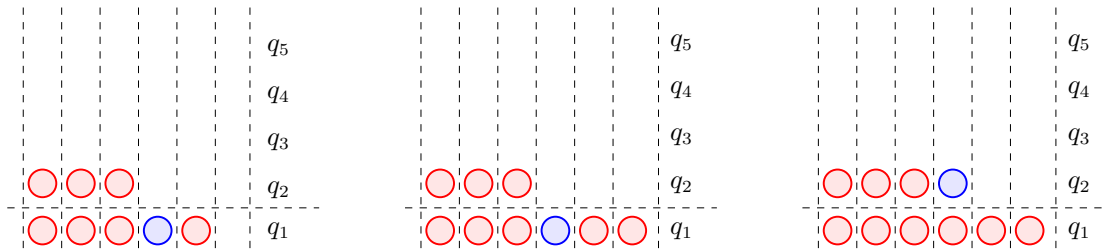

\blue{Figure~\ref{fig:arrivals}  illustrates the effect of arrivals when $(q,\widetilde q) \in \Theta_1^{Q}$.} Namely, when $q_1 \leq n-2$, a new arrival is assigned to the same idle server in both systems. If a customer arrives when $q_1 = n-1$, then  system 1  has only one idle server  and system 2 has none.  In system 1, that customer will be assigned to the last remaining idle server.  Recall that when defining our JSQ model, we allowed for an arbitrary tie-breaking decision in routing arrivals.  Therefore, in system 2, we assign that customer  to the server working on the low-priority customer, causing a service preemption and pushing the low-priority customer to the back of the buffer. An arrival when $q_1 = n-1$ transitions the joint chain from $\Theta_1^{Q}$ to $\Theta_2^{Q}$.

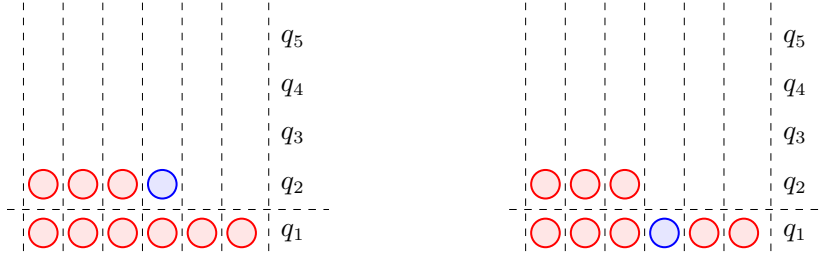
\begin{figure}
\centering 
\scalebox{0.95}{
\begin{tikzpicture}
 [server/.style={circle, inner sep=1.0mm, minimum width=.8cm,
    draw=green,fill=green!10,thick},
    customer/.style={circle, inner sep=0.5mm, minimum width=.4cm,
    draw=red,fill=red!10,thick},
    blue/.style={circle, inner sep=0.5mm, minimum width=.4cm,
    draw=blue,fill=blue!10,thick},
    green/.style={circle, inner sep=0.5mm, minimum width=.4cm,
    draw=green,fill=green!10,thick},
  buffer/.style={rectangle, rounded corners=3pt,
    inner sep=0.0mm, minimum width=.9cm,  minimum height=.8cm,
    draw=orange,fill=blue!10,thick},
 vbuffer/.style={rectangle,rounded corners=3pt,
     inner sep=0.0mm,  minimum width=.6cm, minimum height=.8cm,
     draw=orange,fill=blue!10,thick}
]
	\node (ref){};
  \node[customer]  (C1) {};
  \node[customer]  (C2) [right =.05in of C1]   {};
  \node[customer]  (C3) [right =.05in of C2]   {};
  \node[customer]  (C4) [right =.05in of C3]   {};
  \node[customer]  (C5) [right =.05in of C4]   {};
  \node[customer]  (C6) [right =.05in of C5]   {};
  
  \node[customer]  (C11) [above =.1in of C1]   {};
  \node[customer]  (C21) [above =.1in of C2]   {};
  \node[customer]  (C31) [above =.1in of C3]   {};
  \node[blue]  (C41) [above =.1in of C4]   {};
  \node  (C61) [above =.1in of C6]   {};
  
%  \node[customer]  (D1) [above =.1in of C11]   {};
%  \node[customer]  (D2) [above =.1in of C21]   {};
%  \node[customer]  (D3) [above =.1in of C31]   {};
%  \node[blue]  (D4) [above =.1in of C41]   {};
  \node (D6) [above =.1in of C61]   {};
  
%  \node[customer]  (E1) [above =.1in of D1]   {};
%  \node[customer]  (E2) [above =.1in of D2]   {};
%  \node[customer]  (E3) [above =.1in of D3]   {};
  \node (E6) [above =.1in of D6]   {};
  
%  \node[customer]  (F1) [above =.1in of E1]   {};
  \node (F6) [above =.1in of E6]   {};

 \draw[dashed]  ($(C1.north)+(-0.2in,.045in)$) edge ($(C5.north)+(0.7in,.045in)$);
 \draw[dashed]  ($(C1.west)+(-0.025in,-.1in)$) edge ($(C1.west)+(-0.025in,1.3in)$);
 \draw[dashed]  ($(C2.west)+(-0.025in,-.1in)$) edge ($(C2.west)+(-0.025in,1.3in)$);
 \draw[dashed]  ($(C3.west)+(-0.025in,-.1in)$) edge ($(C3.west)+(-0.025in,1.3in)$);
 \draw[dashed]  ($(C4.west)+(-0.025in,-.1in)$) edge ($(C4.west)+(-0.025in,1.3in)$);
 \draw[dashed]  ($(C5.west)+(-0.025in,-.1in)$) edge ($(C5.west)+(-0.025in,1.3in)$);
 \draw[dashed]  ($(C6.west)+(-0.025in,-.1in)$) edge ($(C6.west)+(-0.025in,1.3in)$);
 \draw[dashed]  ($(C6.east)+(0.065in,-.1in)$) edge ($(C6.east)+(0.065in,1.3in)$);
  
  \node[left, align=left] (X1) at ($(C6)+(.4in,0)$) {$q_1$};
  \node[left, align=left] (X2) at ($(C61)+(.4in,.025in)$) {$q_2$}; 
  \node[left, align=left] (X3) at ($(D6)+(.4in,.09in)$) {$q_3$}; 
  \node[left, align=left] (X4) at ($(E6)+(.4in,.15in)$) {$q_4$}; 
  \node[left, align=left] (X5) at ($(F6)+(.4in,.21in)$) {$q_5$}; 
 \end{tikzpicture}}
 \hspace{2cm}
\scalebox{0.95}{
\begin{tikzpicture}
 [server/.style={circle, inner sep=1.0mm, minimum width=.8cm,
    draw=green,fill=green!10,thick},
    customer/.style={circle, inner sep=0.5mm, minimum width=.4cm,
    draw=red,fill=red!10,thick},
    blue/.style={circle, inner sep=0.5mm, minimum width=.4cm,
    draw=blue,fill=blue!10,thick},
    green/.style={circle, inner sep=0.5mm, minimum width=.4cm,
    draw=green,fill=green!10,thick},
  buffer/.style={rectangle, rounded corners=3pt,
    inner sep=0.0mm, minimum width=.9cm,  minimum height=.8cm,
    draw=orange,fill=blue!10,thick},
 vbuffer/.style={rectangle,rounded corners=3pt,
     inner sep=0.0mm,  minimum width=.6cm, minimum height=.8cm,
     draw=orange,fill=blue!10,thick}
]
	\node (ref){};
  \node[customer]  (C1) {};
  \node[customer]  (C2) [right =.05in of C1]   {};
  \node[customer]  (C3) [right =.05in of C2]   {};
  \node[blue]  (C4) [right =.05in of C3]   {};
  \node[customer]  (C5) [right =.05in of C4]   {};
  \node[customer]  (C6) [right =.05in of C5]   {};
  
  \node[customer]  (C11) [above =.1in of C1]   {};
  \node[customer]  (C21) [above =.1in of C2]   {};
  \node[customer]  (C31) [above =.1in of C3]   {};
%  \node[blue]  (C41) [above =.1in of C4]   {};
  \node  (C61) [above =.1in of C6]   {};
  
%  \node[customer]  (D1) [above =.1in of C11]   {};
%  \node[customer]  (D2) [above =.1in of C21]   {};
%  \node[customer]  (D3) [above =.1in of C31]   {};
%  \node[blue]  (D4) [above =.1in of C41]   {};
  \node (D6) [above =.1in of C61]   {};
  
%  \node[customer]  (E1) [above =.1in of D1]   {};
%  \node[customer]  (E2) [above =.1in of D2]   {};
%  \node[customer]  (E3) [above =.1in of D3]   {};
  \node (E6) [above =.1in of D6]   {};
  
%  \node[customer]  (F1) [above =.1in of E1]   {};
  \node (F6) [above =.1in of E6]   {};

 \draw[dashed]  ($(C1.north)+(-0.2in,.045in)$) edge ($(C5.north)+(0.7in,.045in)$);
 \draw[dashed]  ($(C1.west)+(-0.025in,-.1in)$) edge ($(C1.west)+(-0.025in,1.3in)$);
 \draw[dashed]  ($(C2.west)+(-0.025in,-.1in)$) edge ($(C2.west)+(-0.025in,1.3in)$);
 \draw[dashed]  ($(C3.west)+(-0.025in,-.1in)$) edge ($(C3.west)+(-0.025in,1.3in)$);
 \draw[dashed]  ($(C4.west)+(-0.025in,-.1in)$) edge ($(C4.west)+(-0.025in,1.3in)$);
 \draw[dashed]  ($(C5.west)+(-0.025in,-.1in)$) edge ($(C5.west)+(-0.025in,1.3in)$);
 \draw[dashed]  ($(C6.west)+(-0.025in,-.1in)$) edge ($(C6.west)+(-0.025in,1.3in)$);
 \draw[dashed]  ($(C6.east)+(0.065in,-.1in)$) edge ($(C6.east)+(0.065in,1.3in)$);
  
  \node[left, align=left] (X1) at ($(C6)+(.4in,0)$) {$q_1$};
  \node[left, align=left] (X2) at ($(C61)+(.4in,.025in)$) {$q_2$}; 
  \node[left, align=left] (X3) at ($(D6)+(.4in,.09in)$) {$q_3$}; 
  \node[left, align=left] (X4) at ($(E6)+(.4in,.15in)$) {$q_4$}; 
  \node[left, align=left] (X5) at ($(F6)+(.4in,.21in)$) {$q_5$}; 
 \end{tikzpicture}}
 
\caption{\blue{From left to right, the server containing the blue customer in its buffer completes service, resulting in a transition from $\Theta_2^{Q}$ to $\Theta_1^{Q}$.}}
\label{fig:departure}
\end{figure}

If $2 \leq i \leq b$, then the low-priority customer is in the back of some server's buffer.   A service completion by any of the $q_1$ servers working in system 1 results in a customer departure from both systems. If, however, the service completion happens at the server containing the low-priority customer, then the chain transitions from $\Theta_i^{Q}$ to $\Theta_{i-1}^{Q}$  because the low-priority customer is now assigned to a server with $i-1$ customers; \blue{see Figure~\ref{fig:departure} for a depiction of such a transition}. All new arrivals get assigned to the same server in each system. Note that if an arrival happens when $q_i = n-1$ and $q_1=\cdots=q_{i-1}=n$, then the system transitions from $\Theta_{i}^{Q}$ to $\Theta_{i+1}^{Q}$.

The final case is when $i = b+1$. All transitions are identical to the $2 \leq i \leq b$ case, except for a customer arrival to a system where $q_{b+1}=n-1$ and $q_{1}=\ldots=q_{b} = n$. In that case, system 1 assigns the customer to the last available slot, but system 2 blocks the customer because it is already full. This transition causes the two systems to couple.  Note that our construction immediately implies the three claims in Lemma~\ref{lem:coupling}.
%Let us now verify the claims of the lemma. Note that  $\{\widetilde Q(t)\}$ evolves just like a JSQ system, and therefore \eqref{eq:marginal} is true. By construction, $\widetilde Q(t) =  Q(t)$ for all $t \geq \tau_C$ and provided the joint chain starts in $\bigcup_{i=1}^{b+1} \Theta_i$, it remains in $\bigcup_{i=1}^{b+1} \Theta_i$ for all times until coupling. Lastly, \eqref{eq:coup_time_alternative} is true because coupling can only happen in one of two ways. Either both systems become full, or the low-priority customer is served. The latter event happens once the joint chain spends an exponentially distributed amount of time in $\Theta_1^{Q}$.
\finishproof
\noindent Let $\widetilde X(t) = \big(\delta(n-\widetilde Q_1(t)), \delta \widetilde Q_2(t), \ldots, \delta \widetilde Q_{b+1}(t)\big) $  be the scaled version of $\widetilde Q(t)$. For any $x^q \in S$ with $x^q_1 > 0$,   and any $h \in \mathcal{M}_{disc,2}(C)$,  
\begin{align}
\abs{\int_{0}^{\infty} \E_{x-\delta e^{(1)}}  h(X(t)) - \E_{x}  h(X(t)) dt} =&\ \abs{\int_{0}^{\infty} \E_{(x, x-\delta e^{(1)})} \big( h(\widetilde X(t)) -  h(X(t)) \big) dt} \notag \\
\leq&\ \abs{\int_{0}^{\infty} \E_{(x, x-\delta e^{(1)})} \big(  \delta  1(t\leq \tau_C)\big) dt} = \delta \E_{(x, x-\delta e^{(1)})} \tau_C, \label{eq:first_diff_integral}
\end{align}
where $\E_{(x, x-\delta e^{(1)})}(\cdot)$ denotes the expectation given $(X(0),\widetilde X(0)) = (x, x-\delta e^{(1)})$. The inequality above is true because  the gap between $\{X(t)\}$ and $\{\widetilde X(t)\}$ never increases beyond one customer. The same argument implies that $\abs{ \Delta_{i} f_h(x^q)} \leq  \delta \E_{(x, x+\delta e^{(i)})} \tau_C$ for $i \geq 2$, and we see that bounding the first-order Stein factors amounts to bounding the expected coupling time $\tau_C$. The following lemma provides the necessary bound. It is worth highlighting that proving  this  result requires a large amount of effort  and JSQ-model-specific insight.
\begin{lemma}
\label{lem:coupling_time_bound}
For any $(q,\widetilde q) \in \bigcup_{i=1}^{b+1} \Theta_i^{Q}$,
\begin{align*}
\E_{(q, \widetilde q)} \tau_C \leq C(b,\beta)(1+\delta q_2).
\end{align*}
\end{lemma}
Before proving the lemma, we note that the first-order bounds in Proposition~\ref{prop:diffbounds} are a consequence of \eqref{eq:first_diff_integral} and  Lemma~\ref{lem:coupling_time_bound}; i.e., 
\begin{align}
\big|\Delta_{i}  f_{h}(x^q)\big|  \leq C(b,\beta)\delta(1+x^q_2), \quad   i = 1,\ldots, b+1. \label{eq:bound_first_difference}
\end{align}
Furthermore,  note that for any $x^q \in S$ with $x^q_3 = 0$, 
\begin{align*}
  f_h(x^q)  =&\  f_h(0) + \sum_{j_1=0}^{x^q_1/\delta-1} \Delta_1 f_h\big(\delta j_1,0,\ldots,0\big)+ \sum_{j_2=0}^{x^q_2/\delta-1} \Delta_2 f_h\big(x^q_1,\delta j_2,0,\ldots,0\big).
\end{align*}
Recall that $f_h(0) = 0$, and that the definition of $S$  implies that $\delta(j_1,j_2,0,\ldots,0) \in S$ for any $0 \leq j_1 \leq x^q_1/\delta$ and $0 \leq j_2 \leq x^q_2/\delta$. Combining these facts with \eqref{eq:bound_first_difference} yields
\begin{align}
 \abs{f_h(x^q)}  \leq&\  C(b,\beta)(1+x^q_2) (x^q_1 + x^q_2)/\delta, \quad x^q \in S,\ x^q_3 = 0, \label{eq:fbound}
\end{align}
which proves   one of the claims from Proposition~\ref{prop:diffbounds}. 

 We now describe the main idea and introduce several auxiliary lemmas  used to prove Lemma~\ref{lem:coupling_time_bound}. Our discussion communicates the main intuition behind the proof, leaving the technical details  to Appendix~\ref{sec:appendixsteinfactors}.
%
%To bound $\E_{(q, \widetilde q)} \tau_C$, we set up a  sequence of cycles, where the expected length of a cycle is bounded uniformly in $n$, and the chance of the joint chain coupling during a single cycle is bounded from below uniformly in $n$. The number of cycles until coupling occurs is then upper bounded by a geometric random variable whose parameter does not depend on $n$.   
Let $\gamma > 0$ be a constant independent of $n$ whose precise value will be specified later, and define 
\begin{align*}
\theta_{1} = n-\lfloor \sqrt{n}\beta/2 \rfloor, \quad \text{ and } \quad \theta_{2} =  \lfloor \gamma\sqrt{n}\rfloor.
\end{align*}
 Additionally, we define the stopping times
\begin{align*} 
\tau_{i}(q_i) =&\ \inf \{ t \geq 0 : Q_i(t) = q_i\}, \quad q_i\in \{0,1,\ldots,n\},\ i = 1,2.
\end{align*}
We now describe a sequence of cycles, or attempts, such that in each cycle, the probability of the joint chain coupling is bounded from below by a constant independent of $n$.  Given an initial state  $(Q(0), \widetilde Q(0))  = (q,\widetilde q)$ belonging to some $\Theta_i^{Q}$,  we wait until $\tau_{2}(\theta_{2})$, which marks the start of the first cycle.   From that point,  we wait until   $\min (\tau_{1}(\theta_{1}), \tau_{2}(2\theta_{2}))$. If $\tau_{1}(\theta_{1}) \geq \tau_{2}(2\theta_{2})$, then we give up trying to couple this cycle, and wait until $\tau_2(\theta_{2})$ to start a fresh cycle. If  $\tau_{1}(\theta_{1}) < \tau_{2}(2\theta_{2})$, then there are $\lfloor \sqrt{n}\beta/2 \rfloor$ idle servers and at most $2\theta_{2}$ non-empty buffers. From such a state, we are guaranteed that coupling happens if the joint CTMC enters $\Theta_{1}^{Q}$ and spends an exponentially distributed amount of time there  before all servers in $\{Q(t)\}$ become busy; i.e., $\tau_C < \tau_{1}(n)$. If $\tau_C \geq \tau_{1}(n)$, we give up trying to couple this cycle and wait until $\tau_2(\theta_{2})$ for the next cycle to restart the coupling attempt. Note that this cycle sequence resembles a renewal sequence, but the new cycle times are not renewal times because the values of $Q_3(\cdot), \ldots, Q_{b+1}(\cdot)$ can vary at the start of each new cycle. 

From our discussion, it follows that coupling is guaranteed in any given cycle if, starting from a state with $q_2 = \theta_{2}$, the events $\{ \tau_{1}(\theta_{1}) < \tau_{2}(2\theta_{2}) \}$ and  $ \{\tau_C < \tau_{1}(n) \}$ occur. In Appendix~\ref{sec:appendixsteinfactors} we  derive a lower bound, uniform in $n$, on the probability of coupling in a given cycle, implying that coupling is guaranteed to happen after a geometrically distributed number of cycles. We also derive an upper bound, uniform in $n$, on the expected time until the start of the first cycle, as well as the expected cycle duration, and then combine these bounds and prove   Lemma~\ref{lem:coupling_time_bound}.

\subsection{Higher-Order Bounds}
\label{sec:high_diff}
To prove the higher-order bounds, we first use the Poisson equation to write  $\Delta_1^{2} f_h(x^q)$ in terms of $h(x^q)$, $\E h(X)$, and first-order differences of $f_h(x^q)$. With the help of this expression, we use the dynamics of the JSQ model to relate all the second-order differences  to each other and prove that 
\begin{align}
\abs{\Delta_1^{a_1} \Delta_2^{a_2} f_h(x^q_1,x^q_2,0,\ldots,0)} \leq&\  \delta^2 \sum_{i=1}^{b+1} \E X_i +  C(b,\beta)\delta^{2}(1+x^q_2)^{2} \label{eq:interm_d2}
\end{align}  
for   $\norm{a}_{1} = 2$,  $x^q_1 \leq \delta(n-a_1)$ and $x^q_2 \leq \delta(n-a_2)$, followed by a similar bound for $\abs{(\Delta_1 + \Delta_2) f_h(0,x^q_2,0,\ldots, 0)}$. We then bound $\sum_{i=1}^{b+1} \E X_i$ using the Poisson equation in Section~\ref{sec:moments} and  bound $\abs{\Delta_1^{3} f_h(x^q)}$ and $\abs{(\Delta_1^2 - (\Delta_1+\Delta_2))    f_h(x^q)}$ in Section~\ref{sec:third_bounds}. \blue{In Section~\ref{sec:proof_relate}, we prove two technical lemmas needed to establish \eqref{eq:interm_d2}. We also briefly discuss (see Remark~\ref{rem:prelimgen} there) the advantage of using the prelimit generator approach and working with finite differences of $f_h(x^q)$, as opposed to using the classical generator approach and working with the derivatives of the solution to the Poisson equation for the diffusion.} 

For the following discussion, we assume that $x^q \in S$ with $x^q_3 = 0$.  Recall  from \eqref{eq:gxdef} that
\begin{align}
G_X f_h(x^q) =&\ 1(q_1 < n)  n \lambda  \Delta_1^2 f_h(x^q - \delta e^{(1)} ) +  1(q_1 = n, q_2 < n)  n \lambda ( \Delta_2 + \Delta_1 ) f_h(x^q)  \notag \\
&+  \frac{1}{\delta} (\beta - (x^q_1 + x^q_2))  \Delta_1 f_h(x^q) -  \frac{1}{\delta}x^q_2  \Delta_2 f_h(x^q - \delta e^{(2)}). \label{eq:scaled_gen}
\end{align} 
We rearrange the Poisson equation $G_X f_h(x^q) = \E h(X) - h(x^q)$   to see that when $0 < q_1 < n$, or alternatively $0 < x^q_1 < \delta n$,
\begin{align}
\Delta_1^2  f_h(x^q - \delta e^{(1)}) =&\ \frac{1}{n\lambda}(\E   h(X) -  h(x^q)) - \frac{1}{n\lambda} \frac{1}{\delta}(\beta - (x^q_1+x^q_2)) \Delta_1  f_h(x^q) \notag  \\
&+  \frac{1}{n\lambda}\frac{1}{\delta} x^q_2  \Delta_2  f_h(x^q - \delta e^{(2)}). \label{eq:pre3}
\end{align}
Note that $\E \abs{h(X)} \leq  C\E(X_1 + \cdots + X_{b+1})$ since $h(0) = 0$ and $h \in \mathcal{M}_{disc,2}(C)$.  Together with the bound on $\Delta_i f_h(x^q)$ from \eqref{eq:bound_first_difference}, this implies that 
\begin{align}
 \big|\Delta_1^2  f_h(x^q)\big| \leq&\ \delta^2 C \sum_{i=1}^{b+1} \E X_i + \delta^2 x^q_1 + C(b,\beta)\delta^2(1+x^q_2)^2,\quad  x^q_1 < \delta(n-1),\ x^q_3 = 0. \label{eq:bound_second_difference}
\end{align}
Similarly, if $x^q_1 = 0$, 
\begin{align}
(\Delta_2 + \Delta_1 )  f_h(x^q) =&\ \frac{1}{n\lambda}(\E   h(X) -  h(x^q)) - \frac{1}{n\lambda} \frac{1}{\delta}(\beta - x^q_2) \Delta_1  f_h(x^q) +  \frac{1}{n\lambda}\frac{1}{\delta} x^q_2  \Delta_2  f_h(x^q - \delta e^{(2)}), \label{eq:pre4}
\end{align}
and therefore
\begin{align}
\big| (\Delta_2 + \Delta_1 )  f_h(0,x_2,0,\ldots,0)\big|  \leq&\ \delta^2 C \sum_{i=1}^{b+1} \E X_i + C(b,\beta)\delta^2(1+x^q_2)^2, \quad  x^q_2 < \delta n. \label{eq:bound_second_difference2}
\end{align}
Not all second-order differences can be bounded like this. For example, the equation for $\Delta_2^2 f_h(x^q)$ would involve the third-order difference $\Delta_2 \Delta_1^2 f_h(x^q)$, which we have not  bounded. Instead, the following lemma relates the remaining second-order differences   to $\Delta_1^2  f_h(x^q)$ and $(\Delta_2 + \Delta_1 )  f_h(0,x^q_2,0,\ldots,0)$ using the structure of the JSQ system. The proof is postponed to Section~\ref{sec:proof_relate}.
\begin{lemma}
\label{lem:relate}
Fix $h \in \mathcal{M}_{disc,2}(C)$. Then for any $x^q \in S$ with $x^q_3 = 0$, 
\begin{align*}
\abs{\Delta_{1}^{2} f_h(x^q)} \leq&\ C\delta^2 + \max_{0 \leq y^q_2 \leq x^q_2}\abs{\Delta_{1}^{2} f_h(0,y^q_2,0,\ldots,0)}, \quad \text{ provided } \quad  x^q+2\delta e^{(1)} \in S,\\
\abs{\Delta_{2} \Delta_1  f_h(x^q)} \leq&\ C\delta^2 + \max_{\substack{0 \leq y^q_2 \leq x^q_2 \\ j = 1,2}}  \abs{\Delta_{j}^{2} f_h(0,y^q_2,0,\ldots,0)}, \quad  \text{ provided } \quad   x^q+\delta e^{(1)}+\delta e^{(2)}  \in S,\\
\abs{\Delta_{2}^{2} f_h(x^q)} \leq&\ C\delta^2 + \max_{\substack{0 \leq y^q_2 \leq x^q_2\\ j = 1,2}}  \abs{\Delta_{j}^{2} f_h(0,y^q_2,0,\ldots,0)}, \quad  \text{ provided } \quad   x^q+2\delta e^{(2)} \in S.
%(\Delta_{2}^{+})^{3} \widetilde f_h(q) \leq&\ \delta^3 + \max_{0 \leq k \leq q_2}(\Delta_{2}^{+})^{3} \widetilde f_h(n,k)
\end{align*}
\end{lemma}
We see from Lemma~\ref{lem:relate} that to bound the second-order differences, we only need bounds on  $\abs{\Delta_{1}^{2} f_h(0,x^q_2,0,\ldots,0)}$ and $\abs{\Delta_{2}^{2} f_h(0,x^q_2,0,\ldots,0)}$. The former is bounded in \eqref{eq:bound_second_difference}, and for the latter term, we note that for any $x^q \in S$ with $x^q_1=x^q_3=0$, 
\begin{align}
 \big| \Delta_{2}^{2} f_h(x^q)\big| =&\ \big| \Delta_{2}  f_h(x^q+\delta e^{(2)}) - \Delta_{2} f_h(x^q)\big| \notag  \\
=&\ \big|  \big(\Delta_2  + \Delta_1\big) f_h(x^q+\delta e^{(2)}) -  \Delta_1   f_h(x^q+\delta e^{(2)})- \Delta_{2} f_h(x^q)\big|  \notag \\
=&\ \big|  \big(\Delta_2  + \Delta_1\big) f_h(x^q+\delta e^{(2)}) + (f_h(x^q) - f_h(x^q+\delta e^{(1)}+\delta e^{(2)})) \big|  \notag \\
\leq&\  \delta^2 C\sum_{i=1}^{b+1} \E X_i + C(b,\beta)\delta^2(1+x^q_2)^2 + \big| f_h(0,x^q_2,0,\ldots,0) - f_h(\delta,x^q_2+\delta,0,\ldots,0) \big|, \label{eq:diagpre}
\end{align}
where  the  inequality follows from \eqref{eq:bound_second_difference2}. The following lemma bounds the last term on the right-hand side, implying that $ \big| \Delta_{2}^{2} f_h(0,x^q_2,0,\ldots,0)\big| \leq \delta^2 C\sum_{i=1}^{b+1} \E X_i + C(b,\beta)\delta^2(1+x^q_2)^2$, and, consequently, \eqref{eq:interm_d2}.  It is proved in Section~\ref{sec:proof_relate}. 
\begin{lemma}
\label{lem:diagonals}
For all $n \geq 1$,
\begin{align}
& \big| f_h(0,x^q_2,0,\ldots,0) - f_h(\delta,x^q_2+\delta,0,\ldots,0) \big| \leq C(b,\beta) \delta^2 (1+x^q_2), \quad 0 \leq x^q_2 < \delta n. \label{eq:diag1}
%&\big|- \widetilde f_h(n-1,q_2+2) + \widetilde f_h(n,q_2+1) +\widetilde f_h(n-1,q_2+1) - \widetilde f_h(n,q_2)\big| \leq C \delta^3(1+\delta q_2)^2. \label{eq:diag2}
\end{align}
\end{lemma}  

\subsubsection{Bounding $\sum_{i=1}^{b+1} \E X_i$.}
\label{sec:moments}
The bounds in \eqref{eq:interm_d2} and \eqref{eq:bound_second_difference2} do not yet look like the stated bounds in Proposition~\ref{prop:diffbounds} because the term $\sum_{i=1}^{b+1} \E X_i$ is present. However, we can bound this expectation using the Poisson equation as follows. Recall that $\lambda = 1-\beta/\sqrt{n}$, let  $x(\infty) = \big( \delta(n-\lfloor n\lambda \rfloor), 0, \ldots, 0\big) = \big( \beta+  \delta(n\lambda -\lfloor n\lambda \rfloor), 0, \ldots, 0\big)$, and observe that this point is in $S$. In fact, it is the closest point in $S$, when rounded up, to the fluid equilibrium of the JSQ system, which happens to be $(\beta,0,\ldots, 0)$; cf. \cite{Brav2020}. From \eqref{eq:scaled_gen} we have 
\begin{align*}
G_X f_h(x(\infty)) =&\  n \lambda  \Delta_1^2 f_h(x(\infty) - \delta e^{(1)} ) +(n\lambda -\lfloor n\lambda \rfloor) \Delta_1 f_h(x(\infty)) = \E h(X) - h(x(\infty)).
\end{align*}
Choosing $h(x^q) = \sum_{i=1}^{b+1} x^q_i$ and noting that $h(x(\infty)) = \beta+  \delta(n\lambda -\lfloor n\lambda \rfloor)$ yields
\begin{align}
 n \lambda  \Delta_1^2 f_h(x(\infty)- \delta e^{(1)})   + (n\lambda -\lfloor n\lambda \rfloor) \Delta_1 f_h(x(\infty))   -  \beta-  \delta(n\lambda -\lfloor n\lambda \rfloor) = \sum_{i=1}^{b+1} \E X_i . \label{eq:mompoisson}
\end{align} 
To bound $\sum_{i=1}^{b+1} \E X_i$ we need only bound $ \Delta_1^2 f_h(x(\infty)- \delta e^{(1)})$, because $\abs{\Delta_1 f_h(x(\infty))} \leq \delta C(b,\beta)$  due to \eqref{eq:bound_first_difference}.  Note that we cannot use \eqref{eq:bound_second_difference} for the second-order difference bound because  $\sum_{i=1}^{b+1} \E X_i$ is present on the right-hand side there. Instead, we  exploit the structure of the JSQ model to bound $ \Delta_1^2 f(x(\infty)- \delta e^{(1)})$ as follows.

Define $\tau^{-}(x^q_1) = \inf_{t \geq 0} \{X(t) = (x^q_1-\delta,0, \ldots, 0) |  X(0)=(x^q_1,0, \ldots, 0) \}$, let $(X(t),\widetilde X(t))$ be the scaled version of the  coupling defined in Lemma~\ref{lem:coupling}, and let $V$ be the unit-rate exponentially distributed random variable defined in the same lemma. Fix $x^q = (x^q_1,0,\ldots,0)$ with $x^q_1  \geq 2\delta$, and suppose $X(0) = x^q$ and $\widetilde X(0) = x^q - \delta e^{(1)}$. Consider the evolution of $(X(t),\widetilde X(t))$ for $t \in [0, V \wedge \tau^{-}(x^q_1)]$.   If $V < \tau^{-}(x^q_1)$, the two processes couple and become identical. Otherwise, the joint process is in state $(x^q- \delta e^{(1)},x^q - 2\delta e^{(1)})$. Using the strong Markov property, we conclude that
\begin{align*}
 \Delta_1 f_h(x^q - \delta e^{(1)})=&\  \int_{0}^{\infty} \E_{x^q}\Big[ \Big( h\big(X(t)\big) - h\big(X(t) -\delta e^{(1)}\big)\Big) 1(t \leq (V \wedge \tau^{-}(x^q_1)))\Big] dt \\
&+  \Prob(V \geq \tau^{-}(x^q_1))\Delta_{1} f_h(x^q - 2\delta e^{(1)}).
\end{align*}
Choosing $x^q_1 = x_1(\infty)$, we see that
\begin{align*}
&\Delta_{1}  f_h\big(x(\infty)- \delta e^{(1)}\big) - \Delta_{1}   f_h\big(x(\infty)- 2\delta e^{(1)}\big) \\
=&\ \int_{0}^{\infty} \E_{x(\infty)}\Big[ \Big( h\big(X(t) \big) -  h\big(X(t)- \delta   e^{(1)}\big)\Big) 1\Big(t \leq \big(V \wedge \tau^{-}(x_1(\infty))\big)\Big)\Big] dt \\
&- \Prob\big(V < \tau^{-}(x_1(\infty))\big)\Delta_{1}   f_h\big(x(\infty)- 2\delta e^{(1)}\big).
\end{align*}
Choosing $h(x^q) = \sum_{i=1}^{b+1} x^q_i$ and using $\abs{\Delta_1 f(x(\infty))} \leq \delta C(b,\beta)$, we arrive at
\begin{align}
  \big| \Delta_{1}^2  f_h\big(x(\infty)- \delta e^{(1)}\big) \big|  \leq&\ \delta \E \tau^{-}(x_1(\infty)) +C(b,\beta) \delta  \Prob\big(V < \tau^{-}(x_1(\infty))\big). \label{eq:liphbound}
\end{align}
The quantities involving $\tau^{-}(x_1(\infty))$ are bounded in the following lemma.   
\begin{lemma}
\label{lem:up}
There exists a constant $C(\beta)> 0$ such that for all $n \geq 1$, 
\begin{align}
\E \tau^{-}(x^q_1) \leq C(\beta) \delta, \quad \text{ and } \quad \Prob(V \leq \tau^{-}(x^q_1)) \leq C(\beta) \delta, \quad \text{ for } \quad  x^q_1 \in \{x_1(\infty), \delta, 2\delta\}.  \label{eq:uphitbound}
\end{align} 
\end{lemma}
Lemma~\ref{lem:up} is proved in Appendix~\ref{proof:up}. It implies that  $\big| \Delta_{1}^2  f_h\big(x(\infty)- \delta e^{(1)}\big) \big| \leq C(b,\beta) \delta^2$, and therefore 
\begin{align}
\sum_{i=1}^{b+1} \E X_i \leq C(b,\beta). \label{eq:moment_bound}
\end{align}
Combining \eqref{eq:moment_bound} with \eqref{eq:interm_d2} proves the second-order bounds in Proposition~\ref{prop:diffbounds}. 

Before moving on, let us make a few remarks. The bound in \eqref{eq:moment_bound} implies that the sequence of steady-state distributions $\{X\}_{n=1}^{\infty}$ is tight  and, when combined with process-level convergence of $\{X(t)\}$ to the diffusion $\{Y(t)\}$, tightness can be used to imply convergence of the steady-state distributions   via a limit-interchange argument; for an example of this applied to the JSQ model, see \cite{Brav2020}. \blue{Alternatively, \eqref{eq:moment_bound} can be recast into a result about the convergence rate to the mean-field equilibrium.}

\blue{Let $h(x) = \abs{x_1 + \ldots + x_{b+1}-\beta} - \beta$, noting that $h \in \mathcal{M}_{disc,1}(1)$ and that $h(0) = 0$, and suppose for the sake of exposition that $\lfloor n\lambda \rfloor = n\lambda$. One may check that the bound in \eqref{eq:liphbound} holds even when $h \in \mathcal{M}_{disc,1}(1)$, in which case \eqref{eq:mompoisson} implies that 
\begin{align*}
\E \Big|\sum_{i=1}^{b+1}  X_i - \beta\Big| =&\ \beta +  n \lambda  \Delta_1^2 f_h(x(\infty)- \delta e^{(1)}) \leq  C(b,\beta).
\end{align*}
If we divide both sides by $\sqrt{n}$ to consider the mean-field scaled version of $\sum_{i=1}^{b+1}  X_i$, we get
\begin{align*}
\E \Big|(n-Q_1)/n + \sum_{i=2}^{b+1}  Q_i/n - \beta\Big|  \leq  C(b,\beta)/\sqrt{n}.
\end{align*}
Thus, we recover the  $1/\sqrt{n}$ rate of convergence to the mean field equilibrium that one typically obtains using Stein's method for the mean-field model, like in \cite{Ying2017}. The approach used to show tightness in this section can offer an alternative to the one proposed by \cite{Ying2017}, but the difficulty of implementing our approach is directly related to the difficulty of obtaining the relevant Stein factor bounds. }

As a final remark,  in this section we have shown   that establishing tightness, or rates of convergence to the mean-field equilibrium, is equivalent to  bounding the first- and-second-order differences of $f_h(x^q)$ at a \emph{single} point near the fluid equilibrium of the CTMC. In contrast, establishing rates of convergence to the diffusion requires bounds on the second- and-third-order differences at \emph{all points} in the support of $Y$. 

\subsubsection{Third-Order Bounds.}
\label{sec:third_bounds}
To bound $\Delta_1^3 f_h(x^q)$, we recall \eqref{eq:pre3}, which says that for $0 < x^q_1 < \delta n$ with $x_3 = 0$,
\begin{align*}
\Delta_1^2  f_h(x^q - \delta e^{(1)}) =&\ \frac{1}{n\lambda}(\E   h(X) -  h(x^q)) - \frac{1}{n\lambda} \frac{1}{\delta}(\beta - (x^q_1+x^q_2)) \Delta_1  f_h(x^q) \notag  \\
&+  \frac{1}{n\lambda}\frac{1}{\delta} x^q_2  \Delta_2  f_h(x^q - \delta e^{(2)}).
\end{align*}
%We want to apply $\Delta_1$ to both sides. Note that it $g(x^q) = (\beta - (x^q_1+x^q_2)) \Delta_1  f_h(x^q)$, then 
%\begin{align*}
%\Delta_1 g(x^q) =&\ (\beta - (x^q_1+\delta +x^q_2)) \Delta_1  f_h(x^q+ \delta e^{(1)}) - (\beta - (x^q_1+x^q_2)) \Delta_1  f_h(x^q) \\
%=&\ (\beta - (x^q_1+x^q_2)) \Delta_1^2 f_h(x^q) - \delta \Delta_1  f_h(x^q+ \delta e^{(1)}),
%\end{align*}
%and therefore for
Applying $\Delta_1$ to both sides yields 
\begin{align*}
\Delta_1^3 f_h(x^q- \delta e^{(1)}) =&\ -\frac{1}{n\lambda} \Delta_1 h(x^q)- \frac{1}{n\lambda} \frac{1}{\delta}(\beta - (x^q_1+x^q_2)) \Delta_1^2  f_h(x^q) +\frac{1}{n\lambda}  \Delta_1  f_h(x^q+\delta e^{(1)}) \\
&+ \frac{1}{n\lambda}\frac{1}{\delta} x^q_2  \Delta_1 \Delta_2  f_h(x^q - \delta e^{(2)}), \quad 0 < x^q_1 < \delta (n-1).
\end{align*}
The bounds on the first- and-second-order differences of $f_h(x^q)$, together with the fact that $h \in \mathcal{M}_{disc,2}(C)$, imply that 
\begin{align*}
\abs{\Delta_1^3 f_h(x^q)} \leq C(b,\beta)\delta^3(1+x^q_2)^3, \quad x^q \in S,\ x^q_1 \leq \delta(n-3), 
\end{align*}
which matches the inequality in Proposition~\ref{prop:diffbounds}. The bound on $\abs{(\Delta_1^2 - (\Delta_1+\Delta_2))    f_h(x^q)}$ when $x^q_1 = 0$  is proved identically by subtracting $(\Delta_1+\Delta_2) f_h(x^q)$ in \eqref{eq:pre4} from $\Delta_1^2 f_h(x^q)$ in \eqref{eq:pre3}. This concludes the proof of Proposition~\ref{prop:diffbounds}. 
\finishproof

\subsubsection{Proving  Lemmas \ref{lem:relate} and \ref{lem:diagonals}.}
\label{sec:proof_relate}
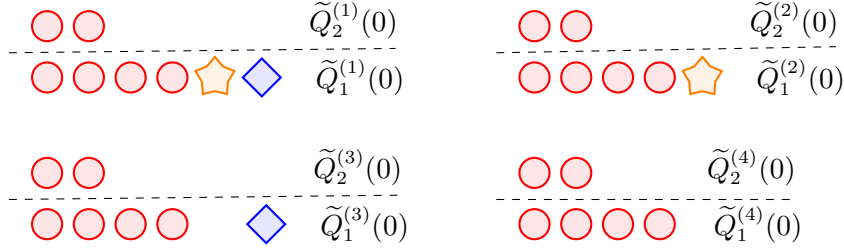
\begin{figure}[h!]
\centering
 \begin{tikzpicture}
  [server/.style={circle, inner sep=1.0mm, minimum width=.8cm,
    draw=green,fill=green!10,thick},
    customer/.style={circle, inner sep=0.5mm, minimum width=.4cm,
    draw=red,fill=red!10,thick},
    blue/.style={diamond,
    draw=blue,fill=blue!10,thick},
    green/.style={rectangle, inner sep=0.5mm, minimum width=.4cm, minimum height=.4cm,
    draw=green,fill=green!10,thick},
    orange/.style={star, draw=orange,fill=orange!10,thick},
  buffer/.style={rectangle, rounded corners=3pt,
    inner sep=0.0mm, minimum width=.9cm,  minimum height=.8cm,
    draw=orange,fill=blue!10,thick},
 vbuffer/.style={rectangle,rounded corners=3pt,
     inner sep=0.0mm,  minimum width=.6cm, minimum height=.8cm,
     draw=orange,fill=blue!10,thick}
]

  \node[customer]  (C1) {};
  \node[customer]  (C2) [right =.05in of C1]   {};
  \node[customer]  (C3) [right =.05in of C2]   {};
  \node[customer]  (C4) [right =.05in of C3]   {};
  \node[orange]  (C5) [right =.05in of C4]   {};
  \node[blue]  (C6) [right =.05in of C5]   {};
  
  \node[customer]  (C11) [above =.1in of C1] {};
  \node[customer]  (C12) [above =.1in of C2] {};
  
  \node[customer]  (B1) [right =1.5in of C5]    {};
  \node[customer]  (B2) [right =.05in of B1]   {};
  \node[customer]  (B3) [right =.05in of B2]   {};
  \node[customer]  (B4) [right =.05in of B3]   {};
  \node[orange]  (B5) [right =.05in of B4]   {};
  
  \node[customer]  (B11) [above =.1in of B1] {};
  \node[customer]  (B12) [above =.1in of B2] {};

  \node[customer]  (A1) [below =.6in of C1]    {};
  \node[customer]  (A2) [right =.05in of A1]   {};
  \node[customer]  (A3) [right =.05in of A2]   {};
  \node[customer]  (A4) [right =.05in of A3]   {};
  \node[blue]  (A6) [right =.3in of A4]   {};
  
  \node[customer]  (A11) [above =.1in of A1] {};
  \node[customer]  (A12) [above =.1in of A2] {};

  \node[customer]  (D1) [below =.6in of B1]    {};
  \node[customer]  (D2) [right =.05in of D1]   {};
  \node[customer]  (D3) [right =.05in of D2]   {};
  \node[customer]  (D4) [right =.05in of D3]   {};
  
  \node[customer]  (D11) [above =.1in of D1] {};
  \node[customer]  (D12) [above =.1in of D2] {};

\draw[dashed]  ($(C1.north)+(-0.2in,.045in)$) edge ($(C6.north)+(0.7in,.045in)$);
  
  \node[left, align=left] (cc1) at ($(C6)+(0.8in,0)$) {$\widetilde Q^{(1)}_1(0)$};
  \node[left, align=left] (cc2) at ($(cc1)+(.25in,0.3in)$) {$\widetilde Q^{(1)}_2(0)$}; 

\draw[dashed]  ($(B1.north)+(-0.2in,.045in)$) edge ($(B5.north)+(0.7in,.045in)$);
  
  \node[left, align=left] (bb1) at ($(B5)+(0.8in,0)$) {$\widetilde Q_1^{(2)}(0)$};
  \node[left, align=left] (bb2) at ($(bb1)+(.25in,0.3in)$) {$\widetilde Q_2^{(2)}(0)$}; 

\draw[dashed]  ($(A1.north)+(-0.2in,.045in)$) edge ($(A6.north)+(0.7in,.045in)$);
  
  \node[left, align=left] (aa1) at ($(A6)+(0.8in,0)$) {$\widetilde Q^{(3)}_1(0)$};
  \node[left, align=left] (aa2) at ($(aa1)+(.25in,0.3in)$) {$\widetilde Q^{(3)}_2(0)$}; 

\draw[dashed]  ($(D1.north)+(-0.2in,.045in)$) edge ($(D4.north)+(0.7in,.045in)$);
  
  \node[left, align=left] (dd1) at ($(D4)+(0.8in,0)$) {$\widetilde Q_1^{(4)}(0)$};
  \node[left, align=left] (dd2) at ($(dd1)+(.25in,0.3in)$) {$\widetilde Q_2^{(4)}(0)$}; 
  
 \end{tikzpicture}
\caption{The initial state of the four systems. The red customers represent those common to all four systems. The diamond and star are the extra customers.}
\label{fig:init_relate1}
\end{figure}

To conclude the section, we prove the auxiliary lemmas from Section~\ref{sec:high_diff}. 
\startproof{Proof of Lemma~\ref{lem:relate}}
Our first task is to bound
\begin{align*}
\Delta_{1}^{2} f_h(x^q) = \int_{0}^{\infty} \Big( \E_{x^q+2\delta e^{(1)}} h(X(t)) - 2 \E_{x^q+\delta e^{(1)}} h(X(t)) + \E_{x^q} h(X(t)) \Big) dt.
\end{align*}
Note that $x^q \in S$ with  $x^q + 2 \delta e^{(1)} \in S$ implies that $q_2 \leq q_1-2$. Working with the unscaled CTMC, we now construct four processes $\{\widetilde Q^{(1)}(t)\},\ldots, \{\widetilde Q^{(4)}(t)\}$ defined on the time interval $[0,\tau_1(n)]$, where
\begin{align}
\tau_1(n) = \inf_{t \geq 0} \big\{\widetilde Q_1^{(1)}(t) = n\big\} = \inf_{t \geq 0} \big\{\widetilde X_1^{(1)}(t) = 0\big\}. \label{eq:tau1def}
\end{align}
We refer to $\{\widetilde Q^{(i)}(t)\}$ as the $i$th process. Process four is a  copy of $\{Q(t)\}$. Numbers two and three are copies of four, but with one extra customer, who is assigned to a server with an empty buffer. The extra customer in two is different from the one in three.  Lastly, process one is a copy of four, but with two extra customers. The extra customers are the same as those in two and three. Figure~\ref{fig:init_relate1} visualizes the initial condition of the  processes.

\noindent Let $\{\widetilde X^{(1)}(t)\},\ldots, \{\widetilde X^{(4)}(t)\}$ be the scaled counterparts of these processes. 
%Assume  $\widetilde Q^{(1)}(0) =q$, and let 
%\begin{align}
%\tau_1(n) = \inf_{t \geq 0} \{\widetilde Q_1^{(1)}(t) = n\} = \inf_{t \geq 0} \{\widetilde X_1^{(1)}(t) = 0\}. \label{eq:tau1def}
%\end{align}
Note that 
\begin{align*}
\Delta_{1}^{2} f_h(x^q) =&\ \int_{0}^{\infty} \Big( \E_{x^q+2\delta e^{(1)}} h(X(t)) - 2 \E_{x^q+\delta e^{(1)}} h(X(t)) + \E_{x^q} h(X(t)) \Big) dt\\
=&\ \int_{0}^{\infty} \E_{\widetilde X^{(1)}(0) = x^q}\Big( \big(h(\widetilde X^{(4)}(t)) -   h(\widetilde X^{(3)}(t)) \big) - \big(  h(\widetilde X^{(2)}(t))   -  h(\widetilde X^{(1)}(t)) \big)   \Big) dt.
\end{align*}
We refer to the different customers according to their shapes in Figure~\ref{fig:init_relate1}.
%To bound the integral above, note that at time $\tau_1(n)$, either both the diamond and star customers are still in the system, or at least one of them has left, in which case the integrand above equals zero. 
Define $\tau_s$ and $\tau_d$ to be the service times of the server with the star  and diamond customer, respectively. Both are exponentially distributed with unit mean. Setting $\tau_{m} = \min \big\{\tau_s,\tau_d,\tau_1(n) \big\}$, we observe that if $\tau_m = \tau_s$, then 
\begin{align*}
&\widetilde X^{(1)}(t) = \widetilde X^{(3)}(t), \quad \widetilde X^{(2)}(t) = \widetilde X^{(4)}(t),  \quad t \geq \tau_m,
\end{align*}
and if $\tau_{m} = \tau_{d}$, then 
\begin{align*}
\widetilde X^{(1)}(t) = \widetilde X^{(2)}(t), \quad \text{ and } \quad   \widetilde X^{(3)}(t) = \widetilde X^{(4)}(t), \quad t \geq \tau_m.
\end{align*}
Therefore,  
\begin{align}
&\int_{0}^{\infty} \E_{\widetilde X^{(1)}(0) = x}\Big( \big(h(\widetilde X^{(4)}(t)) -   h(\widetilde X^{(3)}(t)) \big) - \big(  h(\widetilde X^{(2)}(t))   -  h(\widetilde X^{(1)}(t)) \big)    \Big) dt \notag \\
=&\ \E_{\widetilde X^{(1)}(0) = x}\int_{0}^{\tau_{m}} \Big( \big(h(\widetilde X^{(4)}(t)) -   h(\widetilde X^{(3)}(t)) \big) - \big(  h(\widetilde X^{(2)}(t))   -  h(\widetilde X^{(1)}(t)) \big)     \Big) dt  \notag \\
&+ \Prob_{\widetilde X^{(1)}(0) = x} (\tau_m  = \tau_1(n)) \E_{\widetilde X^{(1)}(0) = x} \Big[ \Delta_{1}^{2} f_h\big(0,\widetilde X^{(1)}_2(\tau_1(n)),0,\ldots, 0\big)  \Big| \tau_m = \tau_1(n)  \Big]. \label{eq:almost}
\end{align}
Since $\widetilde X^{(4)}(t) = \widetilde X^{(3)}(t) + \delta e^{(1)} = \widetilde X^{(2)}(t) + \delta e^{(1)} = \widetilde X^{(1)}(t)+ 2\delta e^{(1)} $ for $0 \leq t \leq \tau_m$, 
\begin{align*}
  \big| \big(h(\widetilde X^{(4)}(t)) -   h(\widetilde X^{(3)}(t)) \big) - \big(  h(\widetilde X^{(2)}(t))   -  h(\widetilde X^{(1)}(t)) \big)   \big| =&\ \big| \Delta_1^2    h(\widetilde X^{(1)}(t))\big| \leq C \delta^2,
\end{align*}
where the last inequality follows from $h \in \mathcal{M}_{disc,2}(C)$.  Combining this with the facts that  $\widetilde X^{(1)}_2(\tau_1(n)) \leq \widetilde X^{(1)}_2(0)$  and $\E_{x} \tau_m  \leq \E \tau_s = 1$, we conclude that the right-hand side of \eqref{eq:almost} is bounded by $C \delta^2  + \big| \Delta_{1}^{2} f_h\big(0,x^q_2,0,\ldots, 0\big) \big|$, which proves the bound on $\abs{\Delta_{1}^{2} f_h(x^q)}$. 

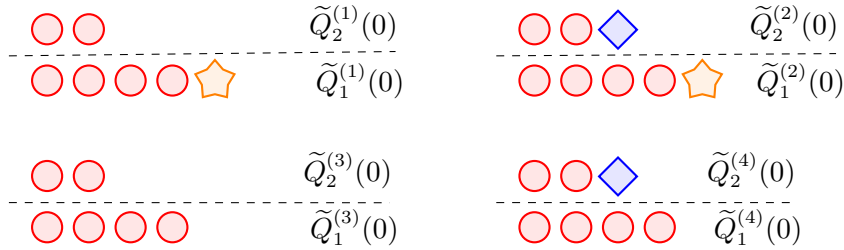
\begin{figure}[h!]
\centering
 \begin{tikzpicture}
  [server/.style={circle, inner sep=1.0mm, minimum width=.8cm,
    draw=green,fill=green!10,thick},
    customer/.style={circle, inner sep=0.5mm, minimum width=.4cm,
    draw=red,fill=red!10,thick},
    blue/.style={diamond,
    draw=blue,fill=blue!10,thick},
    green/.style={rectangle, inner sep=0.5mm, minimum width=.4cm, minimum height=.4cm,
    draw=green,fill=green!10,thick},
    orange/.style={star, draw=orange,fill=orange!10,thick},
  buffer/.style={rectangle, rounded corners=3pt,
    inner sep=0.0mm, minimum width=.9cm,  minimum height=.8cm,
    draw=orange,fill=blue!10,thick},
 vbuffer/.style={rectangle,rounded corners=3pt,
     inner sep=0.0mm,  minimum width=.6cm, minimum height=.8cm,
     draw=orange,fill=blue!10,thick}
]

  \node[customer]  (C1) {};
  \node[customer]  (C2) [right =.05in of C1]   {};
  \node[customer]  (C3) [right =.05in of C2]   {};
  \node[customer]  (C4) [right =.05in of C3]   {};
  \node[orange]  (C5) [right =.05in of C4]   {};
  
  \node[customer]  (C11) [above =.1in of C1] {};
  \node[customer]  (C12) [above =.1in of C2] {};
  
  \node[customer]  (B1) [right =1.5in of C5]    {};
  \node[customer]  (B2) [right =.05in of B1]   {};
  \node[customer]  (B3) [right =.05in of B2]   {};
  \node[customer]  (B4) [right =.05in of B3]   {};
  \node[orange]  (B5) [right =.05in of B4]   {};
  
  \node[customer]  (B11) [above =.1in of B1] {};
  \node[customer]  (B12) [above =.1in of B2] {};
  \node[blue]  (B13) [above =.075in of B3] {};

  \node[customer]  (A1) [below =.6in of C1]    {};
  \node[customer]  (A2) [right =.05in of A1]   {};
  \node[customer]  (A3) [right =.05in of A2]   {};
  \node[customer]  (A4) [right =.05in of A3]   {};
  \node  (A6) [right =.3in of A4]   {};
  
  \node[customer]  (A11) [above =.1in of A1] {};
  \node[customer]  (A12) [above =.1in of A2] {};

  \node[customer]  (D1) [below =.6in of B1]    {};
  \node[customer]  (D2) [right =.05in of D1]   {};
  \node[customer]  (D3) [right =.05in of D2]   {};
  \node[customer]  (D4) [right =.05in of D3]   {};
  
  \node[customer]  (D11) [above =.1in of D1] {};
  \node[customer]  (D12) [above =.1in of D2] {};
  \node[blue]  (D13) [above =.075in of D3] {};

\draw[dashed]  ($(C1.north)+(-0.2in,.045in)$) edge ($(C6.north)+(0.7in,.045in)$);
  
  \node[left, align=left] (cc1) at ($(C6)+(0.8in,0)$) {$\widetilde Q^{(1)}_1(0)$};
  \node[left, align=left] (cc2) at ($(cc1)+(.25in,0.3in)$) {$\widetilde Q^{(1)}_2(0)$}; 

\draw[dashed]  ($(B1.north)+(-0.2in,.045in)$) edge ($(B5.north)+(0.7in,.025in)$);
  
  \node[left, align=left] (bb1) at ($(B5)+(0.8in,0)$) {$\widetilde Q_1^{(2)}(0)$};
  \node[left, align=left] (bb2) at ($(bb1)+(.25in,0.3in)$) {$\widetilde Q_2^{(2)}(0)$}; 

\draw[dashed]  ($(A1.north)+(-0.2in,.045in)$) edge ($(A6.north)+(0.7in,.075in)$);
  
  \node[left, align=left] (aa1) at ($(A6)+(0.8in,0)$) {$\widetilde Q^{(3)}_1(0)$};
  \node[left, align=left] (aa2) at ($(aa1)+(.25in,0.3in)$) {$\widetilde Q^{(3)}_2(0)$}; 

\draw[dashed]  ($(D1.north)+(-0.2in,.045in)$) edge ($(D4.north)+(0.7in,.045in)$);
  
  \node[left, align=left] (dd1) at ($(D4)+(0.8in,0)$) {$\widetilde Q_1^{(4)}(0)$};
  \node[left, align=left] (dd2) at ($(dd1)+(.25in,0.3in)$) {$\widetilde Q_2^{(4)}(0)$}; 
  
 \end{tikzpicture}
\caption{The initial state of the four systems. The red customers represent those common to all four systems.}
\label{fig:init_relate2}
\end{figure}

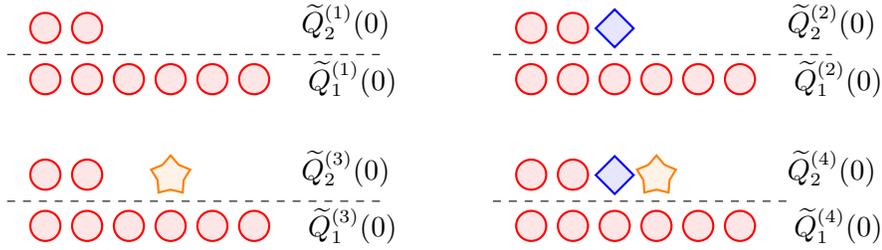
\begin{figure}[h!]
\centering
 \begin{tikzpicture}
  [server/.style={circle, inner sep=1.0mm, minimum width=.8cm,
    draw=green,fill=green!10,thick},
    customer/.style={circle, inner sep=0.5mm, minimum width=.4cm,
    draw=red,fill=red!10,thick},
    blue/.style={diamond,
    draw=blue,fill=blue!10,thick},
    green/.style={rectangle, inner sep=0.5mm, minimum width=.4cm, minimum height=.4cm,
    draw=green,fill=green!10,thick},
    orange/.style={star, draw=orange,fill=orange!10,thick},
  buffer/.style={rectangle, rounded corners=3pt,
    inner sep=0.0mm, minimum width=.9cm,  minimum height=.8cm,
    draw=orange,fill=blue!10,thick},
 vbuffer/.style={rectangle,rounded corners=3pt,
     inner sep=0.0mm,  minimum width=.6cm, minimum height=.8cm,
     draw=orange,fill=blue!10,thick}
]

  \node[customer]  (C1) {};
  \node[customer]  (C2) [right =.05in of C1]   {};
  \node[customer]  (C3) [right =.05in of C2]   {};
  \node[customer]  (C4) [right =.05in of C3]   {};
  \node[customer]  (C5) [right =.05in of C4]   {};
  \node[customer]  (C6) [right =.05in of C5]   {};
  
  \node[customer]  (C11) [above =.1in of C1] {};
  \node[customer]  (C12) [above =.1in of C2] {};
  
  \node[customer]  (B1) [right =1.5in of C5]    {};
  \node[customer]  (B2) [right =.05in of B1]   {};
  \node[customer]  (B3) [right =.05in of B2]   {};
  \node[customer]  (B4) [right =.05in of B3]   {};
  \node[customer]  (B5) [right =.05in of B4]   {};
  \node[customer]  (B6) [right =.05in of B5]   {};
  
  \node[customer]  (B11) [above =.1in of B1] {};
  \node[customer]  (B12) [above =.1in of B2] {};
  \node[blue]  (B13) [above =.075in of B3] {};

  \node[customer]  (A1) [below =.6in of C1]    {};
  \node[customer]  (A2) [right =.05in of A1]   {};
  \node[customer]  (A3) [right =.05in of A2]   {};
  \node[customer]  (A4) [right =.05in of A3]   {};
  \node[customer]  (A5) [right =.05in of A4]   {};
  \node[customer]  (A6) [right =.05in of A5]   {};
  
  \node[customer]  (A11) [above =.1in of A1] {};
  \node[customer]  (A12) [above =.1in of A2] {};
  \node[orange]  (A14) [above =.1in of A4] {};

  \node[customer]  (D1) [below =.6in of B1]    {};
  \node[customer]  (D2) [right =.05in of D1]   {};
  \node[customer]  (D3) [right =.05in of D2]   {};
  \node[customer]  (D4) [right =.05in of D3]   {};
  \node[customer]  (D5) [right =.05in of D4]   {};
  \node[customer]  (D6) [right =.05in of D5]   {};
  
  \node[customer]  (D11) [above =.1in of D1] {};
  \node[customer]  (D12) [above =.1in of D2] {};
  \node[blue]  (D13) [above =.075in of D3] {};
  \node[orange]  (D14) [above =.1in of D4] {};

\draw[dashed]  ($(C1.north)+(-0.2in,.045in)$) edge ($(C6.north)+(0.7in,.045in)$);
  
  \node[left, align=left] (cc1) at ($(C6)+(0.8in,0)$) {$\widetilde Q^{(1)}_1(0)$};
  \node[left, align=left] (cc2) at ($(cc1)+(.25in,0.3in)$) {$\widetilde Q^{(1)}_2(0)$}; 

\draw[dashed]  ($(B1.north)+(-0.2in,.045in)$) edge ($(B5.north)+(0.7in,.045in)$);
  
  \node[left, align=left] (bb1) at ($(B6)+(0.8in,0)$) {$\widetilde Q_1^{(2)}(0)$};
  \node[left, align=left] (bb2) at ($(bb1)+(.25in,0.3in)$) {$\widetilde Q_2^{(2)}(0)$}; 

\draw[dashed]  ($(A1.north)+(-0.2in,.045in)$) edge ($(A6.north)+(0.7in,.045in)$);
  
  \node[left, align=left] (aa1) at ($(A6)+(0.8in,0)$) {$\widetilde Q^{(3)}_1(0)$};
  \node[left, align=left] (aa2) at ($(aa1)+(.25in,0.3in)$) {$\widetilde Q^{(3)}_2(0)$}; 

\draw[dashed]  ($(D1.north)+(-0.2in,.045in)$) edge ($(D4.north)+(0.7in,.045in)$);
  
  \node[left, align=left] (dd1) at ($(D6)+(0.8in,0)$) {$\widetilde Q_1^{(4)}(0)$};
  \node[left, align=left] (dd2) at ($(dd1)+(.25in,0.3in)$) {$\widetilde Q_2^{(4)}(0)$}; 
  
 \end{tikzpicture}
\caption{The  coupling needed to bound $\abs{\Delta_2^2 f_h(x^q)}$.}
\label{fig:init_relate3}
\end{figure}

\begin{figure}[h!]
\centering
\begin{tikzpicture}
 [server/.style={circle, inner sep=1.0mm, minimum width=.8cm,
    draw=green,fill=green!10,thick},
    customer/.style={circle, inner sep=0.5mm, minimum width=.4cm,
    draw=red,fill=red!10,thick},
    blue/.style={diamond,
    draw=blue,fill=blue!10,thick},
    green/.style={rectangle, inner sep=0.5mm, minimum width=.4cm, minimum height=.4cm,
    draw=green,fill=green!10,thick},
    orange/.style={star, draw=orange,fill=orange!10,thick},
  buffer/.style={rectangle, rounded corners=3pt,
    inner sep=0.0mm, minimum width=.9cm,  minimum height=.8cm,
    draw=orange,fill=blue!10,thick},
 vbuffer/.style={rectangle,rounded corners=3pt,
     inner sep=0.0mm,  minimum width=.6cm, minimum height=.8cm,
     draw=orange,fill=blue!10,thick}
]

  \node[customer]  (C1) {};
  \node[customer]  (C2) [right =.05in of C1]   {};
  \node[customer]  (C3) [right =.05in of C2]   {};
  \node[customer]  (C4) [right =.05in of C3]   {};
  \node[customer]  (C5) [right =.05in of C4]   {};
  
  \node[customer]  (C11) [above =.1in of C1]   {};
  \node[customer]  (C21) [above =.1in of C2]   {};
  \node[customer]  (C31) [above =.1in of C3]   {};
  \node[blue]  (C41) [above =.075in of C4]   {};
  
  \node[customer]  (B1) [right =1.5in of C5]    {};
  \node[customer]  (B2) [right =.05in of B1]   {};
  \node[customer]  (B3) [right =.05in of B2]   {};
  \node[customer]  (B4) [right =.05in of B3]   {};
  \node[customer]  (B5) [right =.05in of B4]   {};
  \node[orange]  (B6) [right =.05in of B5]   {};
  
  \node[customer]  (B11) [above =.1in of B1]   {};
  \node[customer]  (B21) [above =.1in of B2]   {};
  \node[customer]  (B31) [above =.1in of B3]   {};

\draw[dashed]  ($(C1.north)+(-0.2in,.045in)$) edge ($(C5.north)+(0.7in,.045in)$);
  
  \node[left, align=left] (X1) at ($(C5)+(0.7in,0)$) {$Q^{(2)}_1(0)$};
  \node[left, align=left] (X1) at ($(C41)+(0.92in,0)$) {$Q^{(2)}_2(0)$}; 

\draw[dashed]  ($(B1.north)+(-0.2in,.045in)$) edge ($(B6.north)+(0.7in,.015in)$);
  
  \node[left, align=left] (X1) at ($(B6)+(0.7in,0)$) {$Q_1^{(1)}(0)$};
  \node[left, align=left] (X1) at ($(B31)+(1.36in,0)$) {$Q_2^{(1)}(0)$}; 
 \end{tikzpicture}
\caption{The initial state of the two systems in an example where $n = 6$. The red circles represent customers common to both systems.}
\label{fig:diag1}
\end{figure}
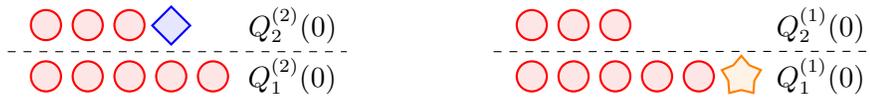

The remaining bounds are proved similarly, starting with $\abs{\Delta_2\Delta_1 f_h(x^q)}$. Fix $x^q \in S$ with $x^q_3=0$, and consider 
\begin{align*}
\Delta_2\Delta_1 f_h(x^q) = \big(f_h(x^q+\delta e^{(1)} + \delta e^{(2)}) -  f_h(x^q+\delta e^{(1)} )\big) - \big(  f_h(x^q+ \delta e^{(2)}) -  f_h(x^q)\big).
\end{align*}
We again construct a coupling $\{\widetilde X^{(1)}(t)\},\ldots, \{\widetilde X^{(4)}(t)\}$ corresponding to the four initial states on the right-hand side above. The initial conditions of the unscaled processes are visualized in Figure~\ref{fig:init_relate2}.
Our construction yields 
\begin{align}
 \Delta_2\Delta_1 f_h(x^q) =&\ \int_{0}^{\infty} \E_{\widetilde X^{(1)}(0) = x^q}\Big( \big(h(\widetilde X^{(4)}(t)) -   h(\widetilde X^{(3)}(t)) \big) - \big(  h(\widetilde X^{(2)}(t))   -  h(\widetilde X^{(1)}(t)) \big)   \Big) dt. \label{eq:integrand2}
\end{align}
Let $\nu_1= \inf_{t \geq 0} \{\widetilde Q_1^{(3)}(t) = n\} $.  We again let $\tau_s$ and $\tau_d$   be the remaining service time of the server with the star and diamond customer, respectively, and set $\tau_m = \min \big\{\tau_s,\tau_d,\nu_1 \big\}$. Just like we argued before, if $\tau_m = \tau_s$, then the integrand in \eqref{eq:integrand2} is zero after $\tau_m$.
If, however, $\tau_m = \tau_d$, then 
\begin{align*}
\widetilde X_2^{(1)}(\tau_m) =&\ \widetilde X_2^{(2)}(\tau_m) = \widetilde X_2^{(3)}(\tau_m) =\widetilde X_2^{(4)}(\tau_m), \\
\widetilde X_1^{(2)}(\tau_m)+ 2\delta =&\ \widetilde X_1^{(1)}(\tau_m)+ \delta = \widetilde X_1^{(4)}(\tau_m)+ \delta =  \widetilde X_1^{(3)}(\tau_m) ,  
\end{align*}
and if $\tau_m = \nu_1$,  then $\widetilde X_1^{(i)}(\tau_m)  = 0$ for $1 \leq i \leq 4$ and
\begin{align*}
\widetilde X_2^{(2)}(\tau_m) =&\ \widetilde X_2^{(1)}(\tau_m) + \delta = \widetilde X_2^{(4)}(\tau_m) + \delta = \widetilde X_2^{(3)}(\tau_m) + 2\delta.
\end{align*}
Therefore,  
\begin{align}
 \Delta_2\Delta_1 f_h(x^q) =&\ \E_{\widetilde X^{(1)}(0) = x^q}\int_{0}^{\tau_{m}} \Big( \big(h(\widetilde X^{(4)}(t)) -   h(\widetilde X^{(3)}(t)) \big) - \big(  h(\widetilde X^{(2)}(t))   -  h(\widetilde X^{(1)}(t)) \big)   \Big)  dt  \notag \\
&+ \Prob_{\widetilde X^{(1)}(0) = x^q} (\tau_m  = \tau_d) \E_{x} \Big[ -\Delta_{1}^{2} f_h\big(\widetilde X^{(2)}(\tau_d)\big)   \Big| \tau_m = \tau_d  \Big]\notag \\
&+ \Prob_{\widetilde X^{(1)}(0) = x^q} (\tau_m  = \nu_1) \E_{x}\Big[ \Delta_{2}^{2} f_h\big(0,\widetilde X^{(3)}_2(\nu_1),0,\ldots, 0\big)  \Big| \tau_m = \nu_1  \Big] \notag \\
\leq&\ C \delta^2+ \big| \Delta_{1}^{2} f_h\big(0,x^q_2,0,\ldots, 0\big) \big|  + \big| \Delta_{2}^{2} f_h\big(0,x^q_2,0,\ldots, 0\big) \big|. \label{eq:almost2}
\end{align}
Figure~\ref{fig:init_relate3} illustrates the coupling needed to bound $\abs{\Delta_2^2 f_h(x^q)}$. The idea of the proof is again to wait until $\tau_1(n)$  and analyze what could happen if one of the servers containing the star or diamond customer completes service before $\tau_1(n)$. We leave the details to the reader.
\finishproof
\begin{remark}
\label{rem:prelimgen}
\blue{Let us say a few words on the advantage of using the prelimit generator comparison approach over the classical generator comparison approach. Lemma~\ref{lem:relate} is proved using a synchronous coupling of four JSQ systems. The four systems are initialized one or two customers apart from one another and due to the discrete state space of the CTMC, all four systems stay one or two customers apart until they couple. Had we used the classical generator comparison approach, we would have  needed to carry out a similar analysis by coupling four copies of the diffusion $\{Y(t)\}$. However, unlike the JSQ coupling, the four diffusions would not maintain their initial spacing relative to each other because $\{Y(t)\}$ takes values in a continuous state space. This would further complicate the analysis as we would now need to keep track of the   positions of the four diffusions relative to each other. }
\end{remark}
\startproof{Proof of Lemma~\ref{lem:diagonals}}   
We want to bound
\begin{align*}
 \abs{f_h(0,x^q_2,0,\ldots,0) - f_h(\delta,x^q_2+\delta,0,\ldots,0)} =&\  \bigg|\int_{0}^{\infty} \big(\E_{(0,x^q_2,0,\ldots,0)} h(X(t)) - \E_{(\delta,x^q_2+\delta,0,\ldots,0)} h(X(t)) \big) dt\bigg|.
\end{align*}
As we are accustomed to doing by now, let us construct a coupling $\{\widetilde Q^{(1)}(t), \widetilde Q^{(2)}(t)\}$  with
\begin{align*}
\widetilde Q^{(1)}(0) = (n,q_2,0,\ldots, 0), \quad \text{ and } \quad \widetilde Q^{(2)}(0) = (n-1,q_2+1,0,\ldots, 0).
\end{align*}
System two has one less idle server and one more customer waiting in a buffer compared to system one, but the total initial customer count is identical across both systems. The initial condition of both systems is visualized in Figure~\ref{fig:diag1}. We assume that the diamond and star customers are independent of each other,  that the systems see identical arrivals, and that the rest of the customers are identical across both systems. 

\noindent Now define $\tau_d$  and $\tau_s$  to be the remaining service times of the server that has the diamond and star customer, respectively;  let $\nu_1 =   \inf_{t \geq 0} \big\{\widetilde Q^{(2)}_1(t) = n\big\}$; and set $\tau_m = \min \big\{\tau_s,\tau_d,\nu_1 \big\}$.  If $\tau_m = \nu_1$ or $\tau_m = \tau_d$, then $\widetilde Q^{(1)}(t) \stackrel{d}{=} \widetilde Q^{(2)}(t)$ for $t \geq \tau_m$. Letting $\{\widetilde X^{(i)}(t)\}$ be the scaled version of $\{\widetilde Q^{(i)}(t)\}$, it follows that 
\begin{align*}
& f_h(0,x^q_2,0,\ldots,0) - f_h(\delta,x^q_2+\delta,0,\ldots,0)  \\
 =&\  \E_{\widetilde X^{(1)}(0) = (0,x^q_2,0,\ldots,0)}\int_{0}^{\tau_{m}} \big( h(\widetilde X^{(1)}(t)) -  h(\widetilde X^{(2)}(t)) \big) dt\\
 &+ \Prob_{\widetilde X^{(1)}(0) =(0,x^q_2,0,\ldots,0)} (\tau_m  = \tau_s) \E_{x^q} \Big[ -\Delta_{2} f_h\big(\widetilde X^{(1)}(\tau_s)\big)  \Big| \tau_m = \tau_s \Big].
\end{align*}
To bound the first term on the right-hand side, note that 
\begin{align*}
\Big| \E_{\widetilde X^{(1)}(0) = (0,x^q_2,0,\ldots,0)}\int_{0}^{\tau_{m}} \big( h(\widetilde X^{(1)}(t)) -  h(\widetilde X^{(2)}(t)) \big) dt \Big| \leq&\   C\delta \E_{\widetilde X^{(2)}(0) = (\delta,x^q_2+\delta,0,\ldots,0)} \nu_1  \leq  C(\beta) \delta^2.
\end{align*}
The first inequality is true because $h \in \mathcal{M}_{disc,2}(C)$, and   the last inequality follows from Lemma~\ref{lem:up} with $x^q_1 = \delta$ there. Furthermore, 
\begin{align*}
&\Prob_{\widetilde X^{(1)}(0) =(0,x^q_2,0,\ldots,0)} (\tau_m  = \tau_s) \Big| \E_{x} \Big[ -\Delta_{2} f_h\big(\widetilde X^{(1)}(\tau_s)\big) \Big| \tau_m = \tau_s  \Big] \Big| \\
& \hspace{3cm} \leq \Prob_{\widetilde X^{(1)}(0) =(0,x^q_2,0,\ldots,0)} ( \tau_s < \nu_1) C(b,\beta) \delta(1 + x^q_2)  \leq  C(b,\beta)\delta^2(1+x^q_2).
\end{align*}
The first inequality follows from the bound on the first-order difference in \eqref{eq:bound_first_difference} together with the fact that $\widetilde X_{2}^{(1)}(t) \leq x^q_2$ for all $t\in [0,\tau_m]$. The second inequality follows by noting that  $\tau_s$ is independent of $\nu_1$  and using Lemma~\ref{lem:up} with    $x^q_1=\delta$, $\tau^{-}(x^q_1)=\nu_1$, and $V = \tau_s$ there.
\finishproof

\section{Conclusion} 
As stated in the introduction, the Stein factor bounds require the bulk of our efforts. Proving the first-order bounds in Section~\ref{sec:first_diff} amounts to considering two coupled JSQ systems, initialized with a difference of one customer, and  bounding the  expected coupling time of this joint chain. We bound the coupling time by considering a sequence of coupling attempts  where the probability of coupling in a single attempt is bounded away from zero uniformly in $n$, and the expected inter-attempt times are also bounded from above, uniformly in $n$. The  coupling time can then be bounded  by a sum of  a geometrically distributed number of random variables representing the inter-attempt durations. This  renewal-like argument applies  more generally to settings where (a) there is a region of the state space where the joint chain is guaranteed to couple provided it spends enough time there and (b) one can control the expected time to reach this region and the probability of coupling in the region before leaving it.

With the first-order Stein factor bounds in hand, the higher-order bounds require less effort. Our proofs of the high-order bounds make heavy use of the transition structure of the JSQ system, and, in particular, that  $Q_2(t), \ldots, Q_{b+1}(t)$ increase  only  at those times when $Q_1(t) = n$.   Readers should not be mislead into thinking that  high-order Stein factor bounds require less effort than first-order bounds for all models. Indeed, in the classical generator comparison approach, high-order bounds require much more effort; e.g., \cite{GorhMack2016, ErdoMackSham2019, Jinetal2021}.

Regarding extending our results, we note that Proposition~\ref{prop:main}, which compares $G_X$ to $G_Y$, can be easily adjusted to hold for other parameter regimes  and  load-balancing policies. The main difficulty would be establishing Stein factor bounds. As mentioned in the introduction, \cite{ZhaoBaneMukh2021} considered the super-Halfin-Whitt regime ($1/2 < \alpha < 1$) and established several hitting-time estimates  similar to the ones we use in the proof of Lemma~\ref{lem:coupling_time_bound} to bound the first-order Stein factors. It may be possible to build on their results and obtain rates of convergence for the super-Halfin-Whitt regime too. 

Furthermore, it seems that the  sub-Halfin-Whitt regime  ($0 < \alpha < 1/2$) should present less of a challenge than our own setting. Recall from the discussion in Section~\ref{sec:first_diff} that coupling of the joint CTMC is guaranteed provided it enters $\Theta_{1}^{Q}$ and spends an exponentially distributed amount of time there  before all servers become busy. Compared to the Halfin-Whitt regime, the rate at which customers arrive in the  sub-Halfin-Whitt regime is much smaller,  so the event that all servers are busy should happen less frequently.   Indeed, \cite{LiuYing2020} showed that the steady-state probability that all servers are busy tends to zero in the sub-Halfin-Whitt regime. Consequently, the Stein factor bounds should be simpler to establish.

\begin{APPENDICES}

\section{Supporting Proofs for Section~\ref{sec:main}}
\label{sec:supportmain}
We first prove Lemma~\ref{lem:ito} and then introduce the operator $A$ in Appendix~\ref{sec:interpolator}. Once $A$ is introduced, we prove Proposition~\ref{prop:main} in Appendix~\ref{sec:prop1proof}.
\startproof{Proof of Lemma~\ref{lem:ito}}
Initialize $Y(0)$ according to $Y$. Since   $\{Y(t)\}$ satisfies \eqref{eq:diffusion2}, for any $f \in C^{2}(\R^{b+1}_+)$ with  $\E \abs{f(Y)} < \infty$,   It\^{o}'s lemma implies that
\begin{align}
0 =&\ \E f(Y(1)) - \E f(Y(0)) \notag \\
 =&\   \E \int_{0}^{1} G_Y f(Y(s)) ds + \E \Big( \int_{0}^{1} \Big(\frac{\partial}{\partial x_1} f(Y(s)) + \frac{\partial}{\partial x_2} f(Y(s)) \Big)1(Y_1(s) = 0)  dU(s) \Big).
\end{align} 
If $\E \abs{G_Y f(Y)} < \infty$, then $\E \int_{0}^{1} G_Y f(Y(s)) ds = \E G_Y f(Y)$ follows from the Fubini-Tonelli theorem.
\finishproof
\subsection{The Interpolator $A$}
\label{sec:interpolator}
The operator $A$ discussed in this section is identical to the one introduced in Appendix~A of \cite{Brav2022}, but we repeat its key properties here as they are needed for the proof of Proposition~\ref{prop:main}. 
Consider a one-dimensional function $f: \delta \Z  \to \R$. We can extend it to $\R$ by defining 
\begin{align*}
A f(x) = \sum_{i=0}^{4} \alpha^{k(x)}_{k(x)+i}(x)f(\delta(k(x)+i)),
\end{align*}
where $k(x) = \lfloor x/\delta\rfloor$ and $\alpha_{k+i}^{k}: \R \to \R$ are weights defined for all $k \in \Z$ and $i = 0, \ldots, 4$. The function $A f(x)$ is  a weighted sum of the five points $f(\delta k(x)), \ldots, f(\delta (k(x)+4))$. We mention the reason for using five points after stating Theorem~\ref{thm:interpolant_def}. Note that if $f(x)$ is  defined only on a subset of $\delta \Z$, then $A f(x)$ can still be defined, provided that $f(\delta k(x)), \ldots, f(\delta (k(x)+4))$ are defined.   \cite{Brav2022} described how to choose these weights to make $Af(x)$  coincide with $f(\cdot)$ on grid points, and also to make it a differentiable function whose  derivatives  behave like the corresponding finite differences of  $f(\cdot)$.   The idea can be applied to multidimensional grid-valued functions  as well.

The following result is Theorem 2 of \cite{Brav2022}. We use this as an interface that contains the important properties of $A$ without delving into the low-level details behind its construction.  
\begin{theorem}
\label{thm:interpolant_def} 
Given a convex set $K \subset \R^{d}$,  define 
\begin{align*}
K_{4} = \{x \in K \cap \delta \Z^{d} : \delta(k(x)+ i) \in K \cap \delta \Z^{d} \text{ for all } 0 \leq i \leq 4e\},
\end{align*} 
let $\text{Conv}(K_4)$ be the convex hull of $K_4$, and, for $x \in \R^{d}$, define $k(x)$ by  $k_{j}(x) = \lfloor x_j/\delta\rfloor$. There exist weights $\big\{\alpha_{k+i}^{k} : \R \to \R,\  k \in \Z,\ i = 0,1,2,3,4\big\}$  such that for any  $f: K \cap  \delta \Z^{d} \to \R$, the function
\begin{align}
A f(x) =&\  \sum_{i_d = 0}^{4} \alpha_{k_d(x)+i_d}^{k_d(x)}(x_d)\cdots \sum_{i_1 = 0}^{4} \alpha_{k_1(x)+i_1}^{k_1(x)}(x_1) f(\delta(k(x)+i)) \notag \\
=&\ \sum_{i_1, \ldots, i_d = 0}^{4} \bigg(\prod_{j=1}^{d} \alpha_{k_j(x) +i_j}^{k_j(x)   }(x_j)\bigg) f(\delta (k(x)+i)) , \quad x \in \text{Conv}(K_4)  \label{eq:af2}
\end{align}
 satisfies $A f(x) \in C^{3}(\text{Conv}(K_4))$, where $i = (i_1, \ldots, i_d)$ in \eqref{eq:af2}. Additionally, $Af(x)$ is infinitely differentiable almost everywhere on $\text{Conv}(K_4)$,
\begin{align}
A f(\delta k) = f(\delta k), \quad \delta k \in K_{4}, \label{eq:interpolates2}
\end{align} 
and  there exists a constant $C(d) > 0$ independent of $f(\cdot)$, $x$, and $\delta$, such that  
\begin{align}
\bigg|  \frac{\partial^{a}}{\partial x^{a}} Af(x) \bigg|   \leq&\  C(d) \delta^{-\norm{a}_{1}}  \max_{\substack{ 0 \leq i_j \leq 4-a_j \\ j = 1,\ldots, d}} \abs{\Delta_{1}^{a_1}\ldots \Delta_{d}^{a_d} f(\delta (k(x)+i))}, \quad x \in \text{Conv}(K_4), \label{eq:multibound2}
\end{align}
for $0 \leq \norm{a}_{1} \leq 3$, and \eqref{eq:multibound2} also holds when $\norm{a}_{1} = 4$  for almost all $x \in  \text{Conv}(K_4)$. Additionally, the weights  $\big\{\alpha_{k+i}^{k} : \R \to \R,\  k \in \Z,\ i = 0,1,2,3,4\big\}$
are degree-$7$ polynomials in $(x-\delta k)/ \delta$ whose coefficients do not depend on $k$ or $\delta$. They satisfy  
\begin{align}
&\alpha_{k}^{k}(\delta k) = 1, \quad \text{ and } \quad \alpha_{k+i}^{k} (\delta k) = 0, \quad &k \in \Z,\ i = 1,2,3,4, \label{eq:alphas_interpolate} \\
&\sum_{i=0}^{4} \alpha^{k}_{k+i}(x) = 1, \quad &k \in \Z,\ x \in \R, \label{eq:weights_sum_one}
\end{align}
and also the following translational invariance property:
\begin{align}
\alpha^{k+j}_{k+j+i}(x+ \delta j)  = \alpha^{k}_{k+i}(x),\quad i,j,k \in \Z, \ x \in \R. \label{eq:weights}
\end{align}
\end{theorem}
\begin{remark}
The bound in \eqref{eq:multibound2} holds almost everywhere when $\norm{a}_{1} = 4$. This bound is the reason we need to use $f(\delta k(x))$ and the four points to the right of it (in each dimension). By using more (fewer) points, one can alter the theorem so that \eqref{eq:multibound2} holds for larger (smaller) values of $\norm{a}_{1}$. It is worth noting that to prove the results in this paper, we do not go beyond $\norm{a}_{1} = 3$.
\end{remark}

Going forward, we let $A$ be the operator described in Theorem~\ref{thm:interpolant_def}. 
Since $Af$ coincides with $f$ on the grid, we refer to $A$ as an interpolator. For the interested reader, $A$ is a degree-7 polynomial spline. From \eqref{eq:interpolates2} we see that $A$ is a linear operator, and  \eqref{eq:weights_sum_one} implies that $A$ applied to a constant simply equals that constant. Before we can prove Proposition~\ref{prop:main}, we require one more lemma. 

\begin{lemma} \label{lem:auxa}
In the setting of Theorem~\ref{thm:interpolant_def}, for any $k \in K_4$ and $1 \leq j \leq d$, 
\begin{align}
\frac{\partial }{\partial x_j} Af(x)\bigg|_{x = \delta k} =  \delta^{-1}\Big(\Delta_j  - \frac{1}{2}\Delta_j^2 + \frac{1}{3}\Delta_j^3 \Big) f(\delta k). \label{eq:derivexplicit}
\end{align}
Furthermore, there exists some $\epsilon: \text{Conv}(K_4) \to \R$ satisfying 
\begin{align*}
\abs{\epsilon(x)}   \leq   C(d) \delta^{-1}  \max_{\substack{ 0 \leq i \leq 4e \\ \norm{a}_{1} = 2 }} \abs{\Delta_{1}^{a_1}\ldots \Delta_{d}^{a_d} f(\delta (k(x)+i))}
\end{align*}
such that for any $x \in \text{Conv}(K_4)$, 
\begin{align*}
\frac{\partial }{\partial x_j} Af(x)   =&\ \delta^{-1}  \Delta_j  f(\delta k(x)) + \epsilon(x).
\end{align*}
%\begin{align*}
%\frac{\partial }{\partial x_1} Af(x) + \frac{\partial}{\partial x_2} Af(x)  = \delta^{-1} (\Delta_1 + \Delta_2) f(\delta k(x)) + \epsilon(x),
%\end{align*} 
\end{lemma}
\startproof{ Proof of Lemma~\ref{lem:auxa}}
 The proof is identical for all indices,  so we assume that $j = 1$. Fix $\delta k \in K_4$ and let  $g(x_1) = Af(x_1, \delta k_2, \ldots, \delta k_{d}  )$ be a function in $x_1$ only.  
The form of $A f(x)$ in \eqref{eq:af2}, together with \eqref{eq:alphas_interpolate}, implies that 
\begin{align*}
\frac{\partial }{\partial x_1} Af(x)\bigg|_{x = \delta k}  = g'(\delta k_1).
\end{align*}
It follows that $g'(\delta k_1) = P_{k_1}'(\delta k_1)$, where $P_{k_1}(x)$ is a polynomial defined in (A.1) of   \cite{Brav2022}. Furthermore, (A.1) implies that
\begin{align*}
P_{k_1}'(\delta k_1) = \delta^{-1}\Big(\Delta_1  - \frac{1}{2}\Delta_1^2 + \frac{1}{3}\Delta_1^3 \Big)  = g(\delta k_1) =\delta^{-1}\Big(\Delta_1  - \frac{1}{2}\Delta_1^2 + \frac{1}{3}\Delta_1^3 \Big) f(\delta k)  ,
\end{align*}
from which \eqref{eq:derivexplicit} follows.  To prove the second claim of the lemma, we write
\begin{align*}
\frac{\partial }{\partial x_j} Af(x) = \delta^{-1}\Big(\Delta_j  - \frac{1}{2}\Delta_j^2 + \frac{1}{3}\Delta_j^3 \Big) f(\delta k(x)) + \frac{\partial }{\partial x_j} Af(x) - \frac{\partial }{\partial x_j} Af(x)\bigg|_{x = \delta k(x)}.
\end{align*}
Now $\abs{\Delta_j^2   f(\delta k(x))} \leq    \max  \big\{ \abs{\Delta_{1}^{a_1}\ldots \Delta_{d}^{a_d} f(\delta (k(x)+i))} :\ 0 \leq i \leq 4e,\ \norm{a}_{1} = 2 \big\}$,
\begin{align*}
\abs{\Delta_j^3   f(\delta k(x))} = \abs{\Delta_j^2   f(\delta (k(x)+ e^{(j)})) - \Delta_j^2   f(\delta k(x))} \leq \max_{\substack{ 0 \leq i \leq 4e \\ \norm{a}_{1} = 2 }} \abs{\Delta_{1}^{a_1}\ldots \Delta_{d}^{a_d} f(\delta (k(x)+i))},
\end{align*} 
and 
\begin{align*}
&\abs{\frac{\partial }{\partial x_j} Af(x) - \frac{\partial }{\partial x_j} Af(x)\bigg|_{x = \delta k(x)}} \\
\leq&\  \sum_{j'=1}^{d} \abs{x_{j'} - \delta k_{j'}(x)} \abs{\frac{\partial^2 }{\partial x_j \partial x_{j'}} A f(\xi)} \leq C(d)  \delta^{-1}\max_{\substack{ 0 \leq i \leq 4e \\ \norm{a}_{1} = 2 }} \abs{\Delta_{1}^{a_1}\ldots \Delta_{d}^{a_d} f(\delta (k(x)+i))},
\end{align*}
where $\xi$ is some point between $\delta k(x)$ and $x$. 
The last inequality follows from \eqref{eq:multibound2} and the fact that $\abs{x_j - \delta k_j(x)} \leq \delta$.
\finishproof 

Note that some of the bounds in Theorem~\ref{thm:interpolant_def} and Lemma~\ref{lem:auxa} have a constant $C(d)$ depending on the dimension $d$ of the function; e.g., \eqref{eq:multibound2}. In the JSQ model $d = b+1$, but  when proving Proposition~\ref{prop:main} in the next section  we can   assume that $d = 2$ because of the following.  Given  a function $f: \delta \N^{b+1} \to \R$, we can use \eqref{eq:af2} and \eqref{eq:alphas_interpolate} of Theorem~\ref{thm:interpolant_def}, and the fact that $Y_i = 0$ for $i > 2$,  to see that
\begin{align*}
 A f(Y)  =&\  \sum_{i_{b+1} = 0}^{4} \alpha_{i_{b+1}}^{0}(0)\cdots \sum_{i_{3} = 0}^{4} \alpha_{i_{3}}^{0}(0)\sum_{i_2 = 0}^{4} \alpha_{k_2(Y)+i_2}^{k_2(Y)}(Y_2) \sum_{i_1 = 0}^{4} \alpha_{k_1(Y)+i_1}^{k_1(Y)}(Y_1) f(\delta(k(Y)+i)) \\
=&\  \sum_{i_2 = 0}^{4} \alpha_{k_2(Y)+i_2}^{k_2(Y)}(Y_2) \sum_{i_1 = 0}^{4} \alpha_{k_1(Y)+i_1}^{k_1(Y)}(Y_1) f\big( \delta(k_1(Y)+i_1), \delta(k_2(Y)+i_2), 0, \ldots, 0 \big).
\end{align*}
Since $k_j(Y)$ depends only on $Y_j$, we see that $A f(Y)$ is actually a bivariate function. In Appendix~\ref{sec:prop1proof},  we treat any function of the form $A f(Y)$ as a function of   two variables. 

\subsection{Proving Proposition~\ref{prop:main}}
\label{sec:prop1proof} 
Fix $h \in \mathcal{M}_{disc,2}(C)$. We recall from \eqref{eq:gxdef} that  for $x^q \in S$,
\begin{align*}
 G_{X} f(x^q)  =&\ -1(q_1 < n) n \lambda  \Delta_1 f(x^q- \delta e^{(1)})  + n \lambda \sum_{j=1}^{b} 1(q_1=\ldots = q_j = n, q_{j+1} < n)  \Delta_{j+1} f(x^q) \notag \\
&+ (q_1 - q_{2})  \Delta_1 f(x^q) - \sum_{j=2}^{b} (q_j - q_{j+1})   \Delta_j f(x^q-\delta e^{(j)}) - q_{b+1}   \Delta_{b+1} f(x^q-\delta e^{(b+1)})
\end{align*}
and $f_h(x^q)$ is the unique solution to the Poisson equation
\begin{align}
G_X f_h(x^q) = \E h(X) - h(x^q), \quad x^q \in S  \label{eq:discrete_poissonapp}
\end{align}
with $f_h(0) = 0$. Also recall that we extended $f_h(x^q)$ by setting $f_h(x^q) = 0$ for $x^q \in \delta \N^{b+1}\setminus S$, and defined 
\begin{align*}
B =&\  \{ (x_1,x_2, 0, \ldots, 0) \in \R^{b+1}_+ :    x_2  + x_1 \leq \delta (n/2 - 8) = (n/2-8)/\sqrt{n} \} \quad \text{ and } \\
I =&\ \big\{i = (i_1,i_2,0,\ldots,0) \in \N^{b+1} :\ 0 \leq i_1,i_2 \leq 4 \big\}.
\end{align*}
We first argue that  $\E \abs{Ah(Y)} < \infty$,  $\E \abs{Af_h(Y)} < \infty$, and $\E \abs{G_Y Af_h(Y)} < \infty$, which together imply that \eqref{eq:prefinalerror} holds. The latter two statements follow immediately from the fact that $f_h(x^q)$, and therefore $A f_h(x)$, have compact support. Since $h \in \mathcal{M}_{disc,2}(C)$, inequality \eqref{eq:multibound2} of Theorem~\ref{thm:interpolant_def} implies that $Ah(Y)$ is Lipschitz and therefore,  $\E \abs{Ah(Y)}< \infty$ due to Lemma~\ref{lem:ymoments}, which states that the moments of $Y_i$ are finite.

Next we argue that $A G_X f_h (x) = \E h(X) - A h(x)$ for all $x \in B$. Given the Poisson equation  \eqref{eq:discrete_poissonapp} and the definition of $A$ in \eqref{eq:af2} of Theorem~\ref{thm:interpolant_def}, it suffices to show that $\delta (k(x)+i) \in S$ for all $i \in I$.    From the definition of $S_Q$  in \eqref{eq:sq} we know that any point $q \in S_Q$ satisfies $0 \leq q_2 \leq q_1 \leq n$. The corresponding points $x^q \in S$ satisfy $x^q_1 \geq 0$, $x^q_2 \geq 0$, and  $x^q_1 + x^q_2 = \delta(n-q_1) + \delta q_2 \leq \delta n$. The latter inequality says that the combined number of idle servers and servers with at least one person waiting in the buffer cannot exceed $n$. Now, provided that $n > 16$,   any point $\delta k$ in
\begin{align*}
  B \cap \delta \N^{b+1} =  \{ (x_1,x_2, 0, \ldots, 0) \in \R^{b+1}_+ :    x_2  + x_1 \leq \delta (n/2 - 8) \}\cap \delta \N^{b+1}, 
\end{align*}
  must satisfy  $\delta(k+i) \in S$ for all $i \in I$  because  $\delta(k_1 + i_1) + \delta(k_2 + i_2) \leq \delta n/2$.  Finally, recall that  
\begin{align*} 
&\varepsilon_1(Y)  =   \big(A G_X f_h(Y) -  G_Y A f_h(Y)\big) 1(Y \in B), \\
&\varepsilon_2(Y)  = \big(\E h(X) - A h(Y)-G_Y Af_h(Y)\big) 1(Y \not \in B), \\
&\varepsilon_3(Y)  =   \Big(\frac{\partial}{\partial x_1} Af_h(Y) + \frac{\partial}{\partial x_2} Af_h(Y)\Big) 1(Y \in B), \text{ and }   \\
&\varepsilon_4(Y)  =   \Big(\frac{\partial}{\partial x_1} Af_h(Y) + \frac{\partial}{\partial x_2} Af_h(Y)\Big) 1(Y \not \in B).
\end{align*} 
We bound $\varepsilon_2(Y)$, $\varepsilon_3(Y)$, and $\varepsilon_4(Y)$ in Appendix~\ref{sec:eps24}  and   bound $\varepsilon_1(Y)$ in Appendix~\ref{sec:eps1}.

\subsubsection{Bounding $\varepsilon_2(Y)$ through $\varepsilon_4(Y)$.} \label{sec:eps24}
We begin with the bound on  
\begin{align*}
 \abs{\varepsilon_2(Y)} \leq&\   \abs{A h(Y)} 1(Y \not \in B)  +   1(Y \not \in B)  \E \abs{h(X)} + \abs{G_Y A f_h(Y)} 1(Y \not \in B).
\end{align*}
The facts that $A h(Y)$ is Lipschitz, that $Ah(0) = h(0) = 0$, and that $h \in \mathcal{M}_{disc,2}(C)$ imply that
\begin{align*}
   \abs{ A h(Y) 1(Y \not \in B)} \leq&\ C (Y_1+Y_2)1(Y \not \in B) \quad \text{ and } \\ 
1(Y \not \in B)\E \abs{h(X)} \leq&\  1(Y \not \in B)C\E(X_1 + \cdots + X_{b+1}) \leq 1(Y \not \in B)C(b,\beta),
\end{align*}
where the last inequality follows from inequality \eqref{eq:moment_bound}. 
% Next, since $h \in \mathcal{M}_{disc,2}(C)$, $h(0) = 0$, and $X_i \geq 0$, it follows that $\abs{h(X)} \leq C(X_1 + \cdots + X_{b+1})$. It was shown in \cite{Brav2020} that $\E (X_1 + \cdots + X_{b+1}) \leq C(\beta)$, but we intentionally avoid relying on this bound  in favor of an argument free of any of the machinery of \cite{Brav2020}. We claim that if $\tilde h: \delta \N^{b+1} \to \R$ is defined by $\tilde h(\delta k) = k_1 + \cdots + k_{b+1}$, then   $A \tilde h(x) = x_1 + \cdots +x_{b+1}$. To see this, note that $\Delta_{i}\Delta_{j} \tilde h(\delta k) = 0$ for any $i,j$, and so Theorem~\ref{thm:interpolant_def} implies that all second order partial derivatives of $A \tilde h(x)$ are zero. Since $A \tilde h \in C^{3}(\R^{b+1}_+)$, this implies that $A \tilde h(x)$ is a linear function, and the only linear function that coincides with $\tilde h(\delta k)$ on the grid is $A \tilde h(x) = x_1 + \cdots +x_{b+1}$. It follows that
%\begin{align*}
% 1(Y \not \in B)  \E \abs{h(X)}   \leq&\   1(Y \not \in B)C\E A \tilde h(Y)  + 1(Y \not \in B) C \big(\E \tilde h(X) - \E A \tilde h(Y)\big) \\
%\leq&\  1(Y \not \in B)C \E (Y_1+Y_2) + 1(Y \not \in B)C  \sup_{h \in \mathcal{M}_{disc,2}(C)} \abs{\E h(X) - \E Ah(Y)}\\
%\leq&\    1(Y \not \in B)C \Big(1 +   \sup_{h \in \mathcal{M}_{disc,2}(C)} \abs{\E h(X) - \E Ah(Y)}\Big).
%  \end{align*} 
To bound the remaining term, we recall  \eqref{eq:multibound2} of Theorem~\ref{thm:interpolant_def}, which says that 
\begin{align}
\bigg|  \frac{\partial^{a}}{\partial x^{a}} Af(Y) \bigg| \leq C \delta^{-\norm{a}_{1}} \max_{ \substack{i \in I \\ 0 \leq i_j \leq 4-a_j  }} \abs{\Delta_{1}^{a_1} \Delta_{2}^{a_2} f(\delta (k(Y)+i))} \leq C \delta^{-\norm{a}_{1}}  \max_{\substack{ i \in I }} \abs{ f(\delta (k(Y)+i))} \label{eq:supfbound}
\end{align}
for $1 \leq \norm{a}_{1} \leq 3$. Combined with this bound, the definition of   $G_Y$ in \eqref{eq:gy} implies that 
\begin{align*}
\abs{G_Y A f_h(Y)}  =&\  \Big|\big(\beta - (Y_1 + Y_2) \big)\frac{\partial}{ \partial x_1} A f_h(Y) - Y_2 \frac{\partial}{\partial x_2} A f_h (Y) + \frac{\partial^2}{\partial x_1^2} A f_h (Y)\Big|  \\
\leq&\ C(\beta) \delta^{-2}     (1 + Y_1+Y_2)     \max_{\substack{ i \in I }} \abs{ f(\delta (k(Y)+i))}.
\end{align*} 
Combining the bounds on the three terms yields the bound on $\varepsilon_2(Y)$.   Lemma~\ref{lem:auxa} implies the bound on $\varepsilon_3(Y)$  and \eqref{eq:supfbound} implies the bound on $\varepsilon_4(Y)$.

\subsubsection{Bounding $\varepsilon_1(Y)$.} \label{sec:eps1}
Bounding $\varepsilon_1(Y)$ requires more  effort. The first thing to note is that the weighted sum representation of  $A G_{X} f_h(Y)$ is difficult to work with. Our first task is therefore to write it in a form that is more amenable to analysis. To this end, we extend the domain of $f_h(x^q)$ to allow either the first or second coordinate to take the value $-\delta$ by defining
\begin{align}
&\widehat f_h(x^q) = f_h(x^q), \qquad &x^q \in \delta \N^{b+1}, \notag \\
&\widehat f_h(-\delta, x^q_2,  \ldots,  x^q_{b+1}) = f_h(0,x^q_2 + \delta, x^q_3, \ldots,  x^q_{b+1}), \qquad &(0,x^q_2 , \ldots,  x^q_{b+1}) \in \delta \N^{b+1}, \notag \\
&\widehat f_h(x^q_1, -\delta, x^q_3, \ldots,  x^q_{b+1}) = (1-\Delta_2)f_h(x^q_1,0, x^q_3, \ldots,  x^q_{b+1}), \qquad &(x^q_1, 0 , x^q_3, \ldots,  x^q_{b+1}) \in \delta \N^{b+1} . \label{eq:fhat}
\end{align}  
The form of $\widehat f_h(x^q)$ is tied  to the transition structure of the JSQ model, and specifically to the ``reflection'' that occurs near the boundaries $\{x^q_1 = 0\}$ and $\{x^q_2 = 0\}$. Furthermore, the definition of $A$ in Theorem~\ref{thm:interpolant_def} implies that $A f_h(x) = A \widehat f_h(x)$ for $x \in \R^{b+1}_+$ because $\widehat f_h  = f_h $ on $\delta \N^{b+1}$.   Having defined $\widehat f_h(x^q)$, we present the following lemma, which is proved in Appendix~\ref{sec:actualinterchangeproof}.
\begin{lemma}
\label{lem:actualinterchange}
For any $x^q \in B \cap \delta \N^{b+1}$, 
\begin{align}
G_X f_h(x^q)  =&\ n \lambda \big(\widehat f_h(x^q-\delta e^{(1)}) - \widehat f_h(x^q)\big) + (n-(x^q_1 +  x^q_2)/\delta) \big( \widehat f_h(x^q + \delta e^{(1)}) - \widehat f_h(x^q) \big) \notag \\
&+ \frac{1}{\delta} x^q_2 \big( \widehat f_h(x^q - \delta e^{(2)}) - \widehat f_h(x^q) \big). \label{eq:gxfhat}
\end{align}
Consequently, for any $x \in B$, 
\begin{align}
A G_X f_h(x) =&\ n\lambda \big( A \widehat f_h(x-\delta e^{(1)}) - A \widehat f_h(x) \big) + (n-(x_1 + x_2)/\delta) \big( A \widehat f_h(x + \delta e^{(1)}) - A\widehat f_h(x) \big)  \notag \\
&+ \frac{1}{\delta} x_2 \big( A \widehat f_h(x - \delta e^{(2)}) - A \widehat f_h(x) \big) + \varepsilon_5(x), \label{eq:agfinal}
\end{align}
where 
\begin{align}
\varepsilon_5(x) =&\ \sum_{i_2 = 0}^{4} \alpha_{k_2(x)+i_2}^{k_2(x)}(x_2) \sum_{i_1 = 0}^{4} \alpha_{k_1(x)+i_1}^{k_1(x)}(x_1) \frac{1}{\delta} \Big(\delta (k_2(x)+i_2) - x_2\Big) \notag \\
& \hspace{2cm} \times \Big(-\Delta_2 \widehat f_h\big(\delta (k(x)+i-e^{(2)})\big) + \Delta_2 \widehat f_h\big(\delta (k(x)-e^{(2)})   \Big)  \notag \\
&+ \sum_{i_2 = 0}^{4} \alpha_{k_2(x)+i_2}^{k_2(x)}(x_2) \sum_{i_1 = 0}^{4} \alpha_{k_1(x)+i_1}^{k_1(x)}(x_1) \frac{1}{\delta} \Big(-\delta (k_1(x)+i_1 + k_2(x)+i_2) +x_1+ x_2\Big) \notag \\
& \hspace{2cm} \times \Big(\Delta_1 \widehat f_h\big(\delta (k(x)+i )\big) - \Delta_1 \widehat f_h\big(\delta k(x)\big) \Big). \label{eq:e5}
\end{align}
\end{lemma}
\noindent We now bound  $\varepsilon_1(Y)$ using Lemma~\ref{lem:actualinterchange}. Applying Taylor expansion to \eqref{eq:agfinal}, we have
\begin{align*}
A G_X f_h(Y) =&\ n\lambda \Big( -\delta \frac{\partial }{\partial x_1} A \widehat f_h(Y) + \frac{1}{2}\delta^2 \frac{\partial^2 }{\partial x_1^2} A \widehat f_h(Y) -\frac{1}{6}\delta^3 \frac{\partial^3 }{\partial x_1^3} A \widehat f_h(\xi^1)  \Big) \\
&+ (n-(Y_1+Y_2)/\delta)\Big( \delta \frac{\partial }{\partial x_1} A \widehat f_h(Y) + \frac{1}{2}\delta^2 \frac{\partial^2 }{\partial x_1^2} A \widehat f_h(Y) + \frac{1}{6}\delta^3 \frac{\partial^3 }{\partial x_1^3} A \widehat f_h(\xi^2) \Big)  \notag \\
&+ \frac{1}{\delta} Y_2 \Big(-\delta \frac{\partial }{\partial x_2} A \widehat f_h(Y) + \frac{1}{2}\delta^2 \frac{\partial^2 }{\partial x_2^2} A \widehat f_h(\xi^3) \Big) + \varepsilon_5 (Y),
\end{align*}
where $\xi^1, \xi^2$, and $\xi^3$ are points strictly between $Y-\delta e^{(1)}$ and $Y$, $Y$ and $Y + \delta e^{(1)}$, and $Y- \delta e^{(2)}$ and $Y$, respectively. Recall that $\delta^2 = 1/n$, $\delta (n-n\lambda) = \beta$, and   $G_Y$ from \eqref{eq:gy}, which imply that
\begin{align*}
 A G_X f_h(Y)- G_Y A \widehat f_h(Y) =&\ -  \frac{1}{6}\delta\lambda \frac{\partial^3 }{\partial x_1^3} A \widehat f_h(\xi^1) + \frac{1}{6}\delta(1-\delta(Y_1+Y_2) ) \frac{\partial^3 }{\partial x_1^3} A \widehat f_h(\xi^2) + \delta Y_2  \frac{\partial^2 }{\partial x_2^2} A \widehat f_h(\xi^3) + \varepsilon_5(Y).
\end{align*}
Note that $ A \widehat f_h(Y)  =  A  f_h(Y) $  because $Y \geq 0$,  so $G_Y A \widehat f_h(Y)  = G_Y A  f_h(Y)$. We now prove the following four bounds, which together imply the bound on $\varepsilon_1(Y)$:
\begin{align}
&\abs{ \frac{1}{6}\delta(1-\delta(Y_1+Y_2) ) \frac{\partial^3 }{\partial x_1^3} A \widehat f_h(\xi^2)} \leq C \delta^{-2} \max_{\substack{i \in I }} \abs{\Delta_{1}^{3}    f_h\big(\delta (k(Y) +i )\big)}, \label{eq:blist1}\\
&\abs{\frac{1}{6}\delta\lambda\frac{\partial^3 }{\partial x_1^3} A \widehat f_h(\xi^1)} \leq  C \delta^{-2} \max_{\substack{i \in I }} \abs{\Delta_{1}^{3}   f_h\big(\delta (k(Y)  +i )\big)} \notag \\
& \hspace{5cm} +   C \delta^{-2}  1(Y_1 \leq \delta) \max_{\substack{i \in I \\ i_1 = 0}} \abs{(\Delta_1^2 - (\Delta_1+\Delta_2))    f_h\big(\delta (k(Y)  +i )\big)}, \label{eq:blist2}\\
&\abs{\delta Y_2\frac{\partial^2 }{\partial x_2^2} A \widehat f_h(\xi^3)} \leq  C \delta^{-1}Y_2 \max_{\substack{i \in I  }} \abs{\Delta_{2}^{2}    f_h\big(\delta (k(Y)  +i )\big)} , \text{ and }  \label{eq:blist3}\\
&\abs{\varepsilon_{5}(Y)} \leq  C \max_{\substack{ a_1+a_2 = 2 \\  i \in I }} \abs{\Delta_1^{a_1}\Delta_2^{a_2}  f_h\big(\delta(k(Y)   + i) \big) }. \label{eq:blist4}
\end{align}
 We begin with \eqref{eq:blist1}.   Observe that $(1-\delta(Y_1+Y_2)) \in (0,1/2)$ because $Y \in B$. Furthermore, $Y  < \xi^2 < Y +\delta e^{(1)}$ implies $k(\xi^2) = k(Y) \geq 0$. Combining this with \eqref{eq:multibound2} of Theorem~\ref{thm:interpolant_def}, we get 
\begin{align*}
\abs{ \frac{1}{6}\delta(1-\delta(Y_1+Y_2) ) \frac{\partial^3 }{\partial x_1^3} A \widehat f_h(\xi^2)}  
\leq  C \delta^{-2} \max_{\substack{  i \in I }} \abs{\Delta_{1}^{3}   \widehat f_h\big(\delta (k(Y) +i )\big)} =&\ C \delta^{-2} \max_{\substack{ i \in I }} \abs{\Delta_{1}^{3}    f_h\big(\delta (k(Y) +i )\big)}.
\end{align*} 
We now prove \eqref{eq:blist2}. As before,  $Y - \delta e^{(1)} < \xi^1 < Y$ implies that $k(\xi^1) = k(Y - \delta e^{(1)}) = k(Y )- e^{(1)}$, so  
\begin{align}
\abs{\frac{1}{6}\delta\lambda\frac{\partial^3 }{\partial x_1^3} A \widehat f_h(\xi^1)}   \leq C \delta^{-2} \max_{\substack{ i \in I }} \abs{\Delta_{1}^{3}   \widehat f_h\big(\delta (k(Y)-  e^{(1)} +i )\big)}. \label{eq:midblist2}
\end{align}
Now when $Y \in [0,\delta)$ and $i_1 = 0$, the definition of $\widehat f_h(x^q)$ in \eqref{eq:fhat} implies that  
\begin{align*}
 \widehat f_h\big(\delta (k(Y)-  e^{(1)} +i )\big) =  \widehat f_h\big(-\delta, \delta(k_2(Y) + i_2), 0, \ldots, 0\big) = f_h\big(0, \delta(k_2(Y) + i_2+1), 0, \ldots, 0\big),
\end{align*}
from which we see that $\Delta_1 \widehat f_h\big(-\delta, \delta(k_2(Y) + i_2), 0,\ldots, 0\big) = -\Delta_2 f_h\big(0, \delta(k_2(Y) + i_2), 0,\ldots, 0\big)$, and  therefore
\begin{align*}
 \Delta_1^3 \widehat f_h\big(-\delta, \delta(k_2(Y) + i_2), 0,\ldots, 0\big) =   (\Delta_1^2 - (\Delta_1+\Delta_2)) f_h\big(0, \delta(k_2(Y) + i_2), 0,\ldots, 0\big).
\end{align*} 
Combining this with \eqref{eq:midblist2} implies \eqref{eq:blist2}. To prove \eqref{eq:blist3}, we note that  $k(\xi^3) = k(Y )- e^{(2)}$ because $Y - \delta e^{(2)} < \xi^3 < Y$,   so
\begin{align*}
&\abs{\delta Y_2\frac{\partial^2 }{\partial x_2^2} A \widehat f_h(\xi^3)}   \leq C \delta^{-1}Y_2 \max_{\substack{i \in I }} \abs{\Delta_{2}^{2}   \widehat f_h\big(\delta (k(Y)-  e^{(2)} +i )\big)}.
\end{align*}
The definition of $\widehat f_h(x^q)$ in \eqref{eq:fhat} says that $\Delta_2 \widehat f_h( Y_1,-\delta, 0,\ldots, 0) = \Delta_2 \widehat f_h( Y_1,0,\ldots, 0)$, so  $\Delta_2^2 \widehat f_h( Y_1,-\delta, 0,\ldots, 0) =  0$, implying \eqref{eq:blist3}. Lastly, we prove \eqref{eq:blist4}.  Theorem~\ref{thm:interpolant_def} tells us that $\alpha_{k_j+i_j}^{k_j}(x_j)$  are degree-$7$ polynomials in $(x_j-\delta k_j)/ \delta$ whose coefficients do not depend on $k_j$ or $\delta$,   so there exists a constant $C > 0$ such that $\abs{\alpha_{k_j(x)+i_j}^{k_j(x)}(x_j)} \leq C$ for $j = 1,2$,  so 
\begin{align*}
&\Bigg|\sum_{i_2 = 0}^{4} \alpha_{k_2(x)+i_2}^{k_2(x)}(x_2) \sum_{i_1 = 0}^{4} \alpha_{k_1(x)+i_1}^{k_1(x)}(x_1) \frac{1}{\delta} \Big(\delta (k_2(x)+i_2) - x_2\Big) \notag \\
& \hspace{4cm} \times \Big(-\Delta_2 \widehat f_h\big(\delta (k(x)+i-e^{(2)})\big) + \Delta_2 \widehat f_h\big(\delta (k(x)-e^{(2)})   \Big)\Bigg| \\
\leq&\ C \max_{\substack{ i \in I }} \abs{\Delta_2 \widehat f_h\big(\delta (k(x)+i-e^{(2)})\big) - \Delta_2 \widehat f_h\big(\delta (k(x)-e^{(2)}) }.
\end{align*}
Now
 \begin{align*}
& \Delta_2 \widehat f_h\big(\delta (k(x)+i-e^{(2)})\big) - \Delta_2 \widehat f_h\big(\delta (k(x)-e^{(2)})\big) \\ 
=&\   \Delta_2 \widehat f_h\big(\delta (k(x)+i-e^{(2)})\big) - \Delta_2 \widehat f_h\big(\delta (k(x)+i_2 e^{(2)}-e^{(2)})\big) +  \Delta_2 \widehat f_h\big(\delta (k(x)+i_2 e^{(2)}-e^{(2)})\big) - \Delta_2 \widehat f_h\big(\delta (k(x)-e^{(2)})\big)  \\
 =&\   \sum_{i'_1 = 0}^{i_1-1} \Delta_1\Delta_2 \widehat f_h\big(\delta (k(x) + i'_1 e^{(1)}+ i_2 e^{(2)}-e^{(2)})\big)  +   \sum_{i'_2 = 0}^{i_2-1} \Delta_2^2 \widehat f_h\big(\delta (k(x) + i'_2 e^{(2)}-e^{(2)})\big),
 \end{align*}
implying that 
 \begin{align*}
& C \max_{\substack{ i \in I }} \abs{\Delta_2 \widehat f_h\big(\delta (k(x)+i-e^{(2)})\big) - \Delta_2 \widehat f_h\big(\delta (k(x)-e^{(2)}) }  \\
 \leq&\ C \max_{\substack{i \in I }} \abs{\Delta_2^2 \widehat f_h\big(\delta(k(Y) - e^{(2)} + i) \big) }   + C \max_{\substack{i \in I }} \abs{\Delta_1\Delta_2 \widehat f_h\big(\delta(k(Y) - e^{(2)} + i) \big) }.
 \end{align*}
An identical argument allows us to bound the second term on the right-hand side of \eqref{eq:e5}, yielding 
\begin{align*}
\varepsilon_5(Y) \leq&\ C \max_{\substack{i \in I }} \abs{\Delta_2^2 \widehat f_h\big(\delta(k(Y) - e^{(2)} + i) \big) }   + C \max_{\substack{i \in I }} \abs{\Delta_1\Delta_2 \widehat f_h\big(\delta(k(Y) - e^{(2)} + i) \big) }   \\
&+  C \max_{\substack{i \in I }} \abs{\Delta_1^2 \widehat f_h\big(\delta(k(Y)   + i) \big) } +  C \max_{\substack{i \in I }} \abs{\Delta_1\Delta_2 \widehat f_h\big(\delta(k(Y)   + i) \big) }.
\end{align*}
Using $\Delta_2 \widehat f_h( Y_1,-\delta, 0,\ldots, 0) = \Delta_2 \widehat f_h( Y_1,0,\ldots, 0)$ and  $\Delta_2^2 \widehat f_h( Y_1,-\delta, 0,\ldots, 0) =  0$, we conclude \eqref{eq:blist4}.

\subsubsection{Proving Lemma~\ref{lem:actualinterchange}}
\label{sec:actualinterchangeproof}
To prove Lemma~\ref{lem:actualinterchange}, we need the following result.
\begin{lemma} \label{lem:auxinterchange}
Suppose $P \subset \delta \Z^{d}$ and let $f,g: P \to \R$. Given $\ell \in \Z^{d}$, for those $k$ such that $\delta k \in P$ and $\delta (k+\ell) \in P$, we define
\begin{align*}
F(\delta k) = g(\delta k) \big(f(\delta (k+\ell)) - f(\delta k) \big).
\end{align*}
Then $A F(x)$ is well defined for those $x \in \R^{d}$ such that $\delta (k(x) + i) \in P$ and $\delta(k(x) + \ell + i) \in P$ for all $0 \leq i \leq 4e$, where   $k_{i}(x) = \lfloor x_i/\delta\rfloor$. Furthermore, for all such $x$, 
\begin{align*}
A F(x) =&\ A g(x) \big( A f(x+ \delta \ell) - A f(x) \big) \\
&+ \sum_{i_1, \ldots, i_d = 0}^{4} \bigg(\prod_{j=1}^{d} \alpha_{k_j(x) +i_j}^{k_j(x)   }(x_j)\bigg) \Big(g\big(\delta (k(x)+i)\big) - A g(x)\Big) \notag \\
& \hspace{1.7cm} \times \Big(f\big(\delta (k(x)+\ell+i)\big)-f\big(\delta (k(x)+i)\big) - \big(f\big(\delta (k(x)+\ell )\big)-f (\delta k(x) )\big)\Big). 
\end{align*}
\end{lemma}
\startproof{Proof of Lemma~\ref{lem:auxinterchange}}
The proof is identical to the proof of Proposition 3 of \cite{Brav2022}.
\finishproof

\startproof{Proof of Lemma~\ref{lem:actualinterchange} } First, we prove \eqref{eq:gxfhat}. Any  $x^q=  \in  B \cap \delta \N^{b+1} $ satisfies $x^q_2  \leq \delta (n/2 - 8)$, or $q_2 \leq n/2 - 8$.  It follows from  the definition of $G_X$ in \eqref{eq:gxdef} that for $x^q \in  B \cap \delta \N^{b+1}$,
\begin{align*}
 G_{X} f_h(x^q) =&\ -1(q_1 < n) n \lambda  \Delta_1 f_h(x^q-\delta e^{(1)} )  +  1(q_1= n )  n \lambda  \Delta_2 f(x^q) \notag \\
&+ (q_1 - q_{2})  \Delta_1 f_h(x^q)   - q_{2}   \Delta_{2} f_h( x^q- \delta   e^{(2)} ).
\end{align*}
Note that $q_1 - q_2 = n - (n-q_1) - q_2 = n-(x^q_1+x^q_2)/\delta$, and $q_2 = x^q_2 /\delta$. Although $\Delta_{2} f_h( x^q- \delta   e^{(2)} )$ is technically not defined when $x^q_2= 0$,   we adopt the convention that $1(q_2 = 0)q_{2}   \Delta_{2} f_h( x^q- \delta   e^{(2)} ) = 0$. Using the definition of $\widehat f(x^q)$ in \eqref{eq:fhat}, we have 
\begin{align*}
1(q_2 = 0) q_2  \Delta_{2} f( x^q- \delta   e^{(2)} ) = 0 = 1(q_2 = 0) q_2 \Delta_2 \widehat f(x^q - \delta e^{(2)}). 
\end{align*}
Similarly, since $q_1 = n$ corresponds to $x^q_1 = 0$,   
\begin{align*}
n \lambda   1(q_1= n )  \Delta_2 f(x^q)   = - 1(q_1= n )  n \lambda   \Delta_1 \widehat f(x^q - \delta e^{(1)}),
\end{align*}
which proves \eqref{eq:gxfhat}. To prove \eqref{eq:agfinal}, note that if $g(x^q) = (n-x^q_1 - x^q_2)/\delta$, then $A g(x) = n-(x^q_1+x^q_2)/\delta$. To see why, note that $\Delta_{i}\Delta_{j} g(x^q) = 0$ for any $i,j$,  so Theorem~\ref{thm:interpolant_def} implies that all second-order partial derivatives of $A g(x)$ are zero. Since $Ag(x)$ is twice continuously differentiable, it must be a linear function, and the only linear function that coincides with $g(x^q)$ on the grid is $A g(x) = n-(x^q_1+x^q_2)/\delta$. Similarly, if $g(x^q) = q_2 = x^q_2/\delta$, then $A g(x) = x_2/\delta$. Applying Lemma~\ref{lem:auxinterchange} to each of the three terms on the right-hand side of \eqref{eq:gxfhat} proves \eqref{eq:agfinal}. 
\finishproof

\section{Supporting Proofs for Section~\ref{sec:perturbation}}
\label{sec:appendixsteinfactors}
Apart from the short proof of Lemma~\ref{lem:up} in Appendix~\ref{proof:up}, this appendix is devoted to the proof of Lemma~\ref{lem:coupling_time_bound}. Going forward, we fix $\gamma = 2(17/\beta + \beta + 1)$, and recall from  Section~\ref{sec:first_diff} that 
\begin{align*}
&\theta_{1} = n-\lfloor \sqrt{n}\beta/2 \rfloor,  \quad \theta_{2} =  \lfloor \gamma\sqrt{n}\rfloor, \\
&\tau_{i}(q_i) =  \inf \{ t \geq 0 : Q_i(t) = q_i\}, \quad q_i\in \{0,1,\ldots,n\},\ i = 1,2.
\end{align*}
Following the proof roadmap of Lemma~\ref{lem:coupling_time_bound}, we need an upper bound on the expected start of the first cycle  and the expected duration of a single cycle. The following two lemmas provide the ingredients for these bounds  and are proved in Appendices~\ref{sec:claimtwo} and \ref{proof:hittingtime}, respectively.
\begin{lemma} \label{lem:lem16_2}
For all $n \geq 1$, 
\begin{align}
&\max_{\substack{\theta_{1} < q_1 \leq n \\ q_2 = \theta_{2},\ q \in S_Q }}  \E_{q}  \big( \tau_{2}(2\theta_{2}) \wedge \tau_{1}(\theta_{1}) \big)  \leq C(b,\beta). \label{eq:lem16_2} 
\end{align} 
\end{lemma}
\begin{lemma} \label{prop:hittingtime}
For all $n \geq 1$ and $q \in S_Q$ with $q_2 > \theta_{2}$, 
\begin{align*}
 &\E_{q} \tau_2(\theta_{2}) \leq C(b,\beta)  (1+\delta q_2 ) =C(b,\beta) (1+ x^q_2 ), \quad q \in S_Q \text{ with } q_2 > \theta_2.
\end{align*} 
\end{lemma}
To  bound the probability of coupling in a given cycle, we require the following two lemmas.
\begin{lemma} \label{lem:lem16_3}
There exists a constant $p_1(\beta) \in (0,1)$ such that for all $n \geq 1$, 
\begin{align*}
& \min_{\substack{\theta_{1} < q_1 \leq n \\ q_2 = \theta_{2},\  q \in S_Q}} \Prob_{q} \big( \tau_{1}(\theta_{1}) < \tau_{2}(2\theta_{2}) \big)  \geq p_1(\beta).
\end{align*}
\end{lemma} 
\begin{lemma} \label{lem:lem16_4}
There exists a constant $p_2(b,\beta) \in (0,1)$ such that for all $n \geq 1$, 
\begin{align*}
&  \min_{\substack{ 0 \leq q_1 \leq \theta_{1} \\  0 \leq q_2  \leq 2 \theta_{2}   \\ q \in S_Q}} \Prob\Big(\tau_C < \tau_{1}(n) \ \big|\ Q(0) = q,\ (Q(0),\widetilde Q(0)) \in \bigcup_{i=1}^{b+1}\Theta_{i}^{Q}\Big)\geq p_2(b,\beta).
\end{align*}
\end{lemma} 
Lemmas~\ref{lem:lem16_3} and \ref{lem:lem16_4} are proved in Appendices~\ref{sec:claimthree} and \ref{sec:claimfour}, respectively.   
\startproof{Proof of Lemma~\ref{lem:coupling_time_bound} }  
Throughout the proof, we use $C$ to denote a positive constant that may change from line to line  but depends only  on $\beta$ and $b$.  Given any initial condition $(q,\widetilde q) \in \bigcup_{i=1}^{b+1} \Theta_i^{Q}$,
\begin{align*}
\E_{(q, \widetilde q)} \tau_C \leq C(b,\beta)(1+\delta q_2).
\end{align*}
For  convenience,  we abuse notation and   adopt the convention that
\begin{align*}
\E_{q} \tau_C  =&\ \max_{\widetilde q: (q,\widetilde q) \in \bigcup_{i=1}^{b+1}\Theta_{i}^{Q}}  \E \big[ \tau_C \big| (Q(0),\widetilde Q(0)) = (q,\widetilde q)   \big], \quad q \in S_Q,
\end{align*}
but $\E_{q}(W) = \E (W | Q(0)=q)$ for any random variable $W$ other than $\tau_C$. We also assume that every $\max$ operator in this proof automatically considers the maximum over all $q \in S_Q$; i.e.,
\begin{align*}
 \max_{ \substack{  q_2 =\theta_{2}  }} \E_{q}\tau_C  =  \max_{ \substack{q \in S_Q \\ q_2 =\theta_{2}  }} \E_{q}\tau_C.
\end{align*}
Lemma~\ref{prop:hittingtime} implies  that  for any $q \in S_Q$,
\begin{align}
\E_{q} \tau_C \leq  \E_{q} \tau_2(\theta_{2})  + \max_{ \substack{ q_2 =\theta_{2}  }} \E_{q}\tau_C   
 \leq&\ C(1+\delta q_2) + \max_{ \substack{  q_2 =\theta_{2}  }} \E_{q}\tau_C. \label{eq:coupstart}
\end{align} 
We will argue that if $p_1 = p_1(\beta)$ and $p_2 = p_2(b,\beta)$ are the constants from Lemmas~\ref{lem:lem16_3} and \ref{lem:lem16_4}, then
\begin{align}
\max_{\substack{\theta_{1} \leq q_1 \leq n \\ q_2 = \theta_{2}}}\E_{q}  \tau_{C}  \leq&\ C + (1 -p_1   p_2 ) \max_{ \substack{ q_2 =\theta_{2}  }} \E_{q}\tau_C, \text{ and }  \label{eq:toprove_1_tauc} \\
\max_{\substack{0 \leq q_1 < \theta_{1}\\ q_2 = \theta_{2}}}\E_{q} \tau_{C}   \leq&\ C + (1-p_2)  \max_{ \substack{ q_2 =\theta_{2}  }} \E_{q}\tau_C. \label{eq:toprove_2_tauc} 
\end{align}
As a result, choosing $p_3 = \max \{(1 -p_1 p_2 ),(1-p_2) \} \in (0,1)$ implies that 
\begin{align*}
 \max_{ \substack{  q_2 =\theta_{2}  }} \E_{q}\tau_C  =  \max \bigg\{\max_{\substack{0 \leq q_1 < \theta_{1}\\ q_2 = \theta_{2}}}\E_{q}\tau_C, \max_{\substack{\theta_{1} \leq q_1 \leq n \\ q_2 = \theta_{2}}}\E_{q} \tau_C \bigg\} \leq&\ C + p_3  \max_{ \substack{  q_2 =\theta_{2}  }} \E_{q}\tau_C,
\end{align*}
and therefore   $\max_{ \substack{  q_2 =\theta_{2}  }} \E_{q}\tau_C \leq C(1-p_3)^{-1} \leq C$. Combining this with  \eqref{eq:coupstart}   implies the lemma. We now prove \eqref{eq:toprove_1_tauc}, followed by \eqref{eq:toprove_2_tauc}.  Defining $\tau_M = \tau_{2}(2\theta_{2}) \wedge \tau_{1}(\theta_{1})$, we have 
\begin{align}
\max_{\substack{\theta_{1} \leq q_1 \leq n \\ q_2 = \theta_{2}}}\E_{q}  \tau_{C} \leq&\ \max_{\substack{\theta_{1} \leq q_1 \leq n \\ q_2 = \theta_{2}}}\E_{q}  \tau_M + \max_{\substack{\theta_{1} \leq q_1 \leq n \\ q_2 = \theta_{2}}}\E_{q}\big[\E_{ Q (\tau_M) } \tau_C\big]  \leq C + \max_{\substack{\theta_{1} \leq q_1 \leq n \\ q_2 = \theta_{2}}}\E_{q}\big[\E_{ Q (\tau_M) } \tau_C\big], \label{eq:interm_tauc}
\end{align}
where in the second inequality we used \eqref{eq:lem16_2} of Lemma~\ref{prop:hittingtime}. To bound the right-hand side, let us define the events 
\begin{align*}
E_1 =&\ \big\{ \tau_{1}(\theta_{1}) < \tau_{2}(2\theta_{2}) \big\}, \quad \text{ and } \quad  E_2 = \big\{\tau_C < \tau_{1}(n) \big\},
\end{align*}
and their complements $E_1^c$ and $E_2^c$, respectively. Note that if $Q_2(0) < 2\theta_{2}$, then the event   $E_1^c$ implies  that $Q (\tau_M) = (n,2\theta_{2})$  
because $Q_2(t)$ increases only at times when $Q_1(t) = n$.  Using the law of total probability, 
\begin{align}
 \max_{\substack{\theta_{1} \leq q_1 \leq n \\ q_2 = \theta_{2}}}\E_{q}\big[\E_{ Q (\tau_M) } \tau_C\big] \leq&\  \max_{\substack{\theta_{1} \leq q_1 \leq n \\ q_2 = \theta_{2}}}  \Bigg\{ \Prob_{q}(E_1^c)  \max_{\substack{q_1' = n \\ q_2' = 2\theta_{2}  }} \E_{q'}  \tau_{C}  + \Prob_{q}(E_1) \max_{\substack{q_1' = \theta_{1} \\ 0 \leq q_2' \leq 2\theta_{2}    }}  \E_{q'} \tau_C \Bigg\}. \label{eq:partial_tauc0}
\end{align}
We note that
\begin{align}
 \max_{\substack{q_1 = n \\ q_2 = 2\theta_{2}  }} \E_{q}  \tau_{C} \leq \max_{\substack{q_1 = n \\ 0 \leq q_2 \leq 2\theta_{2}   }} \E_{q} \tau_{C} \leq&\ \max_{\substack{q_1 = n \\ 0 \leq q_2 \leq 2\theta_{2}   }} \E_{q} \tau_2(\theta_{2})  + \max_{ \substack{  q_2 =\theta_{2}  }} \E_{q}\tau_C \leq C  + \max_{ \substack{  q_2 =\theta_{2}  }} \E_{q}\tau_C, \label{eq:partial_tauc1}
\end{align}
where we used Lemma~\ref{prop:hittingtime} in the last inequality, so  
\begin{align}
\Prob_{q}(E_1^c)  \max_{\substack{q_1 = n \\ q_2 = 2\theta_{2}  }} \E_{q}  \tau_{C} \leq \Prob_{q}(E_1^c) \big(C  + \max_{ \substack{  q_2 =\theta_{2}  }} \E_{q}\tau_C\big). \label{eq:partial_tauc3}
\end{align}
Provided we can show that 
\begin{align}
\Prob_{q}(E_1) \max_{\substack{q_1  = \theta_{1} \\ 0 \leq q_2  \leq 2\theta_{2}  }}  \E_{q} \tau_C \leq \Prob_{q}(E_1)  \big( C + (1-p_2)\max_{ \substack{  q_2 =\theta_{2}  }} \E_{q}\tau_C\big), \label{eq:partial_tauc2}
\end{align}
we can combine \eqref{eq:partial_tauc3} and \eqref{eq:partial_tauc2} with \eqref{eq:partial_tauc0} to get
\begin{align*}
\max_{\substack{\theta_{1} \leq q_1 \leq n \\ q_2 = \theta_{2}}}\E_{q}\big[\E_{ Q (\tau_M) } \tau_C\big] \leq&\  \max_{\substack{\theta_{1} \leq q_1 \leq n \\ q_2 = \theta_{2}}}  \Bigg\{ \Prob_{q}(E_1^c)  \big( C + \max_{ \substack{  q_2' =\theta_{2}  }} \E_{q'}\tau_C\big) + \Prob_{q}(E_1) \big( C + (1-p_2)\max_{ \substack{  q_2' =\theta_{2}  }} \E_{q'}\tau_C\big) \Bigg\}\\
=&\ C + \Big( \max_{ \substack{  q_2 =\theta_{2}  }}   \E_{q}\tau_C\Big)   \max_{\substack{\theta_{1} \leq q_1 \leq n \\ q_2 = \theta_{2}}} \big\{ 1 - p_2 \Prob_{q}(E_1)  \big\}  \notag \\
\leq&\ C + (1 -p_2  p_1 )\max_{ \substack{  q_2 =\theta_{2}  }}   \E_{q}\tau_C, 
\end{align*}
where the last inequality follows from the lower bound on $\Prob_{q}(E_1)$ in Lemma~\ref{lem:lem16_3}. Combining this bound with \eqref{eq:interm_tauc} proves \eqref{eq:toprove_1_tauc}. We now prove \eqref{eq:partial_tauc2}. Recall that $E_2 = \big\{\tau_C < \tau_{1}(n) \big\}$ and observe that 
\begin{align*}
 &\max_{\substack{q_1  = \theta_{1} \\ 0 \leq q_2  \leq 2\theta_{2}  }}  \E_{q} \tau_C  \\
 \leq&\ \max_{\substack{q_1  = \theta_{1} \\ 0 \leq q_2  \leq 2\theta_{2}  }}  \E_{q} \Big(\big[\tau_C \wedge \tau_1(n) \big] 1(E_2) \Big) + \max_{\substack{q_1  = \theta_{1} \\ 0 \leq q_2  \leq 2\theta_{2}  }}    \E_{q} \Big(\big[\tau_C \wedge \tau_1(n) + \E_{Q(\tau_1(n))}\tau_C   \big]1(E_2^c)  \Big)  \\
\leq&\ 2 \max_{\substack{q_1  = \theta_{1} \\ 0 \leq q_2  \leq 2\theta_{2}  }}  \E_{q} \big[\tau_C \wedge \tau_1(n) \big]  + \max_{\substack{q_1  = \theta_{1} \\ 0 \leq q_2  \leq 2\theta_{2}  }}  \Prob\Big( E_2^c \ \big|\ Q(0)=q,\ (Q(0),\widetilde Q(0)) \in \bigcup_{i=1}^{b+1} \Theta_i^Q \Big) \E_{q} \big[ \E_{Q(\tau_1(n))}\tau_C  \big]\\
\leq&\ 2\max_{\substack{q_1  = \theta_{1} \\ 0 \leq q_2  \leq 2\theta_{2}  }}  \E_{q} \big[\tau_C \wedge \tau_1(n) \big] + (1-p_2)  \max_{\substack{q_1=n \\ 0 \leq q_2\leq 2\theta_{2}}} \E_{q}\tau_C,
\end{align*}
where in the last inequality we used Lemma~\ref{lem:lem16_4} and the fact that $Q_2(\tau_1(n)) \leq Q_2(0)$ because $Q_2(t)$ increases only at times when  $Q_1(t) = n$. 
 Applying \eqref{eq:partial_tauc1} to the right-hand side, we arrive at 
\begin{align*}
\max_{\substack{q_1  = \theta_{1} \\ 0 \leq q_2  \leq 2\theta_{2}  }}  \E_{q} \tau_C  \leq&\ 2\max_{\substack{q_1  = \theta_{1} \\ 0 \leq q_2  \leq 2\theta_{2}  }}  \E_{q} \big[\tau_C \wedge \tau_1(n) \big] +C + (1-p_2) \max_{\substack{ q_2 = \theta_{2} }} \E_{q}\tau_C.  
\end{align*}
To conclude, we argue that   
\begin{align}
\max_{\substack{q_1  = \theta_{1} \\ 0 \leq q_2  \leq 2\theta_{2}  }}  \E_{q} \big[\tau_C \wedge \tau_1(n) \big] \leq b+1. \label{eq:bbound}
\end{align}
 If $(Q(0),\widetilde Q(0)) \in \Theta_1^Q$, then $(Q(t),\widetilde Q(t)) \in \Theta_1^Q$ for all $t\in [0,\tau_1(n)]$ by construction. The joint CTMC couples before $\tau_1(n)$  if   $\tau_1(n) > V$, where $V$ is as in \eqref{eq:coup_time_alternative}. If $(Q(0),\widetilde Q(0)) \in \Theta_i^{Q}$ for $i \geq 2$, coupling will happen before $\tau_1(n)$ if the joint CTMC  transitions to $\Theta_1^Q$  and then spends $V$ time units there, all before $\tau_1(n)$. From the construction of $\widetilde Q(\cdot)$, we know that the time taken to get from $\Theta_i^Q$ to $\Theta_1^Q$ equals the sum of $i-1$ unit-mean exponentially distributed random variables,  so the worst case is when $i = b+1$. Letting $\Gamma_{b+1}$ represent this sum, it follows that
\begin{align*}
\max_{\substack{q_1  = \theta_{1} \\ 0 \leq q_2  \leq 2\theta_{2}  }}  \E_{q} \big[\tau_C \wedge \tau_1(n) \big] \leq \max_{\substack{q_1  = \theta_{1} \\ 0 \leq q_2  \leq 2\theta_{2}  }}  \E_{q} \big[\Gamma_{b+1} \wedge \tau_1(n) \big] \leq \E(\Gamma_{b+1}) \leq b+1,
\end{align*}
which proves \eqref{eq:bbound}. Our argument for \eqref{eq:partial_tauc2} can be repeated to prove \eqref{eq:toprove_2_tauc}.
\finishproof

\subsection{Proving Lemma~\ref{lem:lem16_2}} 
\label{sec:claimtwo}
\startproof{Proof of  Lemma~\ref{lem:lem16_2}} 
Define $V(x^q) = \sum_{i=1}^{b+1} q_i$ and observe that
\begin{align*}
G_X V(x^q) = n\lambda 1(q_{b+1}<n) -q_1, \quad x^q \in S.  
\end{align*}
Since $\theta_{1} = n-\lfloor \sqrt{n}\beta/2 \rfloor$, it follows that   for any  $q \in S_Q$ with $ \theta_{1} <  q_1 \leq n $,
\begin{align*}
G_X V(x^q) = n\lambda - q_1 \leq n\lambda - (n-\lfloor \sqrt{n}\beta/2 \rfloor)= -\beta\sqrt{n}+ \lfloor \sqrt{n}\beta/2 \rfloor \leq  -\sqrt{n}\beta/2.
\end{align*}
Let $M> 0$, $t^{(M)} = \min \{\tau_{1}(\theta_{1}), \tau_{2}(2\theta_{2}), M  \}$, and note that $Q_1(t) \geq n-\lfloor \sqrt{n}\beta/2 \rfloor$ for $t \leq t^{(M)}$.   Dynkin's formula, e.g.,  Lemma~17.2 in \cite{Kall2001},  then   implies that for any $q \in S_Q$ with $ \theta_{1} <  q_1 \leq n $ and $q_2 = \theta_{2}$, 
\begin{align*}
\E_{x^q} V(X(t^{(M)})) - V(x^q) = \E_{x^q} \int_{0}^{t^{(M)}} G_X V(X(s)) ds \leq - \frac{\sqrt{n}\beta}{2} \E_{x^q} t^{(M)}.
\end{align*}
Since $Q_1(t^{(M)}) \geq   n-\lfloor \sqrt{n}\beta/2 \rfloor$ and $ \theta_{1} <  q_1 \leq n $, it follows that  $q_1 - Q_1(t^{(M)})  \leq \lfloor \sqrt{n}\beta/2 \rfloor$,  so
\begin{align*}
\frac{\sqrt{n}\beta}{2} \E_{x^q} t^{(M)} \leq  V(x^q) - \E_{x^q} V(X(t^{(M)})) \leq  q_1 - \E_{x^q} Q_1(t^{(M)}) + \sum_{i=2}^{b+1} q_i \leq \lfloor \sqrt{n}\beta/2 \rfloor+ b \theta_2,
\end{align*}
where in the last inequality we used $q_2 \geq q_3 \geq \ldots \geq q_{b+1}$. Dividing both sides by $\sqrt{n}$, and noting that $\theta_2 / \sqrt{n} \leq \gamma = 2(17/\beta + \beta + 1)$, yields $ \E_{x^q} t^{(M)}  \leq C(b,\beta)$.  We conclude by taking $M \to \infty$ and using the monotone convergence theorem.
\finishproof

\subsection{Proving Lemma~\ref{prop:hittingtime}}
\label{proof:hittingtime}
Recall that $\theta_2 =\lfloor \gamma\sqrt{n}\rfloor$ and $\gamma = 2(17/\beta + \beta + 1)$.  In this section  we show that $\E_{q} \tau_2(\theta_2) \leq C(b,\beta)(1+\delta q_2)$ if $q_2 > \theta_2$. Our proof is based on a Lyapunov function characterized by the following proposition, proved in Appendix~\ref{sec:oldpaper}. 
 \begin{lemma}
 \label{lem:fluid_hit}
There exists a function $V : \R^{b+1}_+ \to \R$  such that for  any $n \geq 1$ and any $x^q \in S$ with  $x^q_2 \geq 2(17/\beta + \beta) + \delta $,
\begin{align}
G_X V(x^q) \leq&\ -3/17 + \frac{\delta}{\beta}\big(q_3 1(b>1) -  n\lambda1(q_1 = q_2 = n)\big). \label{eq:fluid_drift}
\end{align} 
Furthermore, there exists a constant $C(\beta) > 0$ such that for  any $n \geq 1$,
\begin{align*}
0 \leq V(x) \leq&\ C(\beta) (1+ x_2), \quad  x \in \R^{b+1} \text{ with } x_2 \geq 2(17/\beta+\beta).
\end{align*}
\end{lemma}
\startproof{Proof of Proposition~\ref{prop:hittingtime}}
Let $V(x)$ be the function in Lemma~\ref{lem:fluid_hit}, fix $X(0) = x^q \in S$ with $x^q_2 \geq \delta \theta_2$, $M > 0$, and define $\tau_{2}^{M}(\theta_2)  =  M \wedge \tau_2(\theta_2)$. Dynkin's formula says that 
\begin{align}
&\E_{x^q} V \big( X\big(\tau_{2}^{M}(\theta_2)\big)\big) - V(x^q) = \E_{x^q}\int_{0}^{\tau_{2}^{M}(\theta_2)} G_X V (X(t)) dt.\label{eq:dynkin}
\end{align}
Since  $X_2(t) \geq \delta \theta_2  \geq 2(17/\beta + \beta )+\delta$  for all $t \in [0,\tau_2^{M}(\theta_2)]$,  Lemma~\ref{lem:fluid_hit}  implies that 
\begin{align*}
G_X V(X(t)) \leq -3/17 +\frac{\delta}{\beta} \big(Q_3(t) 1(b>1) -  n\lambda1(Q_1(t) = Q_2(t) = n)\big) , \quad t \in [0,\tau_2^{M}(\theta_2)].
\end{align*}
Combining this inequality with \eqref{eq:dynkin} and  that $V \big( X\big(\tau_{2}^{M}(\theta_2)\big)\big) \geq 0$ and  $V(x^q) \leq C(\beta) x^q_2$ yields
\begin{align*}
  \frac{3}{17} \E_{x^q} \big( \tau_{2}^{M}(\theta_2)\big) \leq C(\beta) x^q_2+ \frac{\delta}{\beta} \E_{x^q} \int_{0}^{\tau_2^{M}(\theta_2)}\big(Q_3(t) 1(b>1) -  n\lambda1(Q_1(t)=Q_2(t)= n)\big) dt.
\end{align*}
If $b = 1$, the lemma follows trivially, so we assume that $b > 1$.  It suffices to show that 
\begin{align*}
\E_{x^q} \int_{0}^{M}\big(Q_3(t) 1(b>1) -  n\lambda1(Q_1(t)=Q_2(t)= n)\big) dt \leq \sum_{i=3}^{b+1} q_i 
\end{align*}
because $\sum_{i=3}^{b+1} q_i \leq b q_2$. Since $Q_3(t)$ is the number of servers with at least two customers in their buffers, it is also  the number of customers that are second in line at time $t$. Thus, $ \int_{0}^{M} Q_3(t)  dt$ is the cumulative time spent by customers being second in line. This cumulative time is contributed to by customers already in the system at time $t = 0$ and by new arrivals after $t = 0$. Of those customers present in the system at $t = 0$, the number that are, or could at some point become, second in line is $\sum_{i=3}^{b+1} q_i$, and each will spend at most one unit of time being second in line, in expectation.

Let $N$ be the number of customers in the interval $[0,M]$ that arrive  when all servers are busy  and all queues have at least one customer in them; i.e., $Q_1(t) = Q_2(t) = n$. For $1 \leq i \leq N$,   let $\xi_i$ be the time customer $i$ spends being second in line, even if that customer becomes second in line after time $M$. We argue  that conditioned on $\{N \geq i\}$, each $\xi_i$ is exponentially distributed with unit mean. Upon entry into the system, if customer $i$ is routed to a busy server with only one other customer waiting in the buffer, then $\xi_{i}$ is distributed according to the remaining service time of the server, which is exponentially distributed with unit mean. If the buffer has more than one customer waiting, then $\xi_{i}$ equals  the service time of the customer two spots ahead of customer $i$, which is also exponentially distributed with unit mean. Further note that the CTMC can  be constructed in such  a way that the value of $\xi_{i}$ is determined at the instant when customer $i$ enters the system, so 
\begin{align*}
\E_{x^q} \int_{0}^{M} Q_3(t)  dt \leq \sum_{i=3}^{b+1} q_i +  \E_{x^q}\sum_{i=1}^{\infty} \xi_i 1(N \geq i) =\sum_{i=3}^{b+1} q_i +  \sum_{i=1}^{\infty} \E_{x^q} \big(\xi_i| N \geq i\big) \Prob(N \geq i) = \sum_{i=3}^{b+1} q_i +\E_{x^q} N.
\end{align*}
Let $\eta_{i}$ be the time spent by the CTMC in a state with $Q_1(t) = Q_2(t) = n$ before customer $i$'s arrival. Since the arrivals to the JSQ system are governed by a rate-$n\lambda$ Poisson process, the arrival of customer $i$ corresponds to a time when $\eta_{i}$ accumulates to equal an exponentially distributed random variable with rate $n\lambda$, and therefore 
\begin{align*}
 \E_{x^q} \int_{0}^{M} 1(Q_1(t)=Q_2(t)= n)  dt  \geq&\   \E_{x^q} \sum_{i=1}^{\infty} \eta_i 1(N \geq i)= \sum_{i=1}^{\infty} \E_{x^q}\big(\eta_i| N \geq i\big) \Prob(N \geq i)= \frac{1}{n\lambda}\E_{x^q} N.
\end{align*}
\finishproof

\subsubsection{Proving Lemma~\ref{lem:fluid_hit}.}
\label{sec:oldpaper}
The Lyapunov function in Lemma~\ref{lem:fluid_hit} is based on the fluid limit of the JSQ system, studied in \cite{Brav2020}. Lemma~\ref{lem:fluid_hit} was, unfortunately, not proved there, but that paper contains all the necessary ingredients for the proof. We now recall them, using notation from \cite{Brav2020}. 

Consider the two-dimensional process $\{(Q_1(t)-n)/n, Q_2(t)/n\}$. Note that the first coordinate is nonpositive, whereas so far  we have been using a nonnegative first coordinate. Section~4.1 of \cite{Brav2020} described the fluid limit of this process. Letting
\begin{align*}
\Omega = \{x \in \R^{2} : x_1 \leq 0,\ x_2 \geq 0\},
\end{align*}
the fluid limit is a dynamical system $v: \R_+ \to \Omega$ with initial condition $v(0) = x \in \Omega$; we write $v^{x}(t)$ to emphasize the relationship on $x$. Postponing the discussion of the behavior of $v^{x}(t)$, for $\ell,u \in \R$ with $\ell < u$  define the smoothed indicator function $\phi^{(\ell,u)}: \R \to [0,1]$ by
\begin{align}
\phi^{(\ell,u)}(x) = 
\begin{cases}
0, \quad &x \leq \ell,\\
(x-\ell)^2\Big( \frac{-(x-\ell)}{ ((u+\ell)/2-\ell )^2(u-\ell)} + \frac{2}{ ((u+\ell)/2-\ell)(u-\ell)}\Big), \quad &x \in [\ell,(u+\ell)/2],\\
1 - (x-u)^2\Big( \frac{(x-u)}{ ((u+\ell)/2-u )^2(u-\ell)} - \frac{2}{ ((u+\ell)/2-u )(u-\ell)}\Big), \quad &x \in [(u+\ell)/2,u],\\
1, \quad &x \geq u,
\end{cases} \label{eq:phi}
\end{align}
and let 
\begin{align*}
f^{(2)}(x) = \int_{0}^{\infty} \phi^{(\delta \kappa_1,\delta \kappa_2 )}\big(v^{x}(t)\big) dt, \quad x \in \Omega,
\end{align*}
where $\delta = 1/\sqrt{n}$  and  $\kappa_1,\kappa_2 \in \R$ are to be determined. The function $f^{(2)}(x)$ appeared in Section~5.1 of \cite{Brav2020}, where it was used as a  Lyapunov function for the diffusion limit of the JSQ system; i.e., the process $\{Y(t)\}$ in \eqref{eq:diffusion2}. We show that this is also a Lyapunov function for the CTMC. Define
 \begin{align}
 V(x) = f^{(2)}(-\delta x_1,\delta x_2), \quad x \in \R^{b+1}_+. \label{eq:vf}
 \end{align}
The following result proved in Appendix~\ref{sec:fluidjsqbounds} gives us   control over the derivatives of $V(x)$. 
\begin{lemma} \label{lem:fluidjsqbounds}
 For any $x \in \R^{b+1}_+$  with $x_2 \geq \kappa_2$, 
\begin{align}
&\big(  \beta - ( x_1 + x_2) \big) \frac{\partial V(x)}{\partial x_1}  - \delta x_2  \frac{\partial V(x)}{\partial x_2} = -1,\quad \text{ and } \quad 1(x_1 = 0)\Big(  \frac{\partial }{\partial x_1} + \frac{\partial}{\partial x_2}\Big) V(x) = 0. \label{eq:oldlyap}
\end{align} 
Furthermore, if we choose $\kappa_1 = 17/\beta + \beta$ and $\kappa_2 = 2 \kappa_1$, then for any $x \in \R^{b+1}_+$  with $x_2 \geq \kappa_2$, and any $x_2' \geq x_2$, 
\begin{align}
&\frac{\partial^2 }{\partial x_1^2}V(x) \leq 9/17,\label{eq:oldd1} \\
&\frac{\partial}{\partial x_2}  V(x) \leq \frac{\partial}{\partial x_2}  V(0,x_2')  = \frac{1}{\beta},  \quad  \frac{\partial^2}{\partial x_2^2} V(x) \leq 5/17, \label{eq:oldd2}
\end{align} 
and there exists a constant $C(\beta)$ such that $0 \leq V(x) \leq C(\beta) (1+x_2)$. 
\end{lemma}
\startproof{Proof of Lemma~\ref{lem:fluid_hit}}
Let $\kappa_1 = 17/\beta + \beta$ and $\kappa_2 = 2 \kappa_1$, and $V(x)$ be the function from Lemma~\ref{lem:fluidjsqbounds}, and recall  $G_X$ defined in \eqref{eq:gxdef}. Since $V(x)$ depends only on $x_1$ and $x_2$, 
\begin{align*}
G_X V(x^q) =&\ 1(q_1 < n) n \lambda \big(- \Delta_1 V(x^q- \delta e^{(1)}) \big)  + n \lambda 1(q_1= n, q_{2} < n)  \Delta_2 V(x^q) \notag \\
&+ (q_1 - q_{2})  \Delta_1 V(x^q)  +  (q_2 - q_{3}1(b>1))  \big(- \Delta_2 V(x^q-\delta e^{(2)})\big),  \quad x^q \in S.
\end{align*} 
Using Taylor expansion, we get 
\begin{align*}
-\Delta_1 V(x^q- \delta e^{(1)}) = V(x^q - \delta e^{(1)} ) - V(x^q) =&\ -\delta  \frac{\partial}{\partial x_1} V(x^q) + \int_{x^q_1-\delta}^{x^q_1} (u - (x^q_1 -\delta))\frac{\partial^2}{\partial x_1^2} V(u,x^q_2) du, \\
\Delta_1 V(x^q) = V(x + \delta  e^{(1)}) - V(x^q) =&\ \delta  \frac{\partial}{\partial x_1} V(x^q) + \int_{x^q_1}^{x^q_1+\delta} (x^q_1 + \delta - u) \frac{\partial^2}{\partial x_1^2} V(u,x^q_2) du,
\end{align*}
and  a similar expression holds for $\Delta_2 V(x^q)$ and $-\Delta_2 V(x^q-\delta e^{(2)})$. Therefore, 
\begin{align}
G_X V(x^q) =&\ -\delta \big( 1(q_1 < n) n\lambda - (q_1-q_2) \big) \frac{\partial}{\partial x_1} V(x^q) + \delta \big(1(q_1=n,q_2<n) n\lambda - q_2\big) \frac{\partial}{\partial x_2} V(x^q)  \notag \\
&- q_3 1(b>1) \big(- \Delta_2 V(x^q-\delta e^{(2)})\big) + \psi(x^q), \label{eq:gxfluild}
\end{align} 
where 
\begin{align*}
\psi(x^q) =&\ n\lambda 1(q_1 < n) \int_{x^q_1-\delta}^{x^q_1} (u - (x^q_1 -\delta))\frac{\partial^2}{\partial x_1^2} V(u,x^q_2) du  \notag \\
&+ n\lambda 1(q_1 = n, q_2 < n) \int_{x^q_2}^{x^q_2+\delta} (x^q_2 + \delta - u) \frac{\partial^2}{\partial x_2^2} V(x^q_1,u) du  \notag \\
&+ (q_1 - q_2) \int_{x^q_1}^{x^q_1+\delta} (x^q_1 + \delta - u) \frac{\partial^2}{\partial x_1^2} V(u,x^q_2) du + q_2  \int_{x^q_2-\delta}^{x^q_2} (u - (x^q_2 -\delta))\frac{\partial^2}{\partial x_2^2} V( x^q_1,u) du.
\end{align*}
Now suppose  $x^q_2 \geq \kappa_2 + \delta$. The bounds on the second-order derivatives of $V(x)$  from Lemma~\ref{lem:fluidjsqbounds}, together with the facts that $q_1-q_2 \geq 0$, $q_2 \geq 0$, $\delta^2 n\lambda \leq 1$, and $\delta^2 q_i \leq 1$, imply that $ \psi(x^q) \leq 14/17$. Next, we rewrite the   first line on the right-hand side of \eqref{eq:gxfluild}, for which we note that
\begin{align*}
&\lambda = 1-\beta/\sqrt{n}, \quad x^q_1 = \delta(n-q_1), \quad 1(q_1=n) = 1(x^q_1 = 0),\\
&1(q_1<n) = 1-1(x^q_1=0), \quad 1(q_1=n,q_2<n) = 1(x^q_1=0)-1(q_1 = q_2=n),
\end{align*}
so 
\begin{align*}
&-\delta \big( 1(q_1 < n) n\lambda - (q_1-q_2) \big) \frac{\partial}{\partial x_1} V(x^q) + \delta \big(1(q_1=n,q_2<n) n\lambda - q_2\big) \frac{\partial}{\partial x_2} V(x^q) \\
=&\   \big( \beta - (x^q_1 + x^q_2) \big) \frac{\partial}{\partial x_1} V(x^q) - x^q_2\frac{\partial}{\partial x_2} V(x)   -1(q_1 = q_2=n) \delta n\lambda   \frac{\partial}{\partial x_2} V(x^q) + \delta n\lambda 1(x^q_1=0) \Big(\frac{\partial}{\partial x_1} +  \frac{\partial}{\partial x_2}\Big) V(x^q)\\
=&\ -1 -1(q_1 = q_2=n) \delta n\lambda   \frac{\partial}{\partial x_2} V(x^q),
\end{align*} 
where the last equality is due to \eqref{eq:oldlyap} from Lemma~\ref{lem:fluidjsqbounds}. We have thus shown that 
\begin{align*}
G_X V(x^q) \leq&\ -1  + 14/17 -1(q_1 = q_2=n) \delta n\lambda   \frac{\partial}{\partial x_2} V(x^q)- q_3 1(b>1) \big(- \Delta_2 V(x^q-\delta e^{(2)})\big).
\end{align*}
Now $V(x^q) = V(0,\sqrt{n})$ when $q_1 = q_2 =n$,  so \eqref{eq:oldd2} in Lemma~\ref{lem:fluidjsqbounds} tells us that $V(0,\sqrt{n}) = 1/\beta$ provided that $\sqrt{n} \geq \kappa_2 = 2(\beta/17+\beta)$, which we assume,  so 
\begin{align*}
&-1(q_1 = q_2=n) \delta n\lambda   \frac{\partial}{\partial x_2} V(x)- q_3 1(b>1) \big(- \Delta_2 V(x^q-\delta e^{(2)})\big) \\
=&\ -1(q_1 = q_2=n) \delta n\lambda  \frac{1}{\beta} + q_3 1(b>1) \int_{x^q_2-\delta}^{x^q_2} \frac{\partial}{\partial x_2} V(x^q_1,u) du \\ 
\leq&\ \frac{\delta }{\beta} \big( q_3 1(b>1) - n\lambda 1(q_1 = q_2=n)\big),
\end{align*} 
where the inequality follows from \eqref{eq:oldd2} in Lemma~\ref{lem:fluidjsqbounds}.

\finishproof

\subsubsection{Proof of Lemma~\ref{lem:fluidjsqbounds}.}
\label{sec:fluidjsqbounds}
Fix  $\kappa_1 = 17/\beta + \beta$ and $\kappa_2 = 2 \kappa_1$. 
The function $f^{(2)}(x)$ was considered in Lemma~8 of \cite{Brav2020}, which tells us that that  
\begin{align*}
&-(\beta \delta + x_1 - x_2) \frac{\partial}{\partial x_1} f^{(2)}(x) - x_2 \frac{\partial}{\partial x_2} f^{(2)}(x) = -1, \quad &x \in \Omega \text{ with } x_2 > \kappa_2/\sqrt{n},\\
&\frac{\partial}{\partial x_1} f^{(2)}(x) = \frac{\partial}{\partial x_2} f^{(2)}(x), \quad &x \in \Omega \text{ with } x_1 =0.
\end{align*} 
Combining this with  
\begin{align}
\frac{\partial}{\partial x_1} V(x) = -\delta \frac{\partial}{\partial x_1} f^{(2)}(-\delta x_1,\delta x_2), \quad \frac{\partial}{\partial x_2} V(x) =  \delta \frac{\partial}{\partial x_2} f^{(2)}(-\delta x_1,\delta x_2) \label{eq:vfrelationship}
\end{align}
gives us \eqref{eq:oldlyap}.  Going forward, we assume that $x \in \Omega$. Let us bound the derivatives of $V(x)$.  On page~1100 of \cite{Brav2020}, it was shown that 
\begin{align*}
\frac{\partial^2}{\partial x_1^2} f^{(2)}(x) \leq \frac{n}{\beta(\kappa_1-\beta)} + \frac{\kappa_1}{\kappa_1-\beta} \frac{4n}{\beta(\kappa_2-\kappa_1)} = \frac{n}{17} + \frac{17/\beta + \beta}{17/\beta} \frac{4n}{\beta(17/\beta + \beta)} = \frac{5n}{17}, \quad x \in \Omega,
\end{align*}
implying the bound on $\partial^2 V(x) /\partial x_1^2$ in \eqref{eq:oldd1}. We now prove \eqref{eq:oldd2}, followed by the bound on $V(x)$.  Unfortunately, $\partial f^{(2)}(x) /\partial x_2$ and $\partial^2 f^{(2)}(x) /\partial x_2^2$ are not bounded in \cite{Brav2020},   so we must bound these partial derivatives ourselves.

We write the equation for  $\partial f^{(2)}(x) /\partial x_2$ in \eqref{eq:speciald2} below, but writing it requires us to introduce some nontrivial objects from \cite{Brav2020}.     The first object we need is the family of curves $\{\Gamma^{(\kappa)} \subset \Omega\}_{\kappa \geq \beta}$, where $\Gamma^{(\kappa)}$ is the graph of the unique fluid-limit trajectory that intersects the $x_2$ axis at the point $(0,\kappa/\sqrt{n})$. For the purposes of this proof, it suffices to treat $\Gamma^{(\kappa)}$ as a two-dimensional geometric object satisfying the following properties:
\begin{enumerate}
\item $\Gamma^{(\kappa)}$ is a graph of a continuous function; i.e.\ $\Gamma^{(\kappa)} = \{(x_1,f(x_1)\}$ for some continuous function $f: \R_+ \to \R_+$. 
\item $\Gamma^{(\kappa)} \cap \{x \in \Omega : x_1 = 0\} = (0,\kappa/\sqrt{n})$.
\item If $x \in \Gamma^{(\kappa)}$ and $x_1 < 0$, then $x_2 > \kappa/\sqrt{n}$. 
\item  If $\kappa' > \kappa$, then $\Gamma^{(\kappa)} \cap \Gamma^{(\kappa')} = \emptyset$ and $\Gamma^{(\kappa')}$ lies above $\Gamma^{(\kappa)}$.  
\end{enumerate}
The first three properties are implied by Lemma~5 of \cite{Brav2020}, and the fourth one follows from (39) there. 
Since $\Gamma^{(\kappa)}$ is a graph, sets of the form $\{x < \Gamma^{(\kappa)}\}$, $\{x \leq \Gamma^{(\kappa)}\}$, etc.,   are well defined. Let us use $\Gamma^{(\kappa_1)}$ and $\Gamma^{(\kappa_2)}$ to partition $\Omega$ into the four sets
\begin{align*}
&S_0 =  \{x \in \Omega \ :\ x_2 \leq \kappa_1/\sqrt{n} \}, \quad S_1 =  \{x \in \Omega \ :\ x_2 \geq \kappa_1/\sqrt{n},\ x \leq \Gamma^{(\kappa_1)} \},\\
&S_2 = \{x \in \Omega \ :\   \Gamma^{(\kappa_1)} \leq x \leq \Gamma^{(\kappa_2)} \}, \quad S_3 =  \{x \in \Omega \ :\ x \geq \Gamma^{(\kappa_2)}\}.
\end{align*}
The four properties of $\Gamma^{(\kappa)}$ are sufficient to argue that $S_0 \cup S_1\cup S_2\cup S_3 = \Omega$ and that the interiors of $S_i$ and $S_j$ are disjoint when $i \neq j$; we refer the reader to Section~C.2 of \cite{Brav2020} for more details. 

The last object we need  is the function $\tau(x)$,  which represents the first time that the fluid limit hits the $x_2$ axis starting from a state $x > \Gamma^{(\beta)}$. The precise definition of $\tau(x)$ is bulky and involves the Lambert-W function, but we can get by with only a few of its properties. Namely, for any $\kappa > \beta$,  Lemma~6 of \cite{Brav2020} introduces a nonnegative function $\tau: \{x \in \Omega : x \geq \Gamma^{(\kappa)}\} \to \R_+$ with $\tau(0,x_2) = 0$, which is differentiable for all $x \in \{x \in \Omega : x \geq \Gamma^{(\kappa)}\}$ and satisfies
\begin{align}
\frac{\partial}{\partial x_1} \tau(x) = - \frac{e^{-\tau(x)}}{x_2 e^{-\tau(x)} - \beta/\sqrt{n}} \leq 0, \quad \frac{\partial}{\partial x_2} \tau(x) = \tau(x) \frac{\partial}{\partial x_1}\tau(x)  \leq 0, \quad x \in \{x \in \Omega : x \geq \Gamma^{(\kappa)}\}. \label{eq:tauder}
\end{align}
By choosing $\kappa = \kappa_1 = 17/\beta + \beta$, we are assured that $\tau(x)$ is defined on the set $\{x \in \Omega : x \geq \Gamma^{(\kappa_1)}\} = S_2 \cup S_3$. Item~1 of Lemma~6 in \cite{Brav2020} tells us that $\tau(x)$ is tied to $\Gamma^{(\kappa)}$ for  any $\kappa  > \beta$ via 
\begin{align}
x_2 e^{-\tau(x)} \geq \kappa/\sqrt{n}, \quad x \geq \Gamma^{(\kappa)}. \label{eq:taugamma}
\end{align} 
We are now ready to bound the derivatives of $f^{(2)}(x)$.  Equation (C.9) of \cite{Brav2020} tells us that 
\begin{align}
\frac{\partial}{\partial x_2} f^{(2)}(x) = \begin{cases}
0, \quad &x \in S_0, \\
\frac{1}{x_2} \phi(x_2) ,  &x \in S_1,\\
\frac{1}{x_2} \big( \phi(x_2) - \phi(x_2 e^{-\tau(x)})\big)  + \phi(x_2 e^{-\tau(x)})\frac{\sqrt{n}}{\beta}  e^{-\tau(x)} (\tau(x) + 1) ,& x \in S_2, \\ 
\frac{\sqrt{n}}{\beta}  e^{-\tau(x)} (\tau(x) + 1), \quad &x \in S_3,
\end{cases}\label{eq:speciald2}
\end{align}
where $\phi(x) = \phi^{(\delta \kappa_1,\delta \kappa_2)}(x)$ is the smoothed indicator defined in \eqref{eq:phi}. By differentiating both sides of \eqref{eq:phi}, it is straightforward to check that $\phi(x)$ is non-decreasing, and 
\begin{align}
&\ \phi'(x) \leq \frac{4}{\delta(\kappa_2 - \kappa_1)} = \frac{4\sqrt{n}}{17/\beta + \beta}. \label{eq:phider}
\end{align}
Let us now argue that   $\partial f^{(2)}(x) /\partial x_2 \leq \sqrt{n}/\beta$ for any $x \in \Omega$. If $x \in S_3$, this bound is implied by the inequality $e^{-t}(t+1) \leq 1$ for $t \geq 0$. If $x \in S_1$, the bound is implied by the facts that $\phi(x_2) \leq 1$ and $1/x_2 \leq \sqrt{n}/\kappa_1 \leq \sqrt{n}/\beta$. If $x \in S_2$, we note  that $\phi(x_2) - \phi(x_2 e^{-\tau(x)}) \geq 0$, and $1/x_2   \leq \sqrt{n}/\beta$, meaning that
\begin{align*}
\frac{\partial}{\partial x_2} f^{(2)}(x) =&\  \frac{1}{x_2} \big( \phi(x_2) - \phi(x_2 e^{-\tau(x)})\big)  + \phi(x_2 e^{-\tau(x)})\frac{\sqrt{n}}{\beta}  e^{-\tau(x)} (\tau(x) + 1) \\
 \leq&\ \frac{\sqrt{n}}{\beta} \big( \phi(x_2) - \phi(x_2 e^{-\tau(x)})\big)  + \phi(x_2 e^{-\tau(x)})\frac{\sqrt{n}}{\beta}   = \frac{\sqrt{n}}{\beta}.
\end{align*}
Observe that  $\partial f^{(2)}(x) /\partial x_2 = \sqrt{n}/\beta$ when $\tau(x) = 0$, which is true for any $x \in S_2 \cup S_3$ with $x_1 = 0$, implying the claim about $\partial V(x)/\partial x_2$ in \eqref{eq:oldd2}. To conclude the proof, it remains to show $\partial^2 V(x) /\partial x_2^2 \leq 9/17$ by differentiating both sides in \eqref{eq:speciald2}.  Note that $\partial^2 f^{(2)}(x)/\partial x_2^2 = 0$ for $x \in S_0$. When $x \in S_1$, we use the bound on $\phi'(x)$ in \eqref{eq:phider}, as well as the fact that  $1/x_2 \leq \sqrt{n}/\beta$, to see that
\begin{align*}
\frac{\partial^2}{\partial x_2^2} f^{(2)}(x) = -\frac{1}{x_2^2}\phi(x_2) + \frac{1}{x_2} \phi'(x_2) \leq  \frac{1}{x_2} \phi'(x_2) \leq \frac{\sqrt{n}}{\beta} \frac{4\sqrt{n}}{17/\beta + \beta}  \leq \frac{4n}{17}, \quad x \in S_1.
\end{align*}
 When  $x \in S_3$, 
\begin{align*}
\frac{\partial^2}{\partial x_2^2} f^{(2)}(x) =&\  -   \frac{\sqrt{n}}{\beta}  e^{-\tau(x)} (\tau(x) + 1)\frac{\partial}{\partial x_2} \tau(x)  + \frac{\sqrt{n}}{\beta}  e^{-\tau(x)}\frac{\partial}{\partial x_2}  \tau(x)  =  - \frac{\sqrt{n}}{\beta}  e^{-\tau(x)} \tau(x)\frac{\partial}{\partial x_2} \tau(x).
\end{align*} 
Using the expression for $\partial \tau(x) /\partial x_2$ in \eqref{eq:tauder}, we see that 
\begin{align}
\frac{\partial^2}{\partial x_2^2} f^{(2)}(x) =&\ \frac{e^{-\tau(x)}}{x_2 e^{-\tau(x)} - \beta/\sqrt{n}}\tau^2(x) \frac{\sqrt{n}}{\beta}  e^{-\tau(x)}  \leq  \frac{n}{7\beta(\kappa_2-\beta)} \leq \frac{n}{7\beta(34/\beta + \beta)}  \leq \frac{4n}{17}, \quad x \in S_3. \label{eq:s3repeatarg}
\end{align}
The first inequality follows from $x_2 e^{-\tau(x)} \geq \kappa_2 /\sqrt{n}$ due to \eqref{eq:taugamma}  and the fact that $t^2 e^{-2t} \leq 1/7$ for $t \geq 0$.  Lastly, we consider the case when   $x \in S_2$, for which we recall that
\begin{align}
\frac{\partial}{\partial x_2} f^{(2)}(x) = \frac{1}{x_2} \big( \phi(x_2) - \phi(x_2 e^{-\tau(x)})\big)  + \phi(x_2 e^{-\tau(x)})\frac{\sqrt{n}}{\beta}  e^{-\tau(x)} (\tau(x) + 1), \quad x \in S_2. \label{eq:s2recall}
\end{align}
To help organize terms, let $g(x_2) = x_2 e^{-\tau(x)}$  and note from \eqref{eq:tauder}  that 
\begin{align*}
  g'(x_2) =  e^{-\tau(x)} \Big(1  - x_2  \frac{\partial}{\partial x_2}\tau(x)\Big) =&\ e^{-\tau(x)} \Big( 1 +  \frac{x_2 e^{-\tau(x)}}{x_2 e^{-\tau(x)} - \beta/\sqrt{n}}   \tau(x) \Big) \\
  =&\ e^{-\tau(x)} \Big( 1 +  \tau(x) +  \frac{\beta/\sqrt{n}}{x_2 e^{-\tau(x)} - \beta/\sqrt{n}}   \tau(x) \Big).
\end{align*}
We see that $g'(x_2) \geq 0$ because $\tau(x) \geq 0$ and  $x_2 e^{-\tau(x)} \geq \kappa_1/\sqrt{n}$ for $x \in S_2$ due to \eqref{eq:taugamma}. Furthermore, since $e^{-t} t \leq 1$ and  $e^{-t}(t+1) \leq 1$ for $t \geq 0$, we conclude that 
\begin{align}
 0 \leq g'(x_2) \leq 1 + \frac{\beta/\sqrt{n}}{x_2 e^{-\tau(x)} - \beta/\sqrt{n}}  \leq 1 + \frac{\beta}{\kappa_1-\beta} = 1 + \frac{\beta^2}{17}. \label{eq:gbound}
\end{align}
Let us now differentiate and bound each term on the right-hand side of \eqref{eq:s2recall}  individually. First, 
\begin{align*}
\frac{\partial}{\partial x_2}  \Big(\frac{1}{x_2} \big( \phi(x_2) - \phi(x_2 e^{-\tau(x)})\big)\Big) =&\   -\frac{1}{x_2^2}\big( \phi(x_2) - \phi(x_2 e^{-\tau(x)})\big) + \frac{1}{x_2} \big( \phi'(x_2) - g'(x_2) \phi'(x_2 e^{-\tau(x)})\big) \\
\leq&\ \frac{1}{x_2} \phi'(x_2) \leq \frac{\sqrt{n}}{\beta} \frac{4\sqrt{n}}{17/\beta + \beta} = \frac{4n}{17}.
\end{align*}
The first inequality is due to the facts that $\phi(x)$ is non-decreasing and $g'(x_2) \geq 0$, and the second inequality follows from the fact that $1/x_2 \leq \sqrt{n}/\beta$ and  the bound on $\phi'(x)$ in \eqref{eq:phider}. Differentiating the second term in \eqref{eq:s2recall}, we get  
\begin{align*}
\frac{\partial}{\partial x_2}  \Big(\phi(x_2 e^{-\tau(x)})\frac{\sqrt{n}}{\beta}  e^{-\tau(x)} (\tau(x) + 1)\Big) =&\ \phi(x_2 e^{-\tau(x)})\frac{\partial}{\partial x_2}  \Big(\frac{\sqrt{n}}{\beta}  e^{-\tau(x)} (\tau(x) + 1) \Big)  \\
&+ \phi'(x_2 e^{-\tau(x)}) g'(x_2) \frac{\sqrt{n}}{\beta}  e^{-\tau(x)} (\tau(x) + 1).
\end{align*}
To bound the first term, we use the fact that $\phi(x) \leq 1$, and we repeat the argument used to prove \eqref{eq:s3repeatarg} to see that
\begin{align*}
\phi(x_2 e^{-\tau(x)})\frac{\partial}{\partial x_2}  \Big(\frac{\sqrt{n}}{\beta}  e^{-\tau(x)} (\tau(x) + 1) \Big) \leq \frac{n}{7\beta(\kappa_1 - \beta)} = \frac{n}{7\beta(17/\beta)} = \frac{n}{119}.
\end{align*}
Furthermore, the bounds on $\phi'(x)$ and $g'(x_2)$ in \eqref{eq:phider} and \eqref{eq:gbound}, together with the fact that $e^{-t}(t+1) \leq 1$, imply that 
\begin{align*}
\phi'(x_2 e^{-\tau(x)}) g'(x_2) \frac{\sqrt{n}}{\beta}  e^{-\tau(x)} (\tau(x) + 1) \leq  \frac{4\sqrt{n}}{17/\beta + \beta} \Big( 1 + \frac{\beta^2}{17}\Big) \frac{\sqrt{n}}{\beta} = \frac{4\sqrt{n} \beta}{17 + \beta^2}  \frac{17 + \beta^2}{17}    \frac{\sqrt{n}}{\beta} = \frac{4n}{17}.
\end{align*}
Combining the pieces yields $\partial^2 f^{(2)}(x)/\partial x_2^2  \leq 9n/17$, proving \eqref{eq:oldd2}.

 To conclude, we prove that $0 \leq V(x) \leq C(\beta)(1+x_2)$ for $x_2 \geq \kappa_2$  by proving that $0 \leq f^{(2)}(x) \leq C(\beta) (1+\sqrt{n} x_2)$ for $x_2 \geq \kappa_2/\sqrt{n}$.  
The form of $f^{(2)}(x)$ below can be found in Lemma 12 of \cite{Brav2020}:
\begin{align*}
f^{(2)}(x) = \begin{cases}
0, \hfill x \in S_0, \\
  \int_{0}^{\log(\sqrt{n}x_2/\kappa_1) } \phi( x_2 e^{-t}  ) dt,\ \hfill x_2 \leq \kappa_2/\sqrt{n},\   x \in S_1,\\
\log(\sqrt{n}x_2/\kappa_2) + \int_{0}^{\log(\kappa_2/\kappa_1) } \phi\big( \frac{\kappa_2}{\sqrt{n}}e^{-t} \big) dt,\  \hfill x_2 \geq \kappa_2/\sqrt{n},\   x \in S_1,\\
 \int_{0}^{\tau(x)  } \phi(x_2e^{-t})dt + \frac{\sqrt{n}}{\beta} \int_{\kappa_1/\sqrt{n}}^{x_2 e^{-\tau(x)}} \phi(t)dt,\ \hfill x_2 \leq \kappa_2/\sqrt{n},\ x \in S_2, \\
\log(\sqrt{n}x_2/\kappa_2) + \int_{\log(\sqrt{n}x_2/\kappa_2)}^{\tau(x)  } \phi \big( x_2  e^{-t} \big)dt  + \frac{\sqrt{n}}{\beta} \int_{\kappa_1/\sqrt{n}}^{x_2 e^{-\tau(x)}} \phi(t)dt, \quad  \hfill  x_2 \geq \kappa_2/\sqrt{n},\  x \in S_2, \\
\tau(x) + \frac{x_2 e^{-\tau(x)} - \kappa_2/\sqrt{n}}{\beta/\sqrt{n}} + \frac{\sqrt{n}}{\beta} \int_{\kappa_1/\sqrt{n}}^{\kappa_2/\sqrt{n}} \phi(t)dt, \hfill  x \in S_3.
\end{cases}
\end{align*} 
The fact that $f^{(2)}(x) \geq 0$ follows from $\phi(x),\tau(x) \geq 0$, the definitions of $S_1$, $S_2$, and $S_3$,  and \eqref{eq:taugamma}. We  combine all the cases above into the single upper bound  
\begin{align}
f^{(2)}(x) \leq&\ \log(\sqrt{n} x_2/\kappa_1) 1(x \in S_1\cup S_2 \cup S_3) + \log(\kappa_2/\kappa_1) + \tau(x)1(x\in S_2\cup S_3) \notag \\
&+ \frac{\sqrt{n}}{\beta}(x_2 e^{-\tau(x)} - \kappa_1/\sqrt{n}) 1(x\in S_2)  + \frac{\sqrt{n}}{\beta}(x_2e^{-\tau(x)} - \kappa_2/\sqrt{n})1(x \in S_3) + \frac{\kappa_2-\kappa_1}{\beta}. \label{eq:f2ineq}
\end{align}
 Using the inequality $\log(t) \leq 1+t$ for $t \geq 0$, and the facts that $\kappa_1 = 17/\beta + \beta$ and $\kappa_2 = 2\kappa_1$, we see that $\log(\kappa_2/\kappa_1) = \log(2)$, $(\kappa_2-\kappa_1)/\beta = 1+17/\beta^2$, 
 \begin{align*}
&\log(\sqrt{n} x_2/\kappa_1) 1(x \in S_1\cup S_2 \cup S_3) \leq 1 + \frac{\sqrt{n} x_2}{\kappa_1} \leq  1 + \frac{\sqrt{n} x_2}{  \beta},  \\
& \frac{\sqrt{n}}{\beta}(x_2 e^{-\tau(x)} - \kappa_1/\sqrt{n}) 1(x\in S_2) \leq \frac{\sqrt{n}x_2}{\beta}, \quad \text{ and } \quad \frac{\sqrt{n}}{\beta}(x_2 e^{-\tau(x)} - \kappa_2/\sqrt{n}) 1(x\in S_3)\leq \frac{\sqrt{n}x_2}{\beta}.
 \end{align*}
 Furthermore, \eqref{eq:taugamma} and the definitions of $S_2$ and $S_3$ imply that
\begin{align*}
 \tau(x)1(x\in S_2\cup S_3) \leq \log(x_2\sqrt{n}/\kappa_1) \leq 1 + \frac{\sqrt{n} x_2}{\kappa_1}   = 1 + \frac{\sqrt{n} x_2}{ \beta}.
\end{align*} 
We conclude by combining all of these bounds with \eqref{eq:f2ineq}.
\finishproof
 
\subsection{Proof of Lemma~\ref{lem:lem16_3}} 
\label{sec:claimthree}
Assume  without loss of generality  that $Q_1(0) = n$ and $Q_2(0) = \theta_{2}$, because starting from $(q_1,\theta_{2}, q_3, \ldots, q_{b+1}) \in S_Q$, a state with $q_1 = n$ and $q_{2} =  \theta_{2}$ must be visited before  $\tau_{2}(2\theta_{2}) $, so 
\begin{align*}
  \min_{\substack{  q_2 = \theta_{2} \\  q \in S_Q}} \Prob_{q}(\tau_{1}(\theta_{1}) < \tau_{2}(2\theta_{2}))  \geq&\ \min_{\substack{ q_2 = \theta_{2},\  q_1 = n \\ q \in S_Q }} \Prob_{q}(\tau_{1}(\theta_{1}) < \tau_{2}(2\theta_{2})).
\end{align*}
We bound the right-hand side by relating it to the ruin probability in a certain gambler's ruin problem. Namely, we construct a random walk $\{\overline R(t) \}$ with $\overline R(0) = 0$ that satisfies 
\begin{align}
   \min_{\substack{ q_2 = \theta_{2},\  q_1 = n \\ q \in S_Q }} \Prob_{q}(\tau_{1}(\theta_{1}) < \tau_{2}(2\theta_{2}))\geq&\ \Prob\Big( \inf_{t \geq 0} \{ \overline R(t) = \lfloor  \gamma \sqrt{n}\rfloor \} > \inf_{t \geq 0 } \{ \overline R(t)  = -\lfloor  \sqrt{n}\beta/2 \rfloor \}\Big). \label{eq:randomwalkbound}
\end{align}
Jumps in the random walk  are governed by a Poisson process with rate $n\lambda + \theta_{1} -3\theta_{2}$, and the up-step and down-step probabilities are 
\begin{align}
\frac{n\lambda}{n\lambda + \theta_{1} -3\theta_{2}}, \quad \text{ and } \quad \frac{\theta_{1} -3\theta_{2}}{n\lambda + \theta_{1} -3\theta_{2}}, \label{eq:updownwant}
\end{align}
respectively. Note that we implicitly assume $n$ is large enough so that $\theta_{1} - 3 \theta_{2} > 0$.  The right-hand side in \eqref{eq:randomwalkbound} is therefore the ruin probability in a gambler's ruin problem with initial wealth $\lfloor \sqrt{n}\beta/2 \rfloor $ and opponent's wealth $ \lfloor \sqrt{n}\beta/2 \rfloor + \lfloor  \gamma \sqrt{n}\rfloor$. A formula for the ruin probability was given by equation (2.4) in Section XIV.2 of \cite{Fell1968}:
\begin{align*}
\Prob\Big( \inf_{t \geq 0} \{ \overline R(t) = \lfloor  \gamma \sqrt{n}\rfloor \} > \inf_{t \geq 0 } \{ \overline R(t)  = -\lfloor  \sqrt{n}\beta/2 \rfloor \}\Big)  = 1 - \frac{1- \big((\theta_{1} -3\theta_{2})/n\lambda \big)^{\lfloor  \sqrt{n}\beta/2 \rfloor}}{1-\big((\theta_{1} -3\theta_{2})/n\lambda \big)^{\lfloor  \sqrt{n}\beta/2 \rfloor + \lfloor  \gamma\sqrt{n} \rfloor}}.
\end{align*}
Recalling the values of $\theta_1$ and $\theta_2$  and the fact that $\gamma > \beta$, we see that
\begin{align*}
\frac{\theta_{1} -3\theta_{2}}{n\lambda} = \frac{n -\lfloor \sqrt{n}\beta/2 \rfloor -3\lfloor \gamma \sqrt{n}\rfloor}{n\lambda} = 1 - \frac{-\beta \sqrt{n} +\lfloor \sqrt{n}\beta/2 \rfloor +3\lfloor \gamma \sqrt{n}\rfloor}{ n\lambda} < 1,
\end{align*} 
and therefore,
\begin{align*}
 \lim_{n \to \infty} \frac{1- \big((\theta_{1} -3\theta_{2})/n\lambda \big)^{\lfloor  \sqrt{n}\beta/2 \rfloor}}{1-\big((\theta_{1} -3\theta_{2})/n\lambda \big)^{\lfloor  \sqrt{n}\beta/2 \rfloor + \lfloor  \gamma\sqrt{n} \rfloor}} =  \lim_{n \to \infty} \frac{1- \Big(1 - \frac{-\beta \sqrt{n} +\lfloor \sqrt{n}\beta/2 \rfloor +3\lfloor \gamma \sqrt{n}\rfloor}{ n\lambda} \Big)^{\lfloor  \sqrt{n}\beta/2 \rfloor}}{1-\Big(1 - \frac{-\beta \sqrt{n} +\lfloor \sqrt{n}\beta/2 \rfloor +3\lfloor \gamma \sqrt{n}\rfloor}{ n\lambda} \Big)^{\lfloor  \sqrt{n}\beta/2 \rfloor + \lfloor  \gamma\sqrt{n} \rfloor}} < 1,
\end{align*}
implying Lemma~\ref{lem:lem16_3}. It remains to construct $\{ \overline R(t)\}$. 

Recall that $Q_1(0) = n$ and $Q_2(0) = \theta_{2}$, and let $\{\widehat Q(t)\}$ be a copy of $\{Q(t)\}$, but with the modification that any server with a nonempty buffer permanently halts all its their work. Then $\widehat Q_i(t) \geq Q_i(t)$ for all $t \geq 0$ and all $1 \leq i \leq b+1$, because this modified system has the same arrival stream as $\{Q(t)\}$  but  serves fewer customers. It follows that
\begin{align*}
 \tau_{1}(\theta_{1}) = \inf_{t \geq 0} \{  Q_1(t)  =\theta_{1} \} \leq&\  \inf_{t \geq 0} \{ \widehat Q_1(t)  =\theta_{1} \}, \\
 \tau_{2}(2\theta_{2})  =  \inf_{t \geq 0} \{ Q_2(t) = 2\theta_{2} \} \geq&\  \inf_{t \geq 0} \{ \widehat Q_2(t) = 2\theta_{2} \}, \text{ and }  \\
 \min_{\substack{ q_2 = \theta_{2},\  q_1 = n \\ q \in S_Q }} \Prob_{q}(\tau_{1}(\theta_{1}) < \tau_{2}(2\theta_{2})) 
\geq&\ \min_{\substack{ q_2 = \theta_{2},\  q_1 = n \\ q \in S_Q }}\Prob_{q}\Big( \inf_{t \geq 0} \{ \widehat Q_2(t) = 2\theta_{2}  \} > \inf_{t \geq 0} \{ \widehat Q_1(t)  = \theta_{1}  \}\Big).
\end{align*}
Now consider the process 
\begin{align*}
 R(t) =&\ \widehat Q_1(t) +  \widehat Q_2(t) - \widehat Q_1(0) - \widehat Q_2(0)= \widehat Q_1(t) +  \widehat Q_2(t) - (n+\theta_{2} ).
\end{align*}
Note that $R(t) \geq \widehat Q_1(t)  - \widehat Q_1(0)  =  \widehat Q_1(t)  - n$ since $\widehat Q_2(t)$ is non-decreasing in $t$,  which implies that 
\begin{align*}
\inf_{t \geq 0} \{ R(t)  = -\lfloor  \sqrt{n}\beta/2 \rfloor \} \geq \inf_{t \geq 0} \{ \widehat Q_1(t)  = n-\lfloor  \sqrt{n}\beta/2 \rfloor \} = \inf_{t \geq 0} \{ \widehat Q_1(t)  = \theta_{1} \}.
\end{align*}
Note also that $\inf_{t \geq 0} \{ \widehat Q_2(t) = 2\theta_{2} \} = \inf_{t \geq 0} \{ R(t) = \theta_{2}\}$ because  $\widehat Q_2(t)$ is non-decreasing in $t$ and   $\widehat Q_2(t)$  increases only when $\widehat Q_1(t) = n$. Hence,
\begin{align}
&\min_{\substack{ q_2 = \theta_{2},\  q_1 = n \\ q \in S_Q }} \Prob_{q}\Big( \inf_{t \geq 0} \{ R(t) = \theta_{2} \} > \inf_{t \geq 0} \{ R(t)  = -\lfloor  \sqrt{n}\beta/2 \rfloor \}\Big) \notag \\ 
\leq&\ \min_{\substack{ q_2 = \theta_{2},\  q_1 = n \\ q \in S_Q }} \Prob_{q}\Big( \inf_{t \geq 0} \{ \widehat Q_2(t) = 2\theta_{2} \} > \inf_{t \geq 0} \{ \widehat Q_1(t)  = \theta_{1} \}\Big)  \leq \min_{\substack{ q_2 = \theta_{2},\  q_1 = n \\ q \in S_Q }} \Prob_{q}(\tau_{1}(\theta_{1}) < \tau_{2}(2\theta_{2})). \label{eq:intermpbound}
\end{align}
 An arrival to $\{\widehat Q(t)\}$  increases the value of $\{R(t)\}$, and a service completion by a server with an empty buffer decreases its value. However, $\{R(t)\}$ is still not the random walk we desire  because the rate at which it decreases depends on the state of $\widehat Q(t)$. Instead, we want a random walk with a constant downward rate. 
 
 To construct this random walk,  for $0 \leq t \leq \inf_{t \geq 0} \{ \widehat Q_2(t) = 2\theta_{2}\} $ let us define $\big\{\overline Q(t) = (\overline Q_1(t), \overline Q_2(t))   \big\}$ by setting $\overline Q(0) = \widehat Q(0)$  and defining the transitions of the joint process  $\big\{ \big(\widehat Q(t), \overline Q(t) \big)  \big\}$ 
 in Tables~\ref{tab:arriv}, \ref{tab:depart1}, and \ref{tab:depart2} below. Since we are defining $\overline Q(t)$ only until the time $\widehat Q_2(t)$ hits $2\theta_{2}$, we do not need to specify the transitions for states where $\widehat Q_2(t) >2\theta_{2}$. The intuition for the transition structure is as follows. Since arrivals occur at the constant rate of $n \lambda$, we want any arrival to $\{\widehat Q(t)\}$ to also occur in $\{\overline Q(t)\}$. However, we want to keep the rate at which $\{\overline Q(t)\}$ decreases a constant value of $\theta_{1} -3\theta_{2}$. To accomplish this,  when $\widehat Q_1(t) \geq \theta_1 - \theta_2$, the transitions in Table~\ref{tab:depart1} have $\{\overline Q(t)\}$ ignore  some departures from $\{\widehat Q(t)\}$, and when  $\widehat Q_1(t) < \theta_1 - \theta_2$,  we supplement the departures from $\{\widehat Q(t)\}$; e.g., see transition $\# 8$ in Table~\ref{tab:depart2}. 
\begin{table}[h!]
\caption{Arrival transitions for the joint process in  state $ \big( (\widehat u_1, \widehat u_2), (\overline u_1, \overline u_2) \big)$. }
\centering
\label{tab:arriv}
\begin{tabular}{|c|c|c|}
\hline
\# &Rate & Transition      \\ \hline
1& $n\lambda  1(\widehat u_1 < n, \overline u_1 < n)  $ & $\big( (\widehat u_1+1, \widehat u_2), (\overline u_1+1, \overline u_2) \big)$  \\ \hline
2&$n\lambda  1(\widehat u_1 = n, \overline u_1 < n)  $ & $\big( (\widehat u_1, \widehat u_2+1), (\overline u_1+1, \overline u_2) \big)$  \\ \hline
3&$n\lambda  1(\widehat u_1 < n, \overline u_1 = n)  $ & $\big( (\widehat u_1+1, \widehat u_2), (\overline u_1, \overline u_2+1) \big)$  \\ \hline
4&$n\lambda  1(\widehat u_1 = n, \overline u_1 = n)  $ & $\big( (\widehat u_1, \widehat u_2+1), (\overline u_1, \overline u_2+1) \big)$  \\ \hline
\end{tabular}
\end{table}
\begin{table}[h!]
\caption{Departure transitions for the joint process in  state $ \big( (\widehat u_1, \widehat u_2), (\overline u_1, \overline u_2) \big)$ with $\widehat u_2 \leq 2\theta_2$ and $\widehat u_1 \geq \theta_{1} - \theta_{2}$. }
\centering
\label{tab:depart1}
\begin{tabular}{|c|c|c|}
\hline
\# &Rate & Transition      \\ \hline
5& $ \theta_{1} -3\theta_{2}  $  & $\big( (\widehat u_1-1, \widehat u_2), (\overline u_1-1, \overline u_2) \big)$ \\ \hline
6& $ \widehat u_1 - \widehat u_2 - (\theta_{1} -3\theta_{2})  $  & $\big( (\widehat u_1-1, \widehat u_2), (\overline u_1, \overline u_2) \big)$ \\ \hline
\end{tabular}
\end{table}
\begin{table}[h!]
\caption{Departure transitions for the joint process in  state $ \big( (\widehat u_1, \widehat u_2), (\overline u_1, \overline u_2) \big)$ with $\widehat u_2 \leq 2\theta_{2}$ and $\widehat u_1 < \theta_{1} - \theta_{2}$. }
\centering
\label{tab:depart2}
\begin{tabular}{|c|c|c|}
\hline
\# &Rate & Transition      \\ \hline
 7&$ (\widehat u_1  - 2\theta_{2} ) 1(\widehat u_1 \geq 2\theta_{2}) $  & $\big( (\widehat u_1-1, \widehat u_2), (\overline u_1-1, \overline u_2) \big)$ \\  \hline
8& $\theta_{1} -\theta_{2}-  \widehat u_1 \vee  2\theta_{2}  $  & $\big( (\widehat u_1, \widehat u_2), (\overline u_1-1, \overline u_2) \big)$ \\  \hline
9& $ 2\theta_{2} \wedge \widehat u_1 - \widehat u_2  $  & $\big( (\widehat u_1-1, \widehat u_2), (\overline u_1, \overline u_2) \big)$ \\  \hline
\end{tabular}
\end{table}
\newpage
\noindent  Having defined $\overline Q(t)$, let us define 
\begin{align*}
\overline R(t) = \overline Q_1(t) - \overline Q_1(0) + \overline Q_2(t) - \overline Q_2(0), \quad t \leq \inf_{t \geq 0} \{ \widehat Q_2(t) = 2\theta_2 \}.
\end{align*}
To prove that $\{\overline R(t)\}$ satisfies \eqref{eq:randomwalkbound}, we show that
\begin{align}
\overline R(t) \geq R(t) \text{ for all times } t \leq \min \Big\{ \inf_{t \geq 0} \{ \overline R(t) = \lfloor  \gamma \sqrt{n}\rfloor \},\ \inf_{t \geq 0} \{ \overline R(t)  = -\lfloor  \sqrt{n}\beta/2 \rfloor \} \Big\}, \label{eq:tildebigger}
\end{align}
and as a result, 
\begin{align*}
\inf_{t \geq 0} \{ \overline R(t) = \lfloor  \gamma \sqrt{n}\rfloor \} \leq&\ \inf_{t \geq 0} \{  R(t) = \lfloor  \gamma \sqrt{n}\rfloor \},\\
\inf_{t \geq 0} \{ \overline R(t)  = -\lfloor  \sqrt{n}\beta/2 \rfloor \} \geq&\ \inf_{t \geq 0} \{ R(t)  = -\lfloor  \sqrt{n}\beta/2 \rfloor \}.
\end{align*}
Together with \eqref{eq:intermpbound}, these inequalities imply that 
\begin{align*}
&\min_{\substack{ q_2 = \theta_{2},\  q_1 = n \\ q \in S_Q }} \Prob_{q}\Big( \inf_{t \geq 0} \{ \overline R(t) = \lfloor  \gamma \sqrt{n}\rfloor \} > \inf_{t \geq 0} \{ \overline R(t)  = -\lfloor  \sqrt{n}\beta/2 \rfloor \}\Big) \\
\leq&\ \min_{\substack{ q_2 = \theta_{2},\  q_1 = n \\ q \in S_Q }} \Prob_{q}\Big( \inf_{t \geq 0} \{ R(t) = \lfloor  \gamma \sqrt{n}\rfloor \} > \inf_{t \geq 0} \{ R(t)  = -\lfloor  \sqrt{n}\beta/2 \rfloor \}\Big)  \leq \min_{\substack{ q_2 = \theta_{2},\  q_1 = n \\ q \in S_Q }} \Prob_{q}(\tau_{1}(\theta_{1}) < \tau_{2}(2\theta_{2})).
\end{align*}
To see why \eqref{eq:tildebigger} is true, let us study the transitions in Tables~\ref{tab:arriv}--\ref{tab:depart2}.  Table~\ref{tab:arriv} tells us that  $\overline R(t)$ and $R(t)$ increase at the same times. The transitions in Table~\ref{tab:depart1} show that any decrease in $\overline Q_1(t)$, and consequently  $\overline R(t)$, must be accompanied by a decrease in $\widehat Q_1(t)$ and  $ R(t)$, but not vice versa.  The only way $\overline Q_1(t)$ can ever drop below $\widehat Q_1(t)$ is via transition 8, which can happen only if $\widehat Q_1(t) < \theta_{1} - \theta_{2}$, so the first intersection of $\overline Q_1(t)$ and $\widehat Q_1(t)$ has to occur below $\theta_1-\theta_2$. Therefore, $\overline R(t) \geq R(t)$ for all times 
\begin{align*}
\quad t \leq&\ \min \Big\{ \inf_{ t\geq 0} \{ \overline Q_1(t) = \theta_{1} - \theta_{2} \}, \inf_{t \geq 0} \{\widehat Q_2(t) =2 \theta_{2} \} \Big\} = \min \Big\{ \inf_{ t\geq 0} \{ \overline Q_1(t) = \theta_{1} - \theta_{2}  \}, \inf_{t \geq 0} \{R(t) =  \theta_{2}\} \Big\}.
\end{align*}
Let us now prove \eqref{eq:tildebigger} by showing that the right-hand side is greater than 
\begin{align*}
\min \Big\{ \inf_{t \geq 0} \{ \overline R(t) = \lfloor  \gamma \sqrt{n}\rfloor \},\ \inf_{t \geq 0} \{ \overline R(t)  = -\lfloor  \sqrt{n}\beta/2 \rfloor \} \Big\}.
\end{align*}
Since $\overline R(t) \geq R(t)$, 
\begin{align*}
 \min \Big\{ \inf_{ t\geq 0} \{ \overline Q_1(t) =\theta_{1} - \theta_{2} \}, \inf_{t \geq 0} \{R(t) =  \theta_{2}\} \Big\}   \geq&\  \min \Big\{ \inf_{ t\geq 0} \{ \overline Q_1(t) = \theta_{1} - \theta_{2} \}, \inf_{t \geq 0} \{\overline R(t) =  \theta_{2}\} \Big\}.
\end{align*}
Furthermore, since $\overline Q_2(t)$ is non-decreasing and increases only at those times when $\overline Q_1(t) = n$, it follows that   for all $t \leq \inf_{t \geq 0} \{\overline R(t) =\theta_{2} \}$,
\begin{align*}
\overline R(t) =&\ \overline Q_1(t) + \overline Q_2(t) - n - \theta_{2} \leq  \overline Q_1(t) - n + \theta_{2},
\end{align*}
and therefore 
\begin{align*}
 &  \min \Big\{ \inf_{ t\geq 0} \{ \overline Q_1(t) = \theta_{1} - \theta_{2} \}, \inf_{t \geq 0} \{\overline R(t) = \theta_{2} \} \Big\}  \\
 =&\ \min \Big\{ \inf_{ t\geq 0} \{ \overline Q_1(t) = n  -\lfloor \sqrt{n}\beta/2 \rfloor - \theta_{2} \}, \inf_{t \geq 0} \{\overline R(t) = \lfloor  \gamma \sqrt{n}\rfloor \} \Big\}  \\
   \geq&\ \min \Big\{ \inf_{ t\geq 0} \{ \overline R(t) = -\lfloor \sqrt{n}\beta/2 \rfloor \}, \inf_{t \geq 0} \{\overline R(t) = \lfloor  \gamma \sqrt{n}\rfloor \} \Big\}.
\end{align*} 
\finishproof

\subsection{ Proving Lemma~\ref{lem:lem16_4}}  
\label{sec:claimfour}
Central to our argument is a result about the moment-generating function of the duration of a gambler's ruin game. We now describe this result  and then prove Lemma~\ref{lem:lem16_4}. Consider a \textit{discrete-time} gambler's ruin problem where the initial player's wealth is $z$, the win probability is $p$, the loss probability is $q$, and the player keeps playing until they go broke or accumulate a total wealth of $a$. Let $D_z \in \Z_+$ be the number of turns until the game ends, given an  initial wealth of $z$. An expression for the generating function $\E s^{D_z}$ was given in (4.11) and (4.12) in Section XIV.4 of \cite{Fell1968}:
\begin{align}
\E s^{D_z} =  \frac{\lambda_1^{a}(s)\lambda_2^{z}(s) - \lambda_1^{z}(s)\lambda_2^{a}(s)}{\lambda_1^{a}(s) - \lambda_2^{a}(s)} + \frac{\lambda_1^{z}(s) - \lambda_2^{z}(s)}{\lambda_1^{a}(s) - \lambda_2^{a}(s)}, \quad s \in (0,1), \label{eq:genfunctionDz}
\end{align}
where
\begin{align*}
\lambda_1(s) = \frac{1+\sqrt{1-4pqs^2}}{2ps}, \quad \text{ and } \quad \lambda_2(s) = \frac{1-\sqrt{1-4pqs^2}}{2ps}, \quad s \in (0,1). 
\end{align*}
Now consider  the \textit{continuous-time} gambler's ruin problem, where  the durations between turns are governed by an i.i.d.\ sequence $\{E_i\}$ of rate $r$ exponentially distributed random variables. Given initial wealth $z$, the duration of the continuous game equals $\sum_{i=1}^{D_z} E_i$. Since the $E_i$ are independent of $D_z$, it follows that 
\begin{align*}
\E e^{-\sum_{i=1}^{D_z} E_i} = \E  \big(\E e^{- E_1}\big)^{D_z}    = \E \Big( \frac{r}{r+1} \Big)^{D_z},
\end{align*}
 so $\E e^{-\sum_{i=1}^{D_z} E_i}$ is related to \eqref{eq:genfunctionDz}.  The following result  proved in Appendix~\ref{sec:gambler} is needed to   prove Lemma~\ref{lem:lem16_4}.
\begin{lemma}
\label{lem:gamblersruin}
Let $i$ and $q_2$ be integers such that $1 \leq i \leq b+1$ and $0 \leq q_2 \leq \theta_{2}$,  and define 
\begin{align*}
q^{(B,i)} =  n- q_2 - 1- \lfloor \sqrt{n}\beta/2 \rfloor + \big \lfloor \lfloor \sqrt{n}\beta/2 \rfloor (b+1) \big \rfloor.
\end{align*}
 Consider the continuous-time gambler's ruin problem with probabilities
\begin{align*}
p =&\ \frac{n\lambda}{n\lambda + q^{(B,i)} -\lfloor \sqrt{n}\beta/2 \rfloor},  \quad \text{ and } \quad  q = \frac{ q^{(B,i)} -\lfloor \sqrt{n}\beta/2 \rfloor}{n\lambda +  q^{(B,i)} -\lfloor \sqrt{n}\beta/2 \rfloor},
\end{align*} 
rate $r =  n\lambda + q^{(B,i)} -\lfloor \sqrt{n}\beta/2 \rfloor$, initial wealth $z$ and terminal wealth $a$ given by
\begin{align}
z =&\ \lfloor \sqrt{n}\beta/2 \rfloor, \quad \text{ and } \quad a = \lfloor \sqrt{n}\beta/2 \rfloor + \big \lfloor \lfloor \sqrt{n}\beta/2 \rfloor (b+1) \big \rfloor, \label{eq:gamblerwealth}
\end{align}
and game duration $\sum_{i=1}^{D_z} E_i$. Then 
\begin{align}
\lim_{n \to \infty} \max_{\substack{ 0 \leq q_2  \leq 2 \lfloor \gamma\sqrt{n}\rfloor }} \E e^{-\sum_{i=1}^{D_z} E_i} < 1. \label{eq:gamblermgf}
\end{align}
\end{lemma}
\startproof{Proof of Lemma~\ref{lem:lem16_4}}
As discussed below \eqref{eq:bbound},   $\{\tau_C < \tau_1(n)\} \supset \{\Gamma_{b+1} < \tau_1(n)\} $, where   $\Gamma_{b+1}$ is the sum of $b+1$ unit-mean exponentially distributed random variables. The same discussion says that $\Gamma_{b+1}$ represents the time needed by the joint CTMC  $(Q(t),\widetilde Q(t))$ to transition from   $\Theta_{b+1}^{Q}$ to $\Theta_{1}^{Q}$, and to then couple by spending an exponentially distributed amount of time in $\Theta_{1}^{Q}$. Thus, 
\begin{align*}
\min_{\substack{ 0 \leq q_1 \leq \theta_{1} \\  0 \leq q_2  \leq 2 \theta_{2}   \\ q \in S_Q}} \Prob\Big(\tau_C  < \tau_{1}(n)\ \big|\ Q(0) = q,\ (Q(0),\widetilde Q(0)) \in \bigcup_{i=1}^{b+1}\Theta_{i}^{Q}\Big) \geq&\ \min_{\substack{ 0 \leq q_1 \leq \theta_{1} \\  0 \leq q_2  \leq 2 \theta_{2}   \\ q \in S_Q}} \Prob\Big( \Gamma_{b+1} < \tau_{1}(n)\ \big|\ Q(0) = q\Big).
\end{align*} 
Let us analyze the probability above.  At time $t=0$, there are $q_2$ servers with nonempty buffers and another server containing the extra customer in $\{\widetilde Q(t)\}$. We group these $q_2+1$ servers together into group $A$, and the remaining $n-q_2-1$ servers into group $B$. Let $Q^{(A)}_1(t)$ and $Q^{(B)}_1(t)$ be the number of busy group $A$ and $B$ servers, respectively. Since 
\begin{align*}
Q_1^{(A)}(0) =&\ q_2 +1, \quad \text{ and } \quad   Q_1^{(A)}(0)+Q_1^{(B)}(0) = Q_1(0)  \leq  n-\lfloor \sqrt{n}\beta/2 \rfloor,
\end{align*}
it follows that $Q_1^{(B)}(0) \leq  n - q_2 - 1-\lfloor \sqrt{n}\beta/2 \rfloor $. We are implicitly assuming that $n$ is large enough so   $ n - q_2 - 1-\lfloor \sqrt{n}\beta/2 \rfloor \geq 0$.  Note that the buffer of any group $B$ server is empty for all $t \leq \tau_1(n)$.   

%Recall that when introducing the JSQ model in the beginning of this paper, we allowed for arbitrary tie-breaking decisions when choosing where to send the next customer in a system where two or more servers were tied for the smallest number of customers.
If a customer arrives when more than one server is idle, we prioritize assigning this customer to servers in group B over group A. Note that this tie-breaking rule is consistent with the tie-breaking rule we imposed in the proof of Lemma~\ref{lem:coupling}.  Let  $\tau_B = \inf_{t \geq 0} \{Q_1^{(B)}(t) = n-q_2-1 \}$   be the first time that all servers in group B are busy. By construction, $\tau_B \leq \tau_1(n)$, so 
\begin{align*}
 \min_{\substack{ 0 \leq q_1 \leq \theta_{1} \\  0 \leq q_2  \leq 2 \theta_{2}   \\ q \in S_Q}}  \Prob\Big( \Gamma_{b+1} < \tau_{1}(n)\ \big|\ Q(0) = q\Big) \geq&\   \min_{\substack{ 0 \leq q_2  \leq 2\theta_{2} \\  0 \leq q^{(B)} \leq n- q_2 - 1- \lfloor \sqrt{n}\beta/2 \rfloor }} \Prob\Big(\Gamma_{b+1} < \tau_B\ \big|\ Q^{(B)}_1(0) = q^{(B)}\Big)\\
\geq&\ \min_{\substack{ 0 \leq q_2  \leq 2 \lfloor \gamma\sqrt{n}\rfloor }} \Prob\Big(\Gamma_{b+1} < \tau_B\ \big|\ Q^{(B)}_1(0) = n- q_2 - 1- \lfloor \sqrt{n}\beta/2 \rfloor\Big).
\end{align*}
The last inequality is true because increasing the value of the initial condition $Q^{(B)}_1(0)$ does not  increase the chance that $\Gamma_{b+1} < \tau_{B}$.  We now  relate the right-hand side to the moment-generating function considered in Lemma~\ref{lem:gamblersruin} and use that lemma to conclude the proof. We can write $\Gamma_{b+1} = \sum_{i=1}^{b+1} G_i$, where $G_i$ are i.i.d.\ unit-mean exponentially distributed random variables  independent of $Q_1^{(B)}(t)$ for $t \in [0,\tau_B]$, because they correspond to service times of the server containing the additional customer in $\{\widetilde Q(t)\}$, which is a server in group A. 

Fixing  $0 \leq q_2  \leq 2 \theta_{2} $ and $Q_1^{(B)}(0) = n- q_2 - 1- \lfloor \sqrt{n}\beta/2 \rfloor$,   for  $0 \leq i \leq b+1$ we define 
\begin{align*}
q^{(B,i)} =&\  n- q_2 - 1- \lfloor \sqrt{n}\beta/2 \rfloor + \Big \lfloor \lfloor \sqrt{n}\beta/2 \rfloor \frac{i}{b+1} \Big \rfloor,  \quad \text{ and } \\
 \tau_{B,i} =&\ \inf_{t \geq 0}   \Big\{Q_1^{(B)}(t) - Q_1^{(B)}(0)  = \Big \lfloor \lfloor \sqrt{n}\beta/2 \rfloor \frac{i}{b+1} \Big \rfloor \Big\} = \inf_{t \geq 0}   \Big\{Q_1^{(B)}(t)   = q^{(B,i)} \Big\},
\end{align*}
and note that $\tau_{B} = \tau_{B,b+1}$. We are guaranteed that $\Gamma_{b+1} < \tau_{B}$ if for each $1 \leq i \leq b+1$, the exponentially distributed $G_{i}$ is smaller than the time it takes for $Q_{1}^{(B)}(t)$ to reach $q^{(B,i)}$  if started from $q^{(B,i-1)}$,  so  
\begin{align*}
 \Prob\Big(\Gamma_{b+1} < \tau_B\ \big|\ Q^{(B)}_1(0) = n- q_2 - 1- \lfloor \sqrt{n}\beta/2 \rfloor\Big)      \geq&\  \prod_{i=1}^{b+1} \Prob\Big(G_{i}  < \tau_{B,i}  \ \big|\ Q^{(B)}_1(0) =q^{(B,i-1)}\Big).
\end{align*}
We now show that $\tau_{B,i}$ can be bounded from below by the duration of a   gambler's ruin game, which  allows us to    apply Lemma~\ref{lem:gamblersruin}. Fix $1 \leq i \leq b+1$, and consider the time interval $t \in [0,\tau_{B,i}]$, on which we construct the coupling  $\big\{ \big( Q_1^{(B)}(t), \overline Q_1^{(B)}(t)\big) \big\}$  by setting 
\begin{align*}
\overline Q^{(B,i)}_1(0) = Q^{(B)}_1(0) = q^{(B,i-1)} 
\end{align*}
and defining the transitions of the joint process   in Tables~\ref{tab:transitions_random_walk} and \ref{tab:transitions_random_walk2} below. We implicitly assume that $n$ is large enough  that $q^{(B,i)} - \lfloor \sqrt{n}\beta/2 \rfloor > 0$.  
\begin{table}[h!]
\caption{Transition rates in state $(u,\overline u)$ with $u \geq q^{(B,i-1)} - \lfloor \sqrt{n}\beta/2 \rfloor$.}
\centering
\label{tab:transitions_random_walk}
\begin{tabular}{|c|c|}
\hline
Rate & Transition      \\ \hline
$n\lambda   $ & $(u+1,\overline u+1)$  \\ \hline
 $q^{(B,i-1)} -\lfloor \sqrt{n}\beta/2 \rfloor$ & $(u-1,\overline u-1)$ \\ \hline
 $u - (q^{(B,i-1)} -\lfloor \sqrt{n}\beta/2 \rfloor)$ & $(u-1,\overline u)$  \\    \hline
\end{tabular}
\end{table}
\begin{table}[h!]
\caption{Transition rates in state $(u,\overline u)$ with $u < q^{(B,i-1)} -\lfloor \sqrt{n}\beta/2 \rfloor$.}
\centering
\label{tab:transitions_random_walk2}
\begin{tabular}{|c|c|}
\hline
Rate & Transition      \\ \hline
 $u  $ & $(u-1,\overline u-1)$  \\ \hline
$q^{(B,i-1)} -\lfloor \sqrt{n}\beta/2 \rfloor-u$ & $(u,\overline u-1)$ \\ \hline
\end{tabular}
\end{table}

\noindent Note that the only time  $\overline Q_1^{(B,i)}(t)$ decreases but  $ Q_1^{(B)}(t)$ does not is when the latter is smaller than $q^{(B,i-1)} -\lfloor \sqrt{n}\beta/2 \rfloor$,  so we are guaranteed that  
\begin{align}
\overline Q_1^{(B,i)}(t) \geq Q_1^{(B)}(t), \quad \text{ for all } t \leq \min \Big\{ \tau_{B,i},\ \inf_{t \geq 0} \big\{\overline Q_1^{(B,i)}(t) = q^{(B,i-1)}- \lfloor \sqrt{n}\beta/2 \rfloor   \big\} \Big\}. \label{eq:forallt} 
\end{align}
Recalling the definitions of $\tau_{B,i}$ and $q^{(B,i)}$, we have 
\begin{align*}
& \min \Big\{ \tau_{B,i},\ \inf_{t \geq 0} \big\{\overline Q_1^{(B,i)}(t) = q^{(B,i-1)}- \lfloor \sqrt{n}\beta/2 \rfloor   \big\} \Big\} \\
=&\ \min \Big\{ \inf_{t \geq 0} \big\{Q_1^{(B)}(t) = q^{(B,i)} \big\} ,\ \inf_{t \geq 0} \big\{\overline Q_1^{(B,i)}(t) = q^{(B,i-1)}- \lfloor \sqrt{n}\beta/2 \rfloor   \big\} \Big\}\\
=&\ \min \Big\{ \inf_{t \geq 0} \big\{Q_1^{(B)}(t) = q^{(B,i-1)}+ \big \lfloor \lfloor \sqrt{n}\beta/2 \rfloor /(b+1)\big \rfloor \big\} ,\ \inf_{t \geq 0} \big\{\overline Q_1^{(B,i)}(t) = q^{(B,i-1)}- \lfloor \sqrt{n}\beta/2 \rfloor   \big\} \Big\}\\
\geq&\ \min \Big\{ \inf_{t \geq 0} \big\{\overline Q_1^{(B,i)}(t) = q^{(B,i-1)}+ \big \lfloor \lfloor \sqrt{n}\beta/2 \rfloor /(b+1)\big \rfloor \big\} ,\ \inf_{t \geq 0} \big\{\overline Q_1^{(B,i)}(t) = q^{(B,i-1)}- \lfloor \sqrt{n}\beta/2 \rfloor   \big\} \Big\},
\end{align*}
where the last inequality follows from \eqref{eq:forallt}. Let $\overline \tau_{B,i}$ equal the right-hand side and note that  
\begin{align*}
\overline \tau_{B,i} = \inf_{t \geq 0} \Big\{\big(\overline Q_1^{(B,i)}(t) - \overline Q_1^{(B,i)}(0)\big) \in \big\{  - \lfloor \sqrt{n}\beta/2 \rfloor ,  \big \lfloor \lfloor \sqrt{n}\beta/2 \rfloor / (b+1)\big \rfloor  \big\}  \Big\} 
\end{align*}
because $\overline Q^{(B,i)}_1(0)  = q^{(B,i-1)}$. Since  $\overline \tau_{B,i} \leq \tau_{B,i}$, it follows that 
\begin{align*}
\min_{\substack{ 0 \leq q_2  \leq 2 \lfloor \gamma\sqrt{n}\rfloor }} \prod_{i=1}^{b+1} \Prob\Big(G_{i}  < \tau_{B,i}  \ \big|\ Q^{(B)}_1(0) =q^{(B,i-1)}\Big) \geq  \min_{\substack{ 0 \leq q_2  \leq 2 \lfloor \gamma\sqrt{n}\rfloor }} \prod_{i=1}^{b+1} \Prob\Big(G_{i}  < \overline \tau_{B,i}  \ \big|\ Q^{(B)}_1(0) =q^{(B,i-1)}\Big).
\end{align*}
Recall that $G_{i}$ corresponds to the service time of a group A server  and is therefore independent of   $\overline \tau_{B,i}$. Furthermore, since $G_i$ is exponentially distributed with unit mean, conditioning on the value of $\overline \tau_{B,i}$ yields  
\begin{align*}
 \min_{\substack{ 0 \leq q_2  \leq 2 \lfloor \gamma\sqrt{n}\rfloor }} \Prob\Big(G_{i}  < \overline \tau_{B,i}  \ \big|\ \overline Q^{(B,i)}_1(0) = q^{(B,i-1)}\Big) =&\    \min_{\substack{ 0 \leq q_2  \leq 2 \lfloor \gamma\sqrt{n}\rfloor }} \bigg( 1 -  \E\Big(  e^{-\overline \tau_{B,i}} \big|\ \overline Q^{(B,i)}_1(0) = q^{(B,i-1)} \Big) \bigg) \\ 
=&\     1 -  \max_{\substack{ 0 \leq q_2  \leq 2 \lfloor \gamma\sqrt{n}\rfloor }} \E\Big(  e^{-\overline \tau_{B,i}} \big|\ \overline Q^{(B,i)}_1(0) = q^{(B,i-1)} \Big). 
\end{align*}
Applying \eqref{eq:gamblermgf} of Lemma~\ref{lem:gamblersruin} concludes, because our construction of  $\{\overline Q^{(B,i)}_1(t)\}$ implies that $\overline \tau_{B,i}$ is the duration of a gambler's ruin game with initial wealth $z = \lfloor \sqrt{n}\beta/2 \rfloor $,   terminal wealth  $a = \lfloor \sqrt{n}\beta/2 \rfloor + \big \lfloor \lfloor \sqrt{n}\beta/2 \rfloor /(b+1) \big \rfloor$,   rate  $n\lambda + q^{(B,i)} -\lfloor \sqrt{n}\beta/2 \rfloor$,  and up-step and down-step probabilities 
\begin{align*}
\frac{n\lambda}{n\lambda + q^{(B,i)} -\lfloor \sqrt{n}\beta/2 \rfloor}  \quad \text{ and }  \quad \frac{q^{(B,i)} -\lfloor \sqrt{n}\beta/2 \rfloor}{n\lambda + q^{(B,i)} -\lfloor \sqrt{n}\beta/2 \rfloor}. 
\end{align*} 
\finishproof

\subsubsection{Proving the Gambler's Ruin Result.}
\label{sec:gambler}
We require the following auxiliary lemma. 
 \begin{lemma}
 \label{lem:exponential_aux}
 Assume $\{x_n \in \R \}$ is a sequence that converges to $\overline x$. Then 
 \begin{align*}
 \lim_{n\to \infty} \Big(1 + \frac{x_n}{n} \Big)^{n} \to e^{\overline x}.
 \end{align*}
 \end{lemma}
\startproof{Proof of Lemma~\ref{lem:exponential_aux}}
Let $f(x) = e^{x}$ and $f_n(x) = \Big(1 + \frac{x}{n} \Big)^{n}$, and note that for any $n \geq 0$, 
 \begin{align*}
 \abs{f_n(x_n) - e^{\overline x}} \leq \abs{f_n(x_n) - f_n(\overline x)} + \abs{f_n(\overline x) - e^{\overline x}}.
 \end{align*}
From the mean-value theorem, we know that there exists some $c_n$ between $x_n$ and $\overline x$ such that 
 \begin{align*}
 \abs{f_n(x_n) - f_n(\overline x)} \leq \abs{x_n - \overline x} f_n'(c_n) =  \abs{x_n - \overline x}\Big(1 + \frac{c_n}{n} \Big)^{n-1}.
 \end{align*}
 Since $x_n \to \overline x$, it follows that  $\big(1 + c_n/n  \big)^{n-1} \leq \big(1 + 2\abs{\overline x}/n \big)^{n-1}$ for $n$ large enough, and therefore,
  \begin{align*}
 \abs{f_n(x_n) - e^{\overline x}} \leq \abs{x_n - \overline x}\Big(1 + \frac{2\abs{\overline x}}{n} \Big)^{n-1} + \abs{f_n(\overline x) - e^{\overline x}}.
 \end{align*}
 We can make the right-hand side arbitrarily small by increasing $n$. 
\finishproof 

\startproof{Proof of Lemma~\ref{lem:gamblersruin}}
Recall that  $\E e^{-\sum_{i=1}^{D_z} E_i}  = \E (r/(r+1))^{D_z}$, and that 
\begin{align*}
 \E s^{D_z} =  \frac{\lambda_1^{a}(s)\lambda_2^{z}(s) - \lambda_1^{z}(s)\lambda_2^{a}(s)}{\lambda_1^{a}(s) - \lambda_2^{a}(s)} + \frac{\lambda_1^{z}(s) - \lambda_2^{z}(s)}{\lambda_1^{a}(s) - \lambda_2^{a}(s)} = \frac{\lambda_2^{z}(s)(\lambda_1^{a}(s)-1) - \lambda_1^{z}(s)(\lambda_2^{a}(s)-1)}{\lambda_1^{a}(s) - \lambda_2^{a}(s)},
\end{align*}
where
\begin{align*}
\lambda_1(s) = \frac{1+\sqrt{1-4pqs^2}}{2ps}  \quad \text{ and } \quad \lambda_2(s) = \frac{1-\sqrt{1-4pqs^2}}{2ps}, \quad s \in (0,1).
\end{align*}
 Fix $s = r/(r+1)$. To show that $\lim_{n \to \infty}  \E s^{D_z} < 1$,  we  derive expressions for   $\lim_{n\to\infty}\lambda^{z}_j(s)$ and $\lim_{n\to\infty}\lambda^{a}_j(s)$. For notational economy, we let $\theta_{3} = \lfloor \sqrt{n}\beta/2 \rfloor$. We can write $p$ and $q$   as 
%\begin{align*}
%p =&\ \frac{n\lambda}{n\lambda + q^{(B,i)} -\lfloor \sqrt{n}\beta/2 \rfloor},  \quad \text{ and } \quad  q = \frac{ q^{(B,i)} -\lfloor \sqrt{n}\beta/2 \rfloor}{n\lambda +  q^{(B,i)} -\lfloor \sqrt{n}\beta/2 \rfloor},
%\end{align*}
\begin{align*}
p = \frac{n\lambda}{n\lambda + q^{(B,i)} -\theta_{3}} =&\ \frac{1}{2} + \frac{1}{2} \frac{n\lambda - (q^{(B,i)} -\theta_{3})}{n\lambda + q^{(B,i)} -\theta_{3}}, \qquad   q =  \frac{1}{2} - \frac{1}{2} \frac{n\lambda - (q^{(B,i)} -\theta_{3})}{n\lambda + q^{(B,i)} -\theta_{3}},
\end{align*} 
and 
\begin{align*}
pq =&\ \frac{1}{4} - \frac{1}{4} \Big( \frac{n\lambda - (q^{(B,i)} -\theta_{3})}{n\lambda + q^{(B,i)} -\theta_{3}} \Big)^2.
\end{align*}
  Let us first consider  $\lambda_1(s)$, which satisfies
\begin{align}
\lambda_1(s) =&\ \Bigg(1+\sqrt{1-s^2 + \Big( \frac{n\lambda - (q^{(B,i)} -\theta_{3})}{n\lambda + q^{(B,i)} -\theta_{3}} \Big)^2 s^2}\Bigg)  s^{-1}\Bigg(1+\frac{n\lambda - (q^{(B,i)} -\theta_{3})}{n\lambda + q^{(B,i)} -\theta_{3}}\Bigg)^{-1} \notag \\
=&\ \Bigg(1+\frac{1}{\sqrt{n}} \Bigg[ \sqrt{n(1-s^2) + n\Big( \frac{n\lambda - (q^{(B,i)} -\theta_{3})}{n\lambda + q^{(B,i)} -\theta_{3}} \Big)^2 s^2}\Bigg]\Bigg)  s^{-1}\Bigg(1+\frac{1}{\sqrt{n}} \bigg[ \sqrt{n}\frac{n\lambda - (q^{(B,i)} -\theta_{3})}{n\lambda + q^{(B,i)} -\theta_{3}}\bigg]\Bigg)^{-1}. \label{eq:lam1}
\end{align}
We now show that the terms inside the square brackets have limits $\overline x, \overline y \in \R$ as $n \to \infty$; i.e.,  
\begin{align}
\lim_{n\to \infty}    \sqrt{n(1-s^2) + n\Big( \frac{n\lambda - (q^{(B,i)} -\theta_{3})}{n\lambda + q^{(B,i)} -\theta_{3}} \Big)^2 s^2}  = \overline x  \quad \text{ and } \quad \lim_{n\to \infty} \sqrt{n}\frac{n\lambda - (q^{(B,i)} -\theta_{3})}{n\lambda + q^{(B,i)} -\theta_{3}} = \overline y. \label{eq:xnyn}
\end{align}
 Note that $\lim_{n\to \infty} s^2 = 1$, and recall the definition of $r$ to see that $\lim_{n\to \infty} r/n = 1+\lambda$,  so 
\begin{align*}
\lim_{n \to \infty} n(1-s^2) = \lim_{n \to \infty}  \frac{n(2r + 1)}{1 + 2r + r^2}  = \lim_{n \to \infty}   \frac{(2r + 1)/n}{(1 + 2r + r^2)/n^2}   = \lim_{n \to \infty}   \frac{2  }{ r/n }  = \frac{2 }{  1+\lambda }.
\end{align*}
Furthermore, recalling the definition of  $ q^{(B,i)}$, we have 
\begin{align}
 \lim_{n \to \infty}  n\Big( \frac{n\lambda - (q^{(B,i)} -\theta_{3})}{n\lambda + q^{(B,i)} -\theta_{3}} \Big)^2  =&\ \lim_{n \to \infty}  n\bigg( \frac{-\beta \sqrt{n} + q_2 + 1+ 2\lfloor \sqrt{n}\beta/2 \rfloor - \Big \lfloor \lfloor \sqrt{n}\beta/2 \rfloor \frac{i-1}{b+1} \Big \rfloor }{n\lambda + n- q_2 - 1- 2\lfloor \sqrt{n}\beta/2 \rfloor + \Big \lfloor \lfloor \sqrt{n}\beta/2 \rfloor \frac{i-1}{b+1} \Big \rfloor} \bigg)^2  \notag \\
=&\   \bigg( \frac{   \lim_{n \to \infty} q_2/ \sqrt{n}   -   \frac{i-1}{b+1} \beta/2 }{ \lambda + 1     } \bigg)^2. \label{eq:displayabove}
\end{align}
We know that $\lim_{n \to \infty} q_2/ \sqrt{n}$ exists because $q_2$ is fixed between zero and $2\lfloor \gamma \sqrt{n}\rfloor$. This proves \eqref{eq:xnyn}. Recall that $z = \lfloor \sqrt{n}\beta/2 \rfloor$ and $a = \lfloor \sqrt{n}\beta/2 \rfloor + \big \lfloor \lfloor \sqrt{n}\beta/2 \rfloor \frac{1}{b+1} \big \rfloor$. Since $r/n \to 1+\lambda$, it follows that 
\begin{align*}
&\lim_{n \to \infty} s^{a} = \lim_{n \to \infty} \big(1 - 1/(r+1)\big)^{\lfloor \sqrt{n}\beta/2 \rfloor + \big \lfloor \lfloor \sqrt{n}\beta/2 \rfloor / (b+1) \big \rfloor} = 1  \quad \text{ and } \quad \lim_{n \to \infty} s^{z} = \lim_{n \to \infty} s^{\lfloor \sqrt{n}\beta/2 \rfloor } = 1,
\end{align*}
and combined with \eqref{eq:lam1}, \eqref{eq:xnyn}, and Lemma~\ref{lem:exponential_aux}, this implies that   
\begin{align*}
\lim_{n\to\infty}\lambda^{z}_1(s) =&\ \lim_{n\to\infty}\lambda^{\lfloor \sqrt{n}\beta/2 \rfloor}_1(s) = \exp\Big(\frac{\overline x \beta}{2}  \Big)\exp\Big(-\frac{\overline y \beta}{2}  \Big)  ,\\
\lim_{n\to\infty}\lambda^{a}_1(s) =&\ \lim_{n\to\infty}\lambda^{\lfloor \sqrt{n}\beta/2 \rfloor + \big \lfloor \lfloor \sqrt{n}\beta/2 \rfloor \frac{1}{b+1} \big \rfloor}_1(s) =  \exp\Big(\frac{\overline x \beta}{2} \frac{b+2}{b+1} \Big)\exp\Big(-\frac{\overline y \beta}{2} \frac{b+2}{b+1} \Big).
\end{align*}
The expressions for $\lim_{n\to\infty}\lambda^{z}_2(s)$ and $\lim_{n\to\infty}\lambda^{a}_2(s)$ follow similarly. Comparing
\begin{align*}
\lambda_2(s) =&\ \Bigg(1-\sqrt{1-s^2 + \Big( \frac{n\lambda - (q^{(B,i)} -\theta_{3})}{n\lambda + q^{(B,i)} -\theta_{3}} \Big)^2 s^2}\Bigg) s^{-1}\bigg(1+\frac{n\lambda - (q^{(B,i)} -\theta_{3})}{n\lambda + q^{(B,i)} -\theta_{3}}\bigg)^{-1}
\end{align*}
to the form of $\lambda_1(s)$ in \eqref{eq:lam1}, we see that we can use \eqref{eq:xnyn} and  Lemma~\ref{lem:exponential_aux} again to  conclude that
\begin{align*}
\lim_{n\to\infty}\lambda^{z}_2(s) =&\  \exp\Big(-\frac{\overline x \beta}{2}  \Big)\exp\Big(-\frac{\overline y \beta}{2}  \Big), \text{ and } \\
\lim_{n\to\infty}\lambda^{a}_2(s) =&\ \exp\Big(-\frac{\overline x \beta}{2} \frac{b+2}{b+1} \Big)\exp\Big(-\frac{\overline y \beta}{2} \frac{b+2}{b+1} \Big).
\end{align*} 
For convenience, we define   $x = (\overline x - \overline y)\beta/2$ and $y = (\overline x + \overline y)\beta/2$, so that 
\begin{align*}
\lim_{n\to\infty}\lambda^{z}_1(s)= e^{x}, \quad \lim_{n\to\infty}\lambda^{a}_1(s) = e^{x(b+2)/(b+1)}, \quad \lim_{n\to\infty}\lambda^{z}_2(s)= e^{-y}, \quad \lim_{n\to\infty}\lambda^{a}_1(s) = e^{-y(b+2)/(b+1)}. 
\end{align*}
It is straightforward to check that  $x,y > 0$ using \eqref{eq:xnyn}. Let us now prove that  $\lim_{n\to \infty}  \E s^{D_z}  < 1$. Using the definition of $\E s^{D_z}$, we have 
\begin{align*}
\lim_{n\to \infty} \E s^{D_z} =&\  \lim_{n\to \infty}  \frac{\lambda_2^{z}(s)(\lambda_1^{a}(s)-1) - \lambda_1^{z}(s)(\lambda_2^{a}(s)-1)}{\lambda_1^{a}(s) - \lambda_2^{a}(s)} \\
=&\ \frac{e^{-y} ( e^{x(b+2)/(b+1)} - 1 ) -  e^{x} (e^{-y(b+2)/(b+1)} - 1 ) }{e^{x(b+2)/(b+1)} - e^{-y(b+2)/(b+1)}}.
\end{align*} 
Set $c = (b+2)/(b+1)$.  We want to show that for any $x,y > 0$, 
 \begin{align*}
 e^{-y} ( e^{xc} - 1 ) -  e^{x} (e^{-yc} - 1 ) < e^{xc} - e^{-yc}, \quad \text{ or } \quad e^{-y}   e^{xc}  -  e^{x}  e^{-yc}   < e^{xc} - e^{-yc} + e^{-y} - e^{x}.
 \end{align*} 
Rearranging  terms, this is equivalent to 
\begin{align*}
e^{ xc} (e^{-y} -1)    - e^{x}(e^{-yc}-1)  < - e^{-yc} + e^{-y}.
\end{align*}
Fix $ y > 0$ and treat the left-hand side as a function of $x$. Both sides are equal when $x = 0$, so it suffices to show  that the derivative of the left-hand side with respect to $x$  is negative. Now 
\begin{align}
&\frac{\partial}{\partial x} \Big(e^{ xc} (e^{- y} -1)    - e^{  x}(e^{-  y c}-1) \Big) = c e^{ x c} (e^{-   y} -1)    - e^{ x}(e^{-  y c}-1). \label{eq:musthave}
\end{align}
For the right-hand side to be negative, we must have 
\begin{align*}
c e^{ x (c-1)} > \frac{1-e^{-y c}}{1-e^{- y}}.
\end{align*}
Since $c = (b+2)/(b+1) > 1$, the left-hand side is bounded from below by $c$ provided that $x \geq 0$. The right-hand side converges to $c$ as $y \downarrow 0$,  so we must show that the derivative of the right-hand side is negative. Differentiating yields 
\begin{align*}
\frac{\partial }{\partial y} \frac{1-e^{-  y c}}{1-e^{- y}} =&\ \frac{ce^{-  y c}(1-e^{-  y}) -   e^{- y} (1-e^{-  y c})}{(1-e^{- y})^2} = e^{-  y} \times \frac{ce^{-  y(c-1)} - ce^{-  y c} - 1 +  e^{-  y c}}{(1-e^{-  y})^2}.
\end{align*}
The numerator $ce^{-y (c-1)} - (c-1)e^{- y(c-1)} - 1$ equals $0$ when $y= 0$. Its derivative equals
\begin{align*}
-c(c-1)e^{-  y(c-1)} +(c-1)^2 e^{-  y c} < -c(c-1)e^{-  y(c-1)} +(c-1)^2 e^{-  y(c-1)} \frac{c}{c-1} = 0, \quad y \geq 0,
\end{align*}
where the inequality is due to $e^{-y} \leq 1 < c/(c-1)$. Therefore, the numerator is strictly negative for $y > 0$, meaning that \eqref{eq:musthave} holds.  
\finishproof

\subsection{Proof of Lemma~\ref{lem:up}}
\label{proof:up} 
It suffices to show that $\E \tau^{-}(x^q_1) \leq C(\beta) \delta$, because
\begin{align*}
\Prob(V \leq \tau^{-}(x^q_1)) = \int_{0}^{\infty} \Prob(V\leq t) dF(t) = \int_{0}^{\infty} (1-e^{-t}) dF(t) = 1- \E e^{-\tau^{-}(x^q_1)}  
\leq  \E \tau^{-}(x^q_1),
\end{align*}
where $F(t)$ is the distribution function of $\tau^{-}(x^q_1)$. Define
\begin{align*}
\tau^{+}(q_1) = \inf_{t \geq 0} \{Q(t) = (q_1+1,0, \ldots, 0) |  Q(0)=(q_1,0, \ldots, 0) \}, \quad 0 \leq q_1 \leq n-1,
\end{align*}
and note that $\tau^{+}(q_1) = \tau^{-}(x^q_1)$.  If we let $\{\pi_{q}\}_{q \in S_Q}$ be the stationary distribution of the unscaled CTMC, it follows from (2.11) of \cite{BrowXia2001} that
\begin{align*}
\E \tau^{+}(q_1) =&\ \frac{\sum_{i=0}^{q_1}\pi_{i,0,\ldots, 0}}{n\lambda \pi_{q_1,0,\ldots, 0}}.
\end{align*}
Letting $f(x^q) = 1(x^q_1 \leq i)$ and using $\E G_{X} f(X) = 0$ yields   $n\lambda \pi_{i,0,\ldots, 0} = (i+1)\pi_{i+1,0,\ldots, 0}$,  which implies that $\pi_{i,0,\ldots, 0} =\pi_{0,\ldots, 0} (n\lambda)^{i}/i!$, so
\begin{align*}
\E \tau^{+}(q_1)  =  \frac{\sum_{k=0}^{q_1}\frac{(n\lambda)^{k}}{k!}}{n\lambda \frac{(n\lambda)^{q_1}}{q_1!}}=  \frac{q_1!}{(n\lambda)^{q_1+1}}\sum_{k=0}^{q_1}\frac{(n\lambda)^{k}}{k!}.
\end{align*}
 Note that $x^q_1 =\beta+  \delta(n\lambda -\lfloor n\lambda \rfloor)$ is equivalent to $q_1 = \lfloor n\lambda \rfloor$. If $\lfloor n\lambda \rfloor = 0$, we observe that the right-hand side  equals $1/(n\lambda)$, which verifies \eqref{eq:uphitbound} when $x^q_1 =\beta+  \delta(n\lambda -\lfloor n\lambda \rfloor)$. If, however, $\lfloor n\lambda \rfloor > 0$, we may use Stirling's approximation to see that for $q_1 > 0$,
\begin{align*}
\frac{q_1!}{(n\lambda)^{q_1+1}}\sum_{k=0}^{q_1}\frac{(n\lambda)^{k}}{k!} \leq \frac{3q_1^{q_1+1/2} e^{-q_1}}{(n\lambda)^{q_1+1}} \sum_{k=0}^{q_1}\frac{(n\lambda)^{k}}{k!}   \leq \frac{3q_1^{q_1+1/2} e^{-q_1}}{(n\lambda)^{q_1+1}} e^{n\lambda}.
\end{align*}
Setting $q_1 = \lfloor n\lambda \rfloor$ proves \eqref{eq:uphitbound} when $x_1 =\beta+  \delta(n\lambda -\lfloor n\lambda \rfloor)$. To prove \eqref{eq:uphitbound} when $x^q_1 = \delta $ and $x^q_1 = 2\delta$ requires just a little more work. Setting $q_1 = n-1$,
\begin{align*}
 \E \tau_{n-1}^{+}  \leq \frac{3(n-1)^{n-1/2} e^{-(n-1)}}{(n\lambda)^{n}} e^{n\lambda} \leq&\ \frac{3e}{\sqrt{n-1}}   \frac{n^n}{(n-\beta \sqrt{n})^n} e^{-n} e^{n\lambda}  = \frac{3e}{\sqrt{n-1}} \Big(1- \frac{\beta}{\sqrt{n}}\Big)^{-n} e^{-\beta\sqrt{n}}.
\end{align*}
To conclude, we need to bound
\begin{align*}
 \Big(\Big(1- \frac{\beta}{\sqrt{n}}\Big)^{-\sqrt{n}} e^{-\beta }\Big)^{\sqrt{n}} 
=&\ \bigg(\exp \Big(-\sqrt{n}\log\Big(1- \frac{\beta}{\sqrt{n}}\Big) - \beta \Big)\bigg)^{\sqrt{n}}.
\end{align*}
Using Taylor expansion,
\begin{align*}
\log\Big(1- \frac{\beta}{\sqrt{n}}\Big) =  - \frac{\beta}{\sqrt{n}} - \frac{1}{2} \Big(\frac{\beta}{\sqrt{n}}\Big)^2 \frac{1}{(1+\xi(\beta/\sqrt{n}))^2}
\end{align*}
where $\xi(\beta/\sqrt{n}) \in [-\beta/\sqrt{n},0]$. Therefore,
\begin{align*}
\bigg(\exp \Big(-\sqrt{n}\log\Big(1- \frac{\beta}{\sqrt{n}}\Big) - \beta \Big)\bigg)^{\sqrt{n}} =&\ \exp \bigg( \frac{\beta^2/2}{(1+\xi(\beta/\sqrt{n}))^2}\bigg),
\end{align*}
and we conclude that 
\begin{align*}
\sup_{n \geq 0} \Big(\Big(1- \frac{\beta}{\sqrt{n}}\Big)^{-\sqrt{n}} e^{-\beta }\Big)^{\sqrt{n}} < \infty.
\end{align*}
The argument when $q_1 = n-2$ is  identical. This proves \eqref{eq:uphitbound} when $x_1 = \delta, 2\delta$. 
\finishproof

\end{APPENDICES}

% References here (outcomment the appropriate case)

% CASE 1: BiBTeX used to constantly update the references
%   (while the paper is being written).
\bibliographystyle{informs2014} % outcomment this and next line in Case 1
\bibliography{dai20190911} % if more than one, comma separated

%%%%%%%%%%%%%%%%%
\end{document}